\documentclass[11pt,twoside]{amsart}
\date{14 August 2010}
\usepackage{latexsym,amsthm,amsmath,amsfonts,amscd,amssymb}
\usepackage[all]{xy}
\usepackage[utf8]{inputenc}
\usepackage{graphics}
\usepackage{lscape}
\usepackage{array}
\theoremstyle{plain}  % default
\newtheorem{theorem}{Theorem}[section]

\newtheorem*{theorem*}{Theorem}

\newtheorem{corollary}[theorem]{Corollary}
\newtheorem{lemma}[theorem]{Lemma}
\newtheorem{proposition}[theorem]{Proposition}

\newtheorem{tech-lemma}[theorem]{Technical Lemma}
\newtheorem{definition}[theorem]{Definition}
\theoremstyle{remark}

\newtheorem{notation}[theorem]{Notation}
\newtheorem{remark}[theorem]{Remark}
\newtheorem*{remark*}{Remark}

\newtheorem*{claim*}{Claim}
\numberwithin{equation}{section}
\renewcommand{\leq}{\leqslant}

\renewcommand{\geq}{\geqslant}

\newcommand{\R}{\mathbb{R}}
\newcommand{\Z}{\mathbb{Z}}
\newcommand{\C}{\mathbb{C}}

\newcommand{\M}{{\mathcal M}}
\newcommand{\N}{{\mathcal N}}
\newcommand{\cR}{{\mathcal R}}

\newcommand{\vol}{\mathrm{vol}}

\newcommand{\PGL}{\mathrm{PGL}}
\newcommand{\PSL}{\mathrm{PSL}}
\newcommand{\PO}{\mathrm{PO}}
\newcommand{\PSO}{\mathrm{PSO}}

\newcommand{\U}{\mathrm{U}}
\newcommand{\GL}{\mathrm{GL}}

\newcommand{\SL}{\mathrm{SL}}
\newcommand{\SO}{\mathrm{SO}}
\newcommand{\Or}{\mathrm{O}}

\newcommand{\Spin}{\mathrm{Spin}}
\newcommand{\Pin}{\mathrm{Pin}}
\newcommand{\EGL}{\mathrm{EGL}}
\newcommand{\EO}{\mathrm{EO}}
\newcommand{\ESO}{\mathrm{ESO}}

\newcommand{\st}{\;|\;}

\DeclareMathOperator{\Jac}{Jac}

\DeclareMathOperator{\Ad}{Ad}
\DeclareMathOperator{\Aut}{Aut}

\DeclareMathOperator{\rk}{rk}

\DeclareMathOperator{\Hom}{Hom}

\DeclareMathOperator{\ima}{Im}
\newcommand{\liem}{\mathfrak{m}}

\newcommand{\liemc}{\mathfrak{m}^{\mathbb{C}}}
\newcommand{\lieh}{\mathfrak{h}}
\newcommand{\liehc}{\mathfrak{h}^{\mathbb{C}}}
\newcommand{\lieg}{\mathfrak{g}}
\newcommand{\liegc}{\mathfrak{g}^{\mathbb{C}}}
\let\oldmarginpar\marginpar
\renewcommand\marginpar[1]{\oldmarginpar{\tiny\bf\begin{flushleft} #1
\end{flushleft}}}

\begin{document}

%%%%%%%%%%%%%%%%%%%%%%%%%%%%%%%%%%%%%%%%%%%%%%%%%%%%%%%%%%%%%%%%
%
% Title.
%
%%%%%%%%%%%%%%%%%%%%%%%%%%%%%%%%%%%%%%%%%%%%%%%%%%%%%%%%%%%%%%%%

\title[Representations of Surface Groups]
{Representations of Surface Groups in the Projective General Linear Group}
%%%%%%%%%%%%%%%%%%%%%%%%%%%%%%%%%%%%%%%%%%%%%%%%%%%%%%%%%%%%%%%%
%
% Author(s), affiliation(s) and email(s).
%
%%%%%%%%%%%%%%%%%%%%%%%%%%%%%%%%%%%%%%%%%%%%%%%%%%%%%%%%%%%%%%%%

\author[André Gama Oliveira]{André Gama Oliveira}
\address{Centro de Matemática da
  Universidade de Trás-os-Montes e Alto Douro (CMUTAD)\\
  Quinta dos Prados, Apartado 1013 \\ 5000-911 Vila Real \\ Portugal }
\email{agoliv@utad.pt}

\thanks{This work was partially supported by Centro de Matemática da
Universidade do Porto and by the grant SFRH/BD/23334/2005 and the
project POCTI/MAT/58549/2004, financed by Fundação para a Ciência e
a Tecnologia (Portugal) through the programmes POCTI and POSI of the
QCA III (20002006) with European Community (FEDER) and national
funds.}


\begin{abstract}
Given a closed, oriented
surface $X$ of genus $g\geq 2$, and a semisimple Lie group $G$, let $\cR_G$ be the moduli space of
reductive representations of $\pi_1X$ in $G$.
We determine the number of connected components of $\cR_{\PGL(n,\R)}$, for $n\geq 4$ even.
In order to have a first division of connected components, we first classify real projective bundles over such a
surface. Then we achieve our goal, using holomorphic methods through the theory of Higgs bundles over compact Riemann surfaces.

We also show that the complement of the Hitchin component in $\cR_{\SL(3,\R)}$ is homotopically equivalent to $\cR_{\SO(3)}$.
\end{abstract}

\maketitle

%%%%%%%%%%%%%%%%%%%%%%%%%%%%%%%%%%%%%%%%%%%%%%%%%%%%%%%%%%%%%%%%
%
% Text.
%
%%%%%%%%%%%%%%%%%%%%%%%%%%%%%%%%%%%%%%%%%%%%%%%%%%%%%%%%%%%%%%%%

%%%%%%%%%%%%%%%%%%%%%%%%%%%%%%%%%%%%%%%%%%%%%%%%%%%%%%%%%%%%%%%%%%%%%%
\section{Introduction}
%%%%%%%%%%%%%%%%%%%%%%%%%%%%%%%%%%%%%%%%%%%%%%%%%%%%%%%%%%%%%%%%%%%%%%

Let $X$ be a closed, oriented surface of genus $g\geq 2$ and
$\pi_1X$ be its fundamental group. Let
$$\mathcal{R}_{\PGL(n,\R)}=\Hom^\text{red}(\pi_1X,\PGL(n,\R))/\PGL(n,\R)$$
be the quotient space of reductive representations of
$\pi_1X$ in the projective general linear group
$\PGL(n,\R)=\GL(n,\R)/\R^*$, where $\PGL(n,\R)$ acts by
conjugation. In this paper we determine the number of connected
components of $\mathcal{R}_{\PGL(n,\R)}$, for $n\geq 4$ even, applying
the general theory of $G$-Higgs bundles to the $\PGL(n,\R)$ case.

For a semisimple Lie group $G$, this general theory, created among
others by Hitchin \cite{hitchin:1987}, Simpson \cite{simpson:1992,simpson:1994a,simpson:1994b},
Corlette \cite{corlette:1988} and Donaldson \cite{donaldson:1987}, supplies a strong
relation between different subjects such as topology, holomorphic
and differential geometry and analysis. On the one hand, we have the
moduli space $\mathcal{R}_G$ of reductive representations of
$\pi_1X$ in $G$, also known as a character variety.
An element in $\mathcal{R}_G$ is topologically classified by certain
invariants of the isomorphism class of the associated flat principal
$G$-bundle over $X$. If $c$ is a topological class of
principal $G$-bundles, we denote by $\mathcal{R}_G(c)$ the subspace
of $\mathcal{R}_G$ consisting of classes of representations which
belong to the class $c$. On the other hand we fix a complex
structure on $X$ turning it into a Riemann surface, and consider $G$-Higgs bundles over it. A $G$-Higgs bundle is a pair
consisting of a holomorphic bundle, whose structure group depends on $G$, and a section of a certain
associated bundle (see below for precise definitions).
Topologically, a $G$-Higgs bundle is also classified by invariants
taking values in the same set as the representations in
$\mathcal{R}_G$. Again, if $c$ is one topological class, we denote
by $\M_G(c)$ the moduli space of polystable $G$-Higgs bundles in the class
$c$.

Now, the above mentioned authors have proved that the spaces $\mathcal{R}_G(c)$ and $\M_G(c)$
are homeomorphic (see Theorem \ref{fundamental correspondence for semisimple G}). More generally, for a reductive Lie group $G$, there is
a correspondence similar to the previous one, but replacing
$\pi_1X$ by its universal central extension $\Gamma$, defined
in (\ref{extension}) below. We denote the space of such representations, with fixed topological class $c$, by $\mathcal{R}_{\Gamma,G}(c)$.
Related to these two moduli spaces, and essential in the proof of
the existence of the homeomorphism, is a third moduli space:
the moduli space of solutions to the so-called Hitchin's equations on
a fixed $C^\infty$ principal $G$-bundle over $X$.

For $G$ compact and connected, the spaces $\mathcal{R}_G$ and $\M_G$
have been studied in the seminal papers of Narasimhan and Seshadri
\cite{narasimhan-seshadri:1965} and of Ramanathan \cite{ramanathan:1975} from an algebraic viewpoint,
and by Atiyah and Bott \cite{atiyah-bott:1982} from a gauge theoretic point of
view. In this case, the answer about the number of components is known: for each
topological type $c$, each subspace of $\mathcal{R}_G(c)$ is
connected. Since then much has been done to study the geometry and
topology of these spaces. When $G$ is complex, connected and reductive,
the answer to the problem of counting connected components is the
same as in the compact case by the works of Hitchin \cite{hitchin:1987},
Donaldson \cite{donaldson:1987}, Corlette \cite{corlette:1988} and Simpson
\cite{simpson:1992,simpson:1994a,simpson:1994b}. When $G$ is a non-compact real form
of a complex semisimple Lie group, the study of the topology of
$\mathcal{R}_G$ started with the seminal papers of Goldman
\cite{goldman:1988} and Hitchin \cite{hitchin:1992} and, although much work
has been done since then by several people (see, in particular, the paper \cite{gothen:2001} of Gothen, the works \cite{bradlow-garcia-prada-gothen:2003,bradlow-garcia-prada-gothen:2004,bradlow-garcia-prada-gothen:2005} of Bradlow, García-Prada and Gothen and also \cite{garcia-prada-gothen-mundet:2008} by García-Prada, Gothen and Mundet i Riera), it is still far from finished.

In this paper we are interested in studying the components of
$\mathcal{R}_{\PGL(n,\R)}$, for $n\geq 4$ even (when $n\geq 3$ is odd, $\PGL(n,\R)\cong\SL(n,\R)$ hence the components of $\mathcal{R}_{\PGL(n,\R)}$ are known to be $3$ by the work of Hitchin in \cite{hitchin:1992}; the $n=2$ case was studied by Xia in \cite{xia:1999}). Following the ideas of
Hitchin \cite{hitchin:1987,hitchin:1992}, the main tool to reach our goal should be
the $L^2$-norm of the Higgs field in $\M_{\PGL(n,\R)}(c)$, but
in our case another group naturally appears. We will work with the space
$\M_{\EGL(n,\R)}$ of $\EGL(n,\R)$-Higgs bundles ($\EGL(n,\R)=\GL(n,\R)\times\U(1)/\hspace{-0,1cm}\sim$, where $(A,z)\sim (-A,-z)$). This is done mainly for two related reasons. One is that with
this new group we can work with holomorphic vector bundles, rather
than just principal or projective bundles. The other is that we can realize space of $\EGL(n,\R)$-Higgs bundles as
closed subspace of $\M_{\GL(n,\C)}\times\Jac^d(X)$, where
$\M_{\GL(n,\C)}$ is the moduli space of Higgs bundles (see
\cite{hitchin:1987}). In general, when $\M_G(c)$ is smooth, the function $f$ given by
the $L^2$-norm of the Higgs field is
a non-degenerate Morse-Bott function which is also a proper map and, in
some cases, the critical submanifolds are well enough understood to
allow the extraction of topological information such as the Poincaré
polynomial. However, even when $\M_G(c)$ has singularities, the properness
of $f$ allows us to draw conclusions about the connected components,
although one cannot directly apply Morse theory.

The study of the local minima of $f$ is sufficient to obtain the
number of connected components of the space of $\EGL(n,\R)$-Higgs bundles and thus of
$\mathcal{R}_{\Gamma,\EGL(n,\R)}$, for $n\geq
4$. There is a projection from this space to
$\mathcal{R}_{\PGL(n,\R)}$ and using this we compute the components
of $\mathcal{R}_{\PGL(n,\R)}$, obtaining the first of our two main results (see Theorem \ref{princ}):
\begin{theorem}\label{main1}
Let $n\geq 4$ be even. Then the space $\mathcal{R}_{\PGL(n,\R)}$ has
$2^{2g+1}+2$ connected components.
\end{theorem}

Essential in the count of components of $\mathcal{R}_{\PGL(n,\R)}$ is the topological classification of real projective bundles over $X$. This is done in the first part of the paper, where we have found explicit discrete invariants which classify continuous principal $\PGL(n,\R)$-principal bundles over any closed oriented surface. This classification shows for example that, in contrast to the complex case, there are real projective bundles which are not projectivization of real vector bundles. It shows also that, in most cases, there is a collapse of the second Stiefel-Whitney class, when we pass from real vector bundles to projective bundles.

Combining the results of Xia \cite{xia:1999} for $n=2$ and of Hitchin \cite{hitchin:1992} for $n\geq 3$ odd with our Theorem \ref{main1}, we have the number of connected components of $\mathcal{R}_{\PGL(n,\R)}$, for arbitrary $n$, as follows:
\begin{theorem}The number of connected components of $\mathcal{R}_{\PGL(n,\R)}$ is:
\begin{itemize}
 \item $2^{2g+1}+4g-5$ if $n=2$;
 \item $3$ \hspace{2,05cm}if $n\geq 3$ is odd;
 \item $2^{2g+1}+2$ \hspace{0,77cm}if $n\geq 4$ is even.
\end{itemize}
\end{theorem}

Using the results of Hitchin in \cite{hitchin:1992}, we are able to obtain
more topological information of $\mathcal{R}_{\PGL(3,\R)}$ (observe that this is the same as $\cR_{\SL(3,\R)}$ since $\PGL(3,\R)\cong\SL(3,\R)$), because
in this case there are no critical submanifolds of $f$ besides the
local minima and these are of a very special type. The result we obtain is the following (see Theorem \ref{princ2}):
\begin{theorem}
The space $\cR_{\SL(3,\R)}$ has
one contractible component and the space consisting of the other two
components is homotopically equivalent to $\mathcal{R}_{\SO(3)}$.
\end{theorem}

Actually, using a computation of the Poincaré polynomials of $\mathcal{R}_{\SO(3)}$ recently done by Ho and Liu in \cite{ho-liu:2007}, this theorem gives the Poincaré polynomials of $\mathcal{R}_{\SL(3,\R)}$ almost for free.

%%%%%%%%%%%%%%%%%%%%%%%%%%%%%%%%
\section{Representations of $\pi_1X$ in $G$ and $G$-Higgs bundles}\label{Representation spaces}
%%%%%%%%%%%%%%%%%%%%%%%%%%%%%%%%
\subsection{Representations of $\pi_1X$ in $G$}
Let $X$ be a closed oriented surface of genus $g\geq 2$ and let $G$ be a semisimple Lie group.

Consider the space $\Hom(\pi_1X,G)$ of all homomorphisms from the fundamental group of $X$ to $G$. Such a homomorphism $\rho:\pi_1X\to G$ is also called a \emph{representation} of $\pi_1X$ in $G$.
Considering the presentation of $\pi_1X$ given by the usual $2g$ generators 
\begin{equation}\label{pi1X}
\pi_1X=\big\langle a_1,b_1,\dots,a_g,b_g\st\prod_{i=1}^g[a_i,b_i]=1\big\rangle
\end{equation} 
one sees that a representation $\rho\in\Hom(\pi_1X,G)$ is
determined by its values on the set of generators
$a_1,b_1,\dots,a_g,b_g$. The set
$\Hom(\pi_1X,G)$ can thus be embedded in $G^{2g}$ via
$\rho\mapsto(\rho(a_1),\dots,\rho(b_g))$, becoming the subset of $2g$-tuples
$(A_1,B_1,\dots,A_g,B_g)$ of $G^{2g}$ satisfying the algebraic
equation $\prod_{i=1}^{g}[A_i,B_i]=1$, and we consider the induced topology on
$\Hom(\pi_1X,G)$.

 Letting $G$ act on $\Hom(\pi_1X,G)$ by conjugation $$g\cdot\rho=g\rho g^{-1}$$ we obtain the quotient space
$$\Hom(\pi_1X,G)/G.$$ This space may not be Hausdorff because there may exist different orbits with non-disjoint closures, so we consider only \emph{reductive} representations of $\pi_1X$ in $G$, meaning the ones that, when composed with the adjoint representation of $G$ on its Lie algebra, become a sum of irreducible representations. Denote the space of such representations by $\Hom^\text{red}(\pi_1X,G)$. The corresponding quotient is the space we are interested in:

\begin{definition}\label{moduli rep of pi1X}
The \emph{moduli space of representations of $\pi_1X$ in
$G$} is the quotient space $$\cR_G=\Hom^\mathrm{red}(\pi_1X,G)/G.$$ 
\end{definition}

The space $\cR_G$ is also known as the $G$-\emph{character variety} of $X$.

If $G$ acts on $G^{2g}$ through the diagonal adjoint action, the inclusion $j:\Hom(\pi_1X,G)\hookrightarrow G^{2g}$ becomes $G$-equivariant and, from Theorem $11.4$ in \cite{richardson:1988}, a representation $\rho\in\Hom(\pi_1X,G)$ is reductive if and only if the orbit of $j(\rho)$ in $G^{2g}$ is closed, hence it follows that $\cR_G$ is indeed Hausdorff.

\vspace{0,5cm}
If we allow $G$ to be reductive and not just semisimple, then we consider a universal central extension
$\Gamma$ of $\pi_1X$ given by the short exact sequence
$$0\longrightarrow\Z\longrightarrow\Gamma\longrightarrow\pi_1X\longrightarrow0.$$ It is generated by $2g$ generators $a_1,\ldots,b_g$ (which are mapped to the corresponding ones of $\pi_1X$) and by a central element $J$, subject to
the relation $\prod_{i=1}^g[a_i,b_i]=J$:
\begin{equation}\label{extension}
\Gamma=\big\langle
a_1,b_1,\dots,a_g,b_g,J\,|\,\prod_{i=1}^g[a_i,b_i]=J,\ J\in Z(\Gamma)\big\rangle.
\end{equation}

Let $H\subseteq G$ be a maximal compact subgroup of $G$. Analogously to the case of $\pi_1X$, let us consider the reductive representations $\rho$
of $\Gamma$ in $G$ such that $\rho(J)\in(Z(G)\cap H)_0$, the identity component of the centre of $G$ intersected with $H$:
\begin{equation}\label{homgammaG}
\Hom^\text{red}_{\rho(J)\in(Z(G)\cap H)_0}(\Gamma,G)=\{\rho:\Gamma\longrightarrow G\st\rho\text{ is reductive and }\rho(J)\in(Z(G)\cap H)_0\}.
\end{equation}

This definition does not depend on the choice of $H$. Indeed, any other maximal compact subgroup $H'$ of $G$ is conjugate to $H$, so $(Z(G)\cap H)_0=(Z(G)\cap H')_0$.

\begin{definition}\label{Gamma}
The \emph{moduli space of representations of $\Gamma$ in $G$} is the
quotient space
$$\cR_{\Gamma,G}=\Hom^\mathrm{red}_{\rho(J)\in(Z(G)\cap H)_0}(\Gamma,G)/G.$$
\end{definition}

To give a representation $\rho\in\Hom^\text{red}_{\rho(J)\in(Z(G)\cap H)_0}(\Gamma,G)$ is equivalent to give a representation
of $\pi_1(X\setminus\{x_0\})$, the fundamental group of the punctured surface, in $G$ such that the
image of the homotopy class of the loop around the puncture is $\rho(J)\in(Z(G)\cap H)_0$.

Of course, if $G$ is semisimple, $\cR_{\Gamma,G}=\cR_G$.
The main result of this paper is the computation of the number of connected components of $\cR_{\PGL(n,\R)}$, for $n\geq 4$ even.

\vspace{0,5cm}
As is well-known, there is a bijection between isomorphism classes of representations of $\pi_1X$ in $G$ and isomorphism classes of flat $G$-bundles over $X$.
There is as well a one-to-one correspondence between isomorphism classes of representations of $\Gamma$ in $G$ and isomorphism classes of projectively flat $G$-bundles over $X$, i.e., $G$-bundles equipped with connections with constant central curvature in $Z(\lieg)=\mathrm{Lie}(Z(G)_0)$.
Taking these correspondences into account, we make the following definition:
\begin{definition}\label{def top invariant representation}
Let $\rho$ be a representation of $\pi_1X$ in $G$. A \emph{topological invariant} of $\rho$ is a topological invariant of the associated flat $G$-bundle.

Let $\rho$ be a representation of $\Gamma$ in $G$. A \emph{topological invariant} of $\rho$ is a topological invariant of the associated projectively flat $G$-bundle.
\end{definition}

If two representations
$\rho_1,\rho_2\in\Hom^\text{red}(\pi_1X,G)$ are equivalent,
then the associated principal flat $G$-bundles $E_{\rho_1}$ and
$E_{\rho_2}$ are isomorphic and vice-versa. Hence the topological
invariants of $\rho_1$ and of $\rho_2$ are the same. Thus it makes
sense to define a topological invariant of an equivalence class of
representations.
Given a topological class $c$ of $G$-bundles over $X$, denote by $$\cR_G(c)$$ the subspace of $\cR_G$ whose representations belong to the class $c$. Analogously, define $\cR_{\Gamma,G}(c)$.

%%%%%%%%%%%%%%%%%%%%%%%%%%%%%%%
\subsection{$G$-Higgs bundles}\label{G-Higgs bundles}
%%%%%%%%%%%%%%%%%%%%%%%%%%%%%%%

In this section we introduce the main objects which we shall work with. These are called \emph{Higgs bundles} and roughly are pairs consisting of a holomorphic bundle and a section of an associated bundle (see Definition \ref{definition of Higgs bundle} below). Higgs bundles were introduced by Hitchin \cite{hitchin:1987} on compact Riemann surfaces and by Simpson \cite{simpson:1992} on any compact K\"ahler manifold. 

Let $H\subseteq G$ be a maximal compact subgroup of $G$
and $H^\C\subseteq G^\C$ their complexifications. There is a Cartan
decomposition of $\lieg$,
$$\lieg=\lieh\oplus\liem$$
where $\liem$ is the complement of $\lieh$ with
respect to the non-degenerate $\Ad(G)$-invariant bilinear $B$ form on $\lieg$. If
$\theta:\lieg\to\lieg$ is the
corresponding Cartan involution then $\lieh$ and
$\liem$ are its $+1$-eigenspace and $-1$-eigenspace,
respectively. Complexifying, we have the decomposition
$$\liegc=\liehc\oplus\liemc$$ and
$\liemc$ is a representation of $H^\C$ through the
so-called \emph{isotropy representation}
\begin{equation}\label{complex isotropy representation}
 \Ad|_{H^\C}:H^\C\longrightarrow\Aut(\liemc)
\end{equation}
which is induced by the adjoint
representation of $G^\C$ on $\liegc$.
If $E_{H^\C}$ is a principal $H^{\C}$-bundle over $X$, we denote by
$E_{H^\C}(\liemc)=E\times_{H^{\C}}\liemc$ the
vector bundle, with fibre $\liemc$, associated to the
isotropy representation.

Let $K=T^*X^{1,0}$ be the canonical line bundle of $X$. 
\begin{definition}\label{definition of Higgs bundle}
A \emph{$G$-Higgs bundle} over a Riemann surface $X$ is a pair
$(E_{H^\C},\Phi)$ where $E_{H^\C}$ is a principal holomorphic $H^\C$-bundle over
$X$ and $\Phi$ is a global holomorphic section of
$E_{H^\C}(\liemc)\otimes K$, called the \emph{Higgs field}.
\end{definition}

Any continuous $G$-bundle has certain discrete invariants which distinguish bundles which are not isomorphic as continuous (or equivalently $C^\infty$) $G$-bundles. On Riemann surfaces, if $G$ is connected, these invariants take values in $\pi_1G$. For example, complex vector bundles of rank $n$ are classified by their degree $d\in\Z=\pi_1\U(n)$.
For $G$ not necessarily connected, these topological invariants may take values in more complicated sets which depend only on the homotopy type of $G$. If $H$ is a maximal compact subgroup of $G$, then the inclusion $H\subset G$ is a homotopy equivalence so the classification of $G$-bundles is equivalent to that of $H$-bundles.
Now, a $G$-Higgs bundle $(E_{H^\C},\Phi)$ is topologically classified by the topological invariant of the corresponding $H^\C$-bundle $E_{H^\C}$ and, as the maximal compact subgroup of $H^\C$ is $H$, the topological classification of $G$-Higgs bundles is the same as the one of $H$-principal bundles.

\vspace{0,5cm}
Now we consider the moduli space of $G$-Higgs bundles. The notion of (poly)stability, for general $G$ is subtle (see \cite{bradlow-garcia-prada-mundet:2003,garcia-prada-gothen-mundet:2008,schmitt:2004,schmitt:2008}) but for $\GL(n,\C)$ it is easy.
Consider a ($\GL(n,\C)$-)Higgs bundle $(V,\Phi)$ and let $$\mu(V)=\frac{\deg(V)}{\rk(V)}$$ be the slope of $V$. A subbundle $W\subseteq V$ is said \emph{$\Phi$-invariant} if $\Phi(W)\subset W\otimes K$.
\begin{definition}\label{(poly)stability for usual Higgs bundles}
A Higgs bundle $(V,\Phi)$ is:
\begin{itemize}
 \item \emph{stable} if $\mu(W)<\mu(V)$ for all $\Phi$-invariant proper subbundle $W\subset V$;
 \item \emph{semistable} if $\mu(W)\leq\mu(V)$ for all $\Phi$-invariant proper subbundle $W\subset V$;
 \item \emph{polystable} if $V=W_1\oplus\dots\oplus W_k$ and $\Phi=\Phi_1\oplus\dots\oplus\Phi_k$ where, for each $i$, $\Phi_i:W_i\to W_i\otimes K$ and $(W_i,\Phi_i)$ is stable with $\mu(W_i)=\mu(V)$.
\end{itemize}
\end{definition}

\begin{definition}\label{isomorphism of G-Higgs bundles}
 Two $G$-Higgs bundles $(E_{H^\C},\Phi)$ and $(E_{H^\C}',\Phi')$ over $X$ are \emph{isomorphic} if there is an holomorphic isomorphism $f:E_{H^\C}\to E_{H^\C}'$ such that $\Phi'=\tilde f(\Phi)$, where $\tilde f\otimes 1_K:E_{H^\C}(\liemc)\otimes K\to E_{H^\C}'(\liemc)\otimes K$ is the map induced from $f$ and from the isotropy representation $H^\C\to\Aut(\liemc)$.
\end{definition}

In order to construct moduli spaces, we need to consider $S$-equivalence classes of semistable $G$-Higgs bundles (cf. \cite{schmitt:2008}). For a stable $G$-Higgs bundle, its $S$-equivalence class coincides with its isomorphism class and for a strictly semistable $G$-Higgs bundle, its $S$-equivalence contains precisely one (up to isomorphism) representative which is polystable so this class can be thought as the isomorphism class of the unique polystable $G$-Higgs bundle which is $S$-equivalent to the given strictly semistable one.

These moduli spaces have been constructed by Schmitt in \cite{schmitt:2004,schmitt:2005,schmitt:2008}, using methods of Geometric Invariant Theory, showing that they carry a natural structure of algebraic/complex variety.

\begin{definition}\label{moduli of Higgs bundles}
For a reductive Lie group $G$, the \emph{moduli space of $G$-Higgs bundles over a Riemann surface $X$} is the algebraic/complex variety of isomorphism classes of polystable $G$-Higgs bundles. We denote it by $\M_G$:
$$\M_G=\{\text{Polystable }G\text{-Higgs bundles on }X\}/\sim.$$ For a fixed topological class $c$ of $G$-Higgs bundles, denote by $\M_G(c)$ the moduli space of $G$-Higgs bundles which belong to the class $c$.
\end{definition}

The relation between $G$-Higgs bundles over $X$ and representations $\pi_1X\to G$ is given by the following fundamental theorem. 
\begin{theorem}\label{fundamental correspondence for semisimple G}
Let $G$ be a semisimple Lie group. A $G$-Higgs bundle is polystable if and only if it arises from a reductive representation of $\pi_1X$ in $G$. Moreover, this correspondence induces a homeomorphism between the spaces $\cR_G(c)$ and
$\M_G(c)$.

If $G$ is reductive, there is a similar correspondence which induces a homeomorphism between the spaces $\cR_{\Gamma,G}(c)$ and $\M_G(c)$.
\end{theorem}

Strictly speaking, this theorem has been proved for $G=\GL(n,\C)$ and $G=\SL(n,\C)$ by Hitchin in \cite{hitchin:1987} and Simpson in \cite{simpson:1992} (see also the papers \cite{corlette:1988} of Corlette and \cite{donaldson:1987} of Donaldson). The general definition of polystability and the proof of the Hitchin-Kobayashi correspondence for arbitrary $G$-Higgs bundles appears in the preprint \cite{garcia-prada-gothen-mundet:2008} of García-Prada, Gothen and Mundet i Riera.
%%%%%%%%%%%%%%%%%%%%%%%%%%%%%%%%%%%%%%%%%%%%%%%%%%%%%%%%%%%%%%%
\section{Topological invariants for $\PGL(n,\R)$-bundles over closed oriented surfaces}\label{Top}
%%%%%%%%%%%%%%%%%%%%%%%%%%%%%%%%%%%%%%%%%%%%%%%%%%%%%%%%%%%%%%%

In this section we obtain a topological classification of continuous principal $\PGL(n,\R)$-bundles over $X$, with $n\geq 4$ even. We shall, however, start by obtaining a general topological classification for any principal $G$-bundles, with $\pi_0G$ abelian.

\subsection{The case of any topological group $G$ with $\pi_0G$ abelian}

Let $G$ be a topological group.
Denote by $\mathcal{C}(G)$ the sheaf of continuous $G$-valued
functions on $X$ and by $G_0$ the identity component of $G$. We have the short exact sequence of groups
$$\xymatrix{0\ar[r]&G_0\ar[r]&G\ar[r]^(.35){p_1}&\pi_0G\ar[r]&0}$$
and, associated to the corresponding short exact sequence of sheaves
of continuous functions with values in the corresponding groups, we
have the sequence of cohomology sets:
$$\xymatrix{H^1(X,\mathcal{C}(G_0))\ar[r]&H^1(X,\mathcal{C}(G))\ar[r]^{p_{1,*}}&H^1(X,\pi_0G)}.$$
Recall that the cohomology set $H^1(X,\mathcal{C}(G))$ is in natural bijection with the set of isomorphism classes of continuous $G$-principal bundles over $X$. So, from the previous sequence, we define the first topological invariant of a continuous 
$G$-bundle $E$.

\begin{definition}\label{mu1}
The topological invariant $\mu_1$ of $E$ is defined by $$\mu_1(E)=p_{1,*}(E)\in
H^1(X,\pi_0G).$$
\end{definition}

Of course, this invariant yields the obstruction to reducing the structure group of $E$ to $G_0$.
Notice that, if $\pi_0G$ is abelian, $H^1(X,\pi_0G)\cong\Hom(\pi_1X,\pi_0G)\cong\pi_0G^{2g}$. 
From now on we assume that we are on this case: $\pi_0G$ is an abelian group.

% Choose the usual $2$-dimensional CW-complex structure on $X$, with a single $0$-cell, $2g$ $1$-cells and a single $2$-cell. Let $X_1$ denote the corresponding $1$-skeleton of $X$. Relative to a
% covering of $X_1$ by connected open sets, $\mu_1(E)$ tells us in which
% component of $G$ the transition functions of $E|_{ X_1}$ take their
% values, over each one of the $2g$ loops of $X_1$. So, a
% continuous $G$-bundle over $X_1$ is completely classified by
% $\mu_1\in\pi_0G^{2g}$, being trivial over
% $X_1$ if and only if $\mu_1=0$.

Our initial classification of $G$-bundles with $\mu_1$ fixed was much more complicated and was splitted into two assymetric parts: $\mu_1=0$ and $\mu_1\neq 0$. I am gretly indebted to an anonymous referee for providing a much simpler argument for the case $\mu_1\neq 0$ and which allows to study both cases $\mu_1=0$ and $\mu_1\neq 0$ simultaneously. The argument is as follows.

The surface $X$ is homeomorphic to the result of identifying (using orientation reversing homeomorphisms) the sides of a regular $4g$-gon $P$ according to the rule $$A_1B_1A_1^{-1}B_1^{-1}A_2B_2A_2^{-1}B_2^{-1}\cdots A_gB_gA_g^{-1}B_g^{-1}.$$ Let $\pi:P\to X$ be the natural projection. Let $c$ be the centre of $P$, $B(c,\epsilon)\subset P$ be a small disc centred at $c$ of radius $\epsilon$, disjoint from the boundary of $P$, and let $$U=\pi(P\setminus\{c\})\ \text{ and }\ V=\pi(B(c,\epsilon)).$$

A $G$-principal bundle $E$ on $X$ can be described by its restrictions $$E_U=E|_U\ \text{ and }\ E_V=E|_V$$ and by the gluing data $$\rho:E_U|_{U\cap V}\stackrel{\cong}{\longrightarrow} E_V|_{U\cap V}.$$ 

As $V$ is contractible, $E_V$ is isomorphic to the trivial bundle. On the other hand, the fact that the bundle $E_U$ can be extended to $X$ implies that $\pi^*E_U\to P\setminus\{c\}$ can be trivialized. The invariant $\mu_1(E)$ describes the isomorphism type of $E_U$ and can be thought of as specifying, up to homotopy, how to glue the restrictions of $\pi^*E_U\to P\setminus\{c\}$ to the sides of $P$. Choosing a trivialization of $\pi^*E_U$, this is the same as associating, for each $j$, connected components of $G$ to $A_j$ and to $B_j$. This is how $\mu_1$ can be seen as a homomorphism
$$\mu_1:\pi_1X\longrightarrow\pi_0G.$$.

Let $\mathcal{G}(E_U)$ and $\mathcal{G}(E_V)$ be the gauge groups of $E_U$ and $E_V$. Then the relevant gluing information to recover the bundle $E$, up to isomorphism, is the class of $\rho:E_U|_{U\cap V}\to E_V|_{U\cap V}$ in the set of connected components of the double quotient
$\mathcal{G}(E_V)\backslash\mathrm{Isom}(E_U|_{U\cap V},E_V|_{U\cap V})/\mathcal{G}(E_U)$, that is, in
\begin{equation}\label{double quotient}
\pi_0(\mathcal{G}(E_V))\backslash\pi_0(\mathrm{Isom}(E_U|_{U\cap V},E_V|_{U\cap V}))/\pi_0(\mathcal{G}(E_U)).
\end{equation}
Choosing adequate trivializations of $\pi^*(E_U|_{U\cap V})$ and of $E_V|_{U\cap V}$, the map $\rho$ is given by a map 
$$\rho_0:U\cap V\longrightarrow G_0$$
 (note that $U\cap V\cong P\setminus\{c\}\cap B(c,\epsilon)$). Since $U\cap V\sim S^1$ and $\pi_1(G_0)$ is abelian, we can identify $\rho_0$ with an element, still denoted by $\rho_0$, of $\pi_1(G_0)=\pi_1G$ (we define the fundamental group of a topological group as the fundamental group of its identity component): $$\rho_0\in\pi_1G.$$

Now, recall that $\pi_0G$ acts on $\pi_1G$ via the adjoint action (hence $[\alpha_1]=[\alpha_2]$ in $\pi_1G/\pi_0G$ if and only if there is $a\in G$ such that $\alpha_1$ and $a\alpha_2a^{-1}$ are homotopic). Since $\pi_1G$ is an abelian group, we denote the group structure additively.
Given $\mu_1:\pi_1X\to\pi_0G$, define $\Gamma_{\mu_1}\subset\pi_1G$ as the subgroup of $\pi_1G$ generated by the elements of the form $\gamma_2-\gamma_1\cdot\gamma_2$, where $\gamma_2\in\pi_1G$ and $\gamma_1$ lies in the image of $\mu_1$:
\begin{equation}\label{Gammamu1}
\Gamma_{\mu_1}=\left\langle\gamma_2-\gamma_1\cdot\gamma_2\,\st\,\gamma_2\in\pi_1G,\,\gamma_1\in\ima(\mu_1)\subseteq\pi_0G\,\right\rangle.
\end{equation}
Then, since $\pi_0G$ is abelian, the action of $\pi_0G$ on $\pi_1G$ descends to the quotient $\pi_1G/\Gamma_{\mu_1}$.

\begin{definition}\label{mu2}
The topological invariant $\mu_2$ of $E$ is defined as the class of $\rho_0\in\pi_1G$ in $(\pi_1G/\Gamma_{\mu_1})/\pi_0G$:
$$\mu_2(E)=[\rho_0]\in(\pi_1G/\Gamma_{\mu_1})/\pi_0G.$$
\end{definition}

It should be noticed that the values which the invariant $\mu_2$ can take depend on the invariant $\mu_1$.

Similar arguments to the ones used in Proposition 5.1 of \cite{ramanathan:1975} show that the pair $(\mu_1,\mu_2)$ is well-defined and that uniquely characterizes the bundle $E$. In terms of (\ref{double quotient}) this can understood as follows:
\begin{itemize}
 \item $\pi_1G$ corresponds to $\pi_0(\mathrm{Isom}(E_U|_{U\cap V},E_V|_{U\cap V}))$.
\item the action of $\Gamma_{\mu_1}$ on $\pi_1G$ corresponds to the action of $\pi_0(\mathcal{G}(E_U))$ on $\pi_0(\mathrm{Isom}(E_U|_{U\cap V},E_V|_{U\cap V}))$.
\item the action of $\pi_0G$ on $\pi_1G$ corresponds to the action of $\pi_0(\mathcal{G}(E_V))$ on $\pi_0(\mathrm{Isom}(E_U|_{U\cap V},E_V|_{U\cap V}))$.
\end{itemize}

We have therefore the following topological classification of $G$-principal bundles over closed oriented surfaces.
\begin{proposition}\label{classtop}
 Let $X$ be a closed, oriented surface and let $G$ be a topological group such that $\pi_0G$ is abelian.
Given $\mu_1\in H^1(X,\pi_0G)$, there is a bijection between the set of isomorphism classes of continuous $G$-principal bundles $E$ over $X$, with $\mu_1(E)=\mu_1$, and $(\pi_1G/\Gamma_{\mu_1})/\pi_0G$.
 \end{proposition}

\begin{remark}
 In case $G$ is connected, this classification coincides with the well-known topological classification of $G$ bundles over $X$, given $\pi_1G$ (cf. \cite{ramanathan:1975}).
\end{remark}
\begin{remark}
For the case of the sphere $S^2$ (in fact for $S^n$) this was already known (cf. \cite{steenrod:1999}, Section 18).
\end{remark}
\begin{remark} The same result is valid not only for closed, oriented surfaces, but also for any $2$-dimensional connected CW-complex. A proof of this fact, using different methods, can be found in \cite{oliveira:2008}.
\end{remark}

\subsection{The case of $\PGL(n,\R)$}
Now we shall apply the result obtained in the previous section to obtain invariants which classify continuous $\PGL(n,\R)$-principal bundles over our surface $X$.

As $\PGL(n,\R)$ is homotopically equivalent to $\PO(n,\R)=\Or(n,\R)/\Z_2$, its maximal compact subgroup, this is equivalent to classify $\PO(n,\R)$-bundles. From now on, we will write $\PO(n)$ instead of $\PO(n,\R)$ for the real projective orthogonal group, as well as $\Or(n)$ instead of $\Or(n,\R)$ for the real orthogonal group.

For $\PO(n)$, we have that
\begin{equation}\label{mu1 in terms of Cech - PO(n)}
\mu_1\in H^1(X,\pi_0\PO(n))\cong(\Z_2)^{2g}.
\end{equation}
This class is the obstruction to reduce the structure group to $\PSO(n)$.

For $n\geq 4$ even,
$$\pi_1\PO(n)=\begin{cases}
    \Z_2\times\Z_2 & \text{if } n=0\ \text{mod}\ 4 \\
    \Z_4 & \text{if } n=2\ \text{mod}\ 4.
  \end{cases}$$
More precisely, the universal cover of $\PO(n)$ is
$\Pin(n)$ and, if $p:\Pin(n)\to\PO(n)$ is the covering
  projection, then, as a set, $\ker(p)=\{1,-1,\omega_n,-\omega_n\}$ where $\omega_n=e_1\cdots e_n$ is the oriented volume element of $\Pin(n)$ in the standard construction of this group via the Clifford algebra $\mathrm{Cl}(n)$ (see, for example, \cite{lawson-michelson:1989}).

\begin{notation}\label{notadditive}
From now on we shall use the additive notation for $\{1,-1,\omega_n,-\omega_n\}$. Hence, under this notation, $\{1,-1,\omega_n,-\omega_n\}=\{0,1,\omega_n,-\omega_n\}$ (so $1$ becomes $0$ and $-1$ becomes $1$). This is done because we will identify $\{0,1,\omega_n,-\omega_n\}$ with $\pi_1\PO(n)$ which is an abelian group.
\end{notation}

  Recall that $\Pin(n)$ is a group with two connected components, $\Pin(n)^-$ and $\text{Spin}(n)$, where $\Pin(n)^-$ denotes the component which does not contain the identity. We have $\pm\omega_n\notin\text{Z}(\Pin(n))=\{0,1\}$, so the action of $\pi_0\PO(n)$ on $\pi_1\PO(n)$ is not trivial. In fact, $\omega_n$ commutes with elements in $\text{Spin}(n)$ and
  anti-commutes with elements in $\Pin(n)^-$, so $$\pi_1\PO(n)/\pi_0\PO(n)=\{0,1,\omega_n\}$$ where we also write $\omega_n$ for the class of $\omega_n\in\pi_1\PO(n)$ in $\pi_1\PO(n)/\pi_0\PO(n)$, which consists by $\pm\omega$.

For $\PO(n)$-bundles with $\mu_1=0$, we have $\Gamma_0=0$, where $\Gamma_0$ is the subgroup of $\pi_1\PO(n)$ defined in the general setting in (\ref{Gammamu1}).

For $\PO(n)$-bundles with $\mu_1\neq 0$, then it is easy to see that $\Gamma_{\mu_1}=\{0,1\}\cong\Z_2$, therefore
$$\pi_1\PO(n)/\Gamma_{\mu_1}=\{0,\omega_n\}\cong\Z_2,$$ and $\pi_0\PO(n)$ acts trivially on this quotient:
$$(\pi_1\PO(n)/\Gamma_{\mu_1})/\pi_0\PO(n)=\{0,\omega_n\}\cong\Z_2.$$
Hence, we have the invariant $\mu_2$ defined in general in Definition \ref{mu2}, which, for $\PO(n)$-principal bundles over $X$, is such that:
\begin{equation}\label{mu2 PO(n)}
\mu_2\in(\pi_1\PO(n)/\Gamma_{\mu_1})/\pi_0\PO(n)=\begin{cases}
\{0,1,\omega_n\}\ &\text{ if }\ \mu_1=0\\
\{0,\omega_n\}\ &\text{ if }\ \mu_1\neq 0
                                \end{cases}.
\end{equation}

\begin{remark}\label{attention with notation}
When $\mu_1\neq 0$, we also write the possible elements of $\mu_2\in\Z_2$ by $0$ and by $\omega_n$, instead of $[0]$ and $[\omega_n]$. This requires a little attention because, for example, $\mu_2=0$ has different meanings whenever $\mu_1=0$ or $\mu_1\neq 0$. However, it should always be clear in which situation we are.
\end{remark}

\begin{remark}\label{PGL+}
When $\mu_1=0$, we are reduced to the topological classification of
$\PSO(n)$-bundles over $X$ which, for $\PSO(n)$-equivalence, is given by the elements in $\{0,1,\omega_n,-\omega_n\}=\pi_1\PSO(n)$. However, since we are interested in $\PO(n)$-equivalence, the bundles
with invariants $\omega_n$ and $-\omega_n$ become identified.
\end{remark}

The next proposition gives the interpretation of the class $\mu_2$
in terms of obstructions.
\begin{proposition}\label{obs}
Let $n\geq 4$ be even.
\begin{enumerate}
 \item Let $E$ be a continuous $\PO(n)$-bundle over $X$ with $\mu_1(E)=0$. Then:
\begin{itemize}
    \item $E$ lifts to a continuous $\SO(n)$-bundle if and only if $\mu_2(E)\in\{0,1\}$;
    \item $E$ lifts to a continuous $\Spin(n)$-bundle if and only if $\mu_2(E)=0$.
\end{itemize}
\item Let $E$ be a continuous $\PO(n)$-bundle over $X$ with $\mu_1(E)\neq 0$. Then $E$ lifts to a continuous $\Pin(n)$-bundle if and only if $\mu_2(E)=0$.
\end{enumerate}
\end{proposition}
\begin{proof} Suppose $\mu_1(E)=0$, so that $E$ is in fact a $\PSO(n)$-bundle. From the construction of $\mu_2$ in the previous subsection, we have
$$\mu_2(E)=[g]\in\pi_1\PO(n)/\pi_0\PO(n),$$ where
$g:(S^1,y_0)\to(\PO(n),[I_n])$. Let
$p:\Or(n)\to\PO(n)$ be the projection. There is
a lift $g':(S^1,y_0)\to(\mathrm{O}(n),I_n)$ if and
only if $g_*(\pi_1S^1)\subseteq p_*(\pi_1\Or(n))$,
which happens if and only if $[g]\in\{0,1\}$.
The case for the lift to $\Pin(n)$ is completely analogous.

The case of $\mu_1(E)\neq 0$ is proved in a similar way, noticing also that over the $1$-skeleton $X_1$ of $X$ there are no
obstructions to lifting the bundle because there the bundle is trivialized on contractible open sets.
\end{proof}

\begin{remark}
Notice that, when $\mu_1\neq 0$, a $\PO(n)$-bundle lifts to an $\Or(n)$-bundle if and only if ot lifts to a $\Pin(n)$-bundle. This is clear since, when $\mu_1\neq 0$, the $0$ in $\{0,\omega_n\}$ is the class of $0$ and $1$ in the quotient $(\pi_1PO(n)/\Gamma_{\mu_1})/\pi_0\PO(n)$ (cf. Remark \ref{attention with notation}). 

Another way to see that a $\PO(n)$-bundle lifts to a $\Pin(n)$-bundle if it lifts to an $\Or(n)$-bundle is as follows.
Suppose that $E$ is a real projective bundle, with $\mu_1(E)\neq 0$, and which is the projectivization of a real vector bundle $W$. Since the projection from $\Or(n)$ onto $\PO(n)$ preserves components of the groups (because $n$ is even), $w_1(W)=\mu_1(E)\neq 0$ where
$w_1(W)$ is the first Stiefel-Whitney class of $W$. So the first Stiefel-Whitney class of all lifts of $E$ to $\Or(n)$ is the same (another way to see this is to note that $w_1(W\otimes F)=w_1(W)$, for any real line bundle $F$, whenever $\rk(W)$ is even).
Nevertheless, different lifts of $E$ can have different second Stiefel-Whitney class because their first Stiefel-Whitney class is non-zero. In fact, given a real vector bundle $W$ of rank $n$ on $X$ with $w_1(W)\neq 0$, it is easy to see that there exists a real line bundle $F$ such that $w_2(W)\neq w_2(W\otimes F)$ (note that
$w_2(W\otimes F)=w_2(W)+w_1(W)w_1(F)$). This is the reason why the second Stiefel-Whitney class ``disappears'' on projective bundles with $\mu_1\neq 0$.
Hence either $W$ or $W\otimes F$ has $w_2=0$ and therefore lifts to a $\Pin(n)$-bundle. Choosing this lift of $E$, we see that $E$ lifts to a $\Pin(n)$-bundle.
\end{remark}

\begin{remark}\label{SWclasses}
If $E$ is a real projective bundle with $\mu_1(E)=0$ and $\mu_2(E)\in\{0,1\}$, then also the second Stiefel-Whitney class of the lifts is well defined and is equal to $\mu_2$.
\end{remark}

From Proposition \ref{classtop}, we obtain a full topological classification of real projective bundles over $X$.

\begin{theorem}\label{classtoppgl}
Let $n\geq 4$ be even, and let $X$ be a closed oriented surface of genus $g\geq 2$. Then continuous $\PO(n)$-bundles over $X$
are classified by
$$(\mu_1,\mu_2)\in\left(\{0\}\times\{0,1,\omega_n\}\right)\cup\left(\left((\Z_2)^{2g}\setminus\{0\}\right)\times\Z_2\right).$$
\end{theorem}

%%%%%%%%%%%%%%%%%%%%%%%%%%%%%%%%%%%%%%%%%%%%%%%%%%%%%%%%%%%%%%%
\section{Representations and topological classification}\label{Representations and topological classification}
%%%%%%%%%%%%%%%%%%%%%%%%%%%%%%%%%%%%%%%%%%%%%%%%%%%%%%%%%%%%%%%
\subsection{Representations of $\pi_1X$ in $\PGL(n,\R)$}\label{Rep}

In this section we begin our analysis of the space $\cR_{\PGL(n,\R)}$. The first thing to do is to define a topological invariant of a representation $\rho:\pi_1X\to\PGL(n,\R)$. From Definition \ref{def top invariant representation}, we already know that is done via the correspondence between representations and flat bundles.

\begin{definition}\label{inv of rep} Let $\rho$ be a representation $\pi_1X\to
\PGL(n,\R)$ and let $E_\rho=\widetilde{ X}\times_\rho\PGL(n,\R)$, the principal flat
$\PGL(n,\R)$-bundle over $X$ associated to $\rho$, viewed as a continuous bundle. The \emph{topological invariants $\mu_1(\rho)$ and $\mu_2(\rho)$
of $\rho$} are defined by $\mu_1(\rho)=\mu_1(E_\rho)$ and $\mu_2(\rho)=\mu_2(E_\rho)$ where $\mu_1(E_\rho)$ and $\mu_2(E_\rho)$ are the invariants defined in (\ref{mu1 in terms of Cech - PO(n)}) and (\ref{mu2 PO(n)}).
Thus $$\mu_1(\rho)=\mu_1(E_\rho)\in\Z_2^{2g}\text{ and }\mu_2(\rho)=\mu_2(E_\rho)=\begin{cases}
    \{0,1,\omega_n\} & \text{if }\mu_1(\rho)=0\\
\{0,\omega_n\} & \text{if }\mu_1(\rho)\neq 0.
  \end{cases}$$
\end{definition}

Recall that our goal is to determine the number of connected components of
$$\cR_{\PGL(n,\R)}=\Hom^\text{red}(\pi_1X,\PGL(n,\R))/\PGL(n,\R)$$ for $n\geq 4$ even.

For fixed topological invariants,
$$(\mu_1,\mu_2)\in\left(\{0\}\times\{0,1,\omega_n\}\right)\cup\left(\left((\Z_2)^{2g}\setminus\{0\}\right)\times\Z_2\right),$$
we define the subspace $\cR_{\PGL(n,\R)}(\mu_1,\mu_2)$
of $\cR_{\PGL(n,\R)}$ as
$$\cR_{\PGL(n,\R)}(\mu_1,\mu_2)=\{\rho\st\mu_i(\rho)=\mu_i, i=1,2\}.$$

\subsection{Non-emptiness of $\cR_{\PGL(n,\R)}(\mu_1,\mu_2)$}
For fixed invariants $(\mu_1,\mu_2)$, we will now study the non-emptiness of $\cR_{\PGL(n,\R)}(\mu_1,\mu_2)$. To do so, we will see how to detect the classes $\mu_1$ and $\mu_2$ of a flat
$\PGL(n,\R)$-bundle, using only the corresponding representation of $\pi_1X$
in $\PGL(n,\R)$.

Let $\PGL(n,\R)_0$ denote the identity component of $\PGL(n,\R)$ and let $\PGL(n,\R)^-$ denote the component of $\PGL(n,\R)$ which does not contain the identity.

\begin{definition}\label{delta1rho}
Given a representation $\rho:\pi_1X\to\PGL(n,\R)$, let $A_1,B_1,\dots,B_g\in\PGL(n,\R)$ be the
images of the generators of $\pi_1X$ by $\rho$. The \emph{invariant $\delta_1$} of $\rho$ is defined as
$$\delta_1(\rho)\in(\Z_2)^{2g}$$ to be such that:
\begin{itemize}
 \item the $(2i-1)$-th coordinate of $\delta_1(\rho)$ is $0$ if $A_i\in\PGL(n,\R)_0$;
 \item the $(2i-1)$-th coordinate of $\delta_1(\rho)$ is $1$ if $A_i\in\PGL(n,\R)^-$;
 \item the $2i$-th coordinate of $\delta_1(\rho)$ is $0$ if $B_i\in\PGL(n,\R)_0$;
 \item the $2i$-th coordinate of $\delta_1(\rho)$ is $1$ if $B_i\in\PGL(n,\R)^-$.
\end{itemize}
\end{definition}

 Obviously, $\delta_1(\rho)$ is the
obstruction to reducing the representation to $\PGL(n,\R)_0$. So we
have
\begin{equation}\label{delta1}
 \delta_1(\rho)=\mu_1(\rho).
\end{equation}

\vspace{0,5cm}
In order to obtain something similar for the invariant $\mu_2$, we will consider representations in the maximal compact $\PO(n)$. In terms of topological invariants, there is no loss of generality in doing this and has the advantage that these representations are automatically reductive due to the compactness of $\PO(n)$.

Let $p':\Or(n)\to\PO(n)$ be the projection. Choose $A_i'\in p'^{-1}(A_i)$
and $B_i'\in p'^{-1}(B_i)$ in $\Or(n)$, and consider the product $$\prod_{i=1}^g[A_i',B_i'].$$
Since $\ker(p')\subseteq Z(\Or(n))$, the value of this product does not depend on the choice of the lifts $A_i', B_i'$ and it is the obstruction to lifting $\rho:\pi_1X\to\PO(n)$ to a representation $\rho':\pi_1X\to\Or(n)$.
\begin{definition}\label{delta2rho}
Let $\rho:\pi_1X\to\PO(n)$ be a
representation and let $A_1,B_1,\dots,B_g\in\PO(n)$ be the
images of the generators of $\pi_1X$ by $\rho$. The \emph{invariant $\delta_2$} of $\rho$ is defined as
$$\delta_2(\rho)=\prod_{i=1}^g[A_i',B_i']\in\{\pm I_n\}$$ where $A_i'$ and $B_i'$ are lifts of $A_i$ and $B_i$, respectively, to $\Or(n)$.
\end{definition}

\begin{remark}
In this remark (and only here) we will \emph{not} use the additive notation of Notation \ref{notadditive}, since here we are going to work on the $\Pin(n)$ and $\Spin(n)$ group (which are not abelian).
If $\delta_2(\rho)=I_n$, one can ask whether $\rho':\pi_1X\to\Or(n)$ lifts to a representation
$\rho'':\pi_1X\to\Pin(n)$ under the projection $p'':\Pin(n)\to\Or(n)$ and the way to measure the obstruction to the existence of this lift is exactly the same as in the previous case: choose lifts
$\tilde{A_i}\in p''^{-1}(A_i')$ and $\tilde{B_i}\in p''^{-1}(B_i')$, for
all $i\in\{1,\dots,g\}$, and consider
the value
\begin{equation}\label{liftpin}
 \prod_{i=1}^g[\tilde{A_i},\tilde{B_i}]\in\{\pm 1\}=p''^{-1}(I_n).
\end{equation}

Again this is well-defined because $\ker(p'')\subseteq Z(\Pin(n))$ and it is the obstruction
to lifting $\rho'$ to a representation $\rho'':\pi_1X\to\Pin(n)$.

If $\tilde{p}:\Pin(n)\to\PO(n)$ is the universal cover ($\tilde{p}=p'\circ p''$)
then, in the case $\delta_1(\rho)\neq 0$, we could not use the same procedure as in the previous cases to measure directly
the obstruction to lifting $\rho$ to a representation $\tilde{\rho}:\pi_1X\to\Pin(n)$ because $\ker(\tilde{p})=\{\pm1,\pm\omega_n\}\not\subset Z(\Pin(n))=\{\pm 1\}$. In principle, the
above procedure only gives partial information about the possible 
lifts of $\rho$ to $\Pin(n)$: if $\delta_2(\rho)=-I_n$ then clearly $\rho$ does not lift
to $\Pin(n)$; if $\delta_2(\rho)=I_n$ and the lift $\rho'$ of $\rho$ to $\Or(n)$ lifts to $\Pin(n)$ then $\rho$ lifts to $\Pin(n)$; if $\delta_2(\rho)=I_n$ but the lift $\rho'$ of $\rho$ to $\Or(n)$ does not lift to $\Pin(n)$, we cannot conclude that $\rho$ does not lift to $\Pin(n)$ because if we change the lift
of $\rho$ to $\Or(n)$ (or, equivalently, if we change the lifts of some of the generators
$A_i$ and $B_i$) then this new representation of $\pi_1X$ on $\Or(n)$ might
lift to $\Pin(n)$. In fact,
this is always possible, if $\mu_1(\rho)\neq 0$ (i.e., if $\delta_1(\rho)\neq 0$). To see
this, suppose $\rho$ is such that $\delta_1(\rho)\neq 0$ and $\delta_2(\rho)=I_n$. Then $\prod_{i=1}^g[A_i',B_i']=I_n$ for any lifts of $A_i$ and of $B_i$. On the other hand, there is some $A_{i_0}\in\PO(n)^-$, so $A_{i_0}'\in\Or(n)^-$.
If $\prod_{i=1}^g[\tilde{A_i},\tilde{B_i}]=-1$, then $-B_{i_0}'\in\Or(n)$ is other lift of $B_{i_0}$
and, choosing it, we have a new lift of $\rho$ to $\Or(n)$. Lifting $-B_{i_0}'$ to $\Pin(n)$ we obtain
$\pm\omega_n\tilde{B_{i_0}}$ and now, since $\tilde A_{i_0}^{-1}\in\Pin(n)^-$ and since $\tilde{B_{i_0}}$ and 
$\tilde B_{i_0}^{-1}$ belong to the same component of $\Pin(n)$, we have
$$\tilde A_1\tilde B_1\tilde A_1^{-1}\tilde B_1 ^{-1}\cdots\tilde A_{i_0} \omega_n\tilde B_{i_0}
\tilde A_{i_0}^{-1}\tilde B_{i_0}^{-1}\omega_n^{-1}\cdots\tilde A_g\tilde B_g\tilde A_g^{-1}\tilde B_g^{-1}=
-\prod_{i=1}^g[\tilde A_i,\tilde B_i]=1.$$
Thus, if $\delta_1(\rho)\neq 0$, the value of (\ref{liftpin}) does not give any new information.
\end{remark}

If $\delta_1(\rho)=0$, $\rho$ reduces to a representation in $\PSO(n)$ and, as $\ker(\tilde{p})\subseteq Z(\Spin(n))$ where $\tilde{p}:\Spin(n)\to\PSO(n)$,
we have a well defined obstruction $\tilde{\delta}(\rho)$ to lifting $\rho$ to $\Spin(n)$, defined as follows:

\begin{definition}\label{tildedeltarho}
Let $n\geq 4$ be even. Let $\rho:\pi_1X\to\PO(n)$ be a
representation with $\delta_1(\rho)=0$ and let $A_1,B_1,\dots,B_g\in\PO(n)$ be the
images of the generators of $\pi_1X$ by $\rho$. The \emph{invariant $\tilde\delta$} of $\rho$ is defined as
$$\tilde{\delta}(\rho)=\prod_{i=1}^g[\tilde A_i,\tilde B_i]\in\{0,1,\omega_n\}$$
where $\tilde A_i$ and $\tilde B_i$ are lifts of $A_i$ and $B_i$, respectively, to $\Spin(n)$.
\end{definition}
Again, in $\{0,1,\omega_n\}$ of this definition we have identified $\omega_n$ and $-\omega_n$ due to the $\PO(n)$-equivalence.

Recall Definition \ref{inv of rep}. From Proposition \ref{obs} and from what we have seen, we have the following lemma:

\begin{lemma}\label{equiv between inv} Let $n\geq 4$ be even. The following equivalences hold:
$$\delta_2(\rho)=-I_n\Longleftrightarrow\mu_2(\rho)=\omega_n$$
and, if $\delta_1(\rho)\neq 0$,
$$\delta_2(\rho)=I_n\Longleftrightarrow\mu_2(\rho)=0.$$
If $\delta_1(\rho)=0$, we have
$$\tilde{\delta}(\rho)=\mu_2(\rho)\in\{0,1,\omega_n\}.$$
\end{lemma}

\begin{proposition}\label{nonempty}
Let $n\geq 4$ even be given.
Then, the space
$\cR_{\PGL(n,\R)}(\mu_1,\mu_2)$ is non-empty, for each pair
$(\mu_1,\mu_2)\in\left(\{0\}\times\{0,1,\omega_n\}\right)\cup\left(\left((\Z_2)^{2g}\setminus\{0\}\right)\times\Z_2\right)$.
\end{proposition}
\begin{proof}
Let us start by seeing that $\cR_{\PGL(n,\R)}(\mu_1,\omega_n)$ is
non-empty for each $\mu_1\in(\Z_2)^{2g}$. To do so we will find an explicit representation of $\pi_1X$ in $\PO(n)$ (hence in $\PGL(n,\R)$) 
with these invariants.
From (\ref{delta1}) and Lemma \ref{equiv between inv}, in order to show
that $\cR_{\PGL(n,\R)}(\mu_1,\omega_n)$
is non-empty we only need to find a reductive representation $\rho:\pi_1X\to\PO(n)\subset\PGL(n,\R)$
with $$\delta_1(\rho)=\mu_1$$ and $$\delta_2(\rho)=-I_n.$$
In other words, from the definition of $\delta_1(\rho)$, we need to find $n\times n$ invertible
matrices $A_i'$ and $B_i'$ such that
\begin{itemize}
 \item $A_i'\in\SO(n)$ if and only if the $(2i-1)$-th coordinate of $\delta_1(\rho)$ is $0$;
 \item $A_i'\in\Or(n)^-$ if and only if the $(2i-1)$-th coordinate of $\delta_1(\rho)$ is $1$;
 \item $B_i'\in\SO(n)$ if and only if the $2i$-th coordinate of $\delta_1(\rho)$ is $0$;
 \item $B_i'\in\Or(n)^-$ if and only if the $2i$-th coordinate of $\delta_1(\rho)$ is $1$.
\end{itemize}
and, from the definition of $\delta_2(\rho)$, which satisfy the equality
$$\prod_{i=1}^g[A_i',B_i']=-I_n.$$ As we are using the compact group $\PO(n)$,
the reductiveness condition on the representation is automatically satisfied.

Let us start with the
following orthogonal matrices:
$$X_2=\begin{pmatrix}
     0 & 1 \\
     1 & 0
\end{pmatrix},\ X'_2=\begin{pmatrix}
     1 & 0 \\
     0 & -1
\end{pmatrix},\ Y_2=X'_2,\ Y'_2=-X'_2\text{ and }
Z_2=\begin{pmatrix}
     0 & -1 \\
     1 & 0
\end{pmatrix}.$$
Note that $X_2$, $X'_2$ and $Z_2$ are pairwise anti-commuting and
that $Y_2$ and $Y'_2$ commute. For $n\geq 4$ even, define
$$X_n=\begin{pmatrix}
     X_2 & 0 \\
     0 & X_{n-2}
\end{pmatrix},\
X'_n=\begin{pmatrix}
     X'_2 & 0 \\
     0 & X'_{n-2}
\end{pmatrix},\ Y_n=\begin{pmatrix}
     Y_2 & 0 \\
     0 & I_{n-2}
\end{pmatrix},\ Y'_n=\begin{pmatrix}
     Y'_2 & 0 \\
     0 & I_{n-2}
\end{pmatrix}$$
$$Z_n=\begin{pmatrix}
     Z_2 & 0 \\
     0 & X'_{n-2}
\end{pmatrix},\ W_n=\begin{pmatrix}
     X_2 & 0 \\
     0 & Z_{n-2}
\end{pmatrix}\text{ and }\ W'_n=\begin{pmatrix}
     Z_2 & 0 \\
     0 & X_{n-2}
\end{pmatrix}.$$

We have the following facts:
\begin{itemize}
 \item $X_n,X'_n\in\SO(n)\Longleftrightarrow n=0\ \mathrm{mod}\ 4$;
\item $X_n,X'_n\in\Or(n)^-\Longleftrightarrow  n=2\ \mathrm{mod}\ 4$;
\item $Y_n,Y'_n\in\Or(n)^-$ for all $n$ even;
 \item $Z_n\in\Or(n)^-\Longleftrightarrow n=0\ \mathrm{mod}\ 4$;
\item $Z_n\in\SO(n)\Longleftrightarrow n=2\ \mathrm{mod}\ 4$;
\item $W_n,W'_n\in\SO(n)\Longleftrightarrow n=2\ \mathrm{mod}\ 4$;
\item $W_n,W'_n\in\Or(n)^-\Longleftrightarrow n=0\ \mathrm{mod}\ 4$;
\item $X_n$ and $X'_n$ anti-commute for all $n$ even;
 \item $Y_n$ and $Y'_n$ commute for all $n$ even;
 \item $Z_n$ anti-commutes with $X_n$ for all $n$ even;
\item $W_n$ and $W'_n$ anti-commute for all $n\geq 4$ even.
\end{itemize}

Using these orthogonal matrices and the identity $I_n$ it is
possible to construct the required representation. The important thing to note is that for each $n$ we
always have a pair of commuting and anti-commuting matrices both
in $\SO(n)$ or both in $\Or(n)^-$ or one in $\SO(n)$ and the
other in $\Or(n)^-$.

The case of commuting matrices is easy:
if one of the matrices is to be in $\SO(n)$, use the identity $I_n$; if both
must be in $\Or(n)^-$, use $Y_n$ and $Y'_n$:

\vspace{0,5cm}
%\begin{table}[ht]
%\tbl{Commuting matrices.}
\begin{tabular}{|c|c|c|c|}
  \hline
  Commuting matrices & $\SO(n)$, $\SO(n)$ & $\SO(n)$, $\Or(n)^-$ & $\Or(n)^-$, $\Or(n)^-$ \\
  \hline
  $n$ even & $I_n$, any & $I_n$, any & $Y_n$, $Y'_n$ \\
  \hline
\end{tabular}
%\end{table}
\vspace{0,5cm}

The case of anti-commuting matrices is also easy:

\vspace{0,5cm}
%\begin{table}[ht]
%\tbl{Anti-commuting matrices.}
\begin{tabular}{|c|c|c|c|}
  \hline
  Anti-commuting matrices & $\SO(n)$, $\SO(n)$ & $\SO(n)$, $\Or(n)^-$ & $\Or(n)^-$, $\Or(n)^-$ \\
  \hline
  $n=0\ \text{mod}\,4$ & $X_n$, $X'_n$ & $X_n$, $Z_n$ & $W_n$, $W'_n$ \\
\hline
  $n=2\ \text{mod}\,4$ & $W_n$, $W'_n$ & $Z_n$, $X_n$ & $X_n$, $X'_n$ \\
  \hline
\end{tabular}
%\end{table}

\vspace{0,5cm}

 This shows that given
$\mu_1=\delta_1=(x_1,x_2,\dots,x_{2g})\in(\Z_2)^{2g}$, we can
choose $A'_1$ and $B'_1$ in $\SO(n)$ or in $\Or(n)^-$
depending on $x_1$ and on $x_2$ and such that $[A'_1,B'_1]=-I_n$.
Then, for $i\geq 2$, we can also choose $A'_i$ and $B'_i$
accordingly to $x_{2i-1}$ and to $x_{2i}$ respectively, and such
that $[A'_i,B'_i]=I_n$. Hence $\prod_{i=1}^g[A'_i,B'_i]=-I_n$ as
wanted. Putting $\rho(a_i)=p_2(A'_i)$ and $\rho(b_i)=p_2(B'_i)$
gives a representation $\rho:\pi_1X\to\PGL(n,\R)$ with the given
$\mu_1$ and $\mu_2(\rho)=\omega_n$.

For the other cases, the proof is similar but easier. For $\cR_{\PGL(n,\R)}(\mu_1,0)$ with $\mu_1\neq 0$, we have, from (\ref{delta1}) and Lemma \ref{equiv between inv}, to find a representation $\rho$ with $\delta_1(\rho)=\mu_1$ and $\delta_2(\rho)=I_n$ and this is done in same way as above, using the first table.

The cases $\cR_{\PGL(n,\R)}(0,\mu_2)$ with $\mu_2=0,1$ should be dealt with similarly, but now we would need to consider the invariant $\tilde\delta$ of Definition \ref{tildedeltarho} and, hence, elements on the $\Pin(n)$ group. Instead, note that, since $\mu_1=0$, we are looking for representations on the connected group $\PGL(n,\R)_0$. Hence, the non-emptiness of $\cR_{\PGL(n,\R)}(0,\mu_2)$ follows from Proposition $7.7$ of \cite{ramanathan:1975}.
\end{proof}

The map in $\cR_{\PGL(n,\R)}$ which takes a class $\rho$ to $(\mu_1(\rho),\mu_2(\rho))$ is continuous hence, if classes lie in the same connected component of $\cR_{\PGL(n,\R)}$, they must have the same topological invariants.
From this and from Theorem \ref{classtoppgl} and Proposition \ref{nonempty}, we conclude that, for $n\geq 4$ even, $\cR_{\PGL(n,\R)}$ has at least $2^{2g+1}+1$ connected components this being the
number of topological invariants.

It remains to see whether, for each pair $(\mu_1,\mu_2)$,
$\cR_{\PGL(n,\R)}(\mu_1,\mu_2)$ is connected or not, and
it is now that the theory of Higgs bundles comes into play.

%%%%%%%%%%%%%%%%%%%%%%%%%%%%%%%%%%%%%%%%%%%%%%%%%%%%%%%%%%%%%%%%%%%%%%%%%%%%%%%%%%%%%%%%%%%%%%%%%%%%%%%%
\section{$\PGL(n,\R)$-Higgs bundles and $\EGL(n,\R)$-Higgs bundles}\label{section PGL(n,R)-Higgs bundles and quadruples}
%%%%%%%%%%%%%%%%%%%%%%%%%%%%%%%%%%%%%%%%%%%%%%%%%%%%%%%%%%%%%%%%%%%%%%%%%%%%%%%%%%%%%%%%%%%%%%%%%%%%%%%%

In this section we begin the study of $\PGL(n,\R)$-Higgs bundles and explain why and how one wants to work with another group instead of $\PGL(n,\R)$.

We begin by defining $\PGL(n,\R)$-Higgs bundles, using Definition \ref{definition of Higgs bundle}. Recall that $\PO(n)^\C=\PO(n,\C)=\Or(n,\C)/\Z_2$.
\begin{definition}
A \emph{$\PGL(n,\R)$-Higgs bundle} over $X$ is a pair $(E,\Phi)$,
where $E$ is a holomorphic principal $\PO(n,\C)$-bundle and $\Phi\in
H^0(X,E\times_{\PO(n,\C)}\mathfrak{so}(n,\C)^\perp\otimes K)$ where
$\mathfrak{so}(n,\C)^\perp$ is the vector space of
$n\times n$ symmetric and traceless complex matrices.
\end{definition}

We would like to work naturally with holomorphic vector
bundles associated to the corresponding $\PGL(n,\R)$-Higgs
bundles. However, this will not be done directly because
$\PO(n,\C)$ does not have a standard action on $\C^n$, and to fix this we use a standard procedure as follows. Enlarge the
complex orthogonal group $\Or(n,\C)$ so that it still has a canonical action on
$\C^n$ and such that it has a non-discrete centre, and consider the sheaf of holomorphic functions with values in this centre. Then, the second cohomology of $X$ of this sheaf vanishes, so there is no obstruction to lifting a holomorphic $\PO(n,\C)$-bundle to a
holomorphic bundle with this new structure group, and hence to do the
same to Higgs bundles with the corresponding groups.

Let us then consider the group $\GL(n,\R)\times\U(1)$, the normal subgroup $\{(I_n,1),(-I_n,-1)\}\cong\Z_2\vartriangleleft \GL(n,\R)\times\U(1)$ and the corresponding quotient group $$\GL(n,\R)\times_{\Z_2}\U(1)=(\GL(n,\R)\times\U(1))/\Z_2.$$

Its maximal compact is $\Or(n)\times_{\Z_2}\U(1)$, whose complexification is $\Or(n,\C)\times_{\Z_2}\C^*$.

\begin{notation}\label{notation enhanced}
From now on, we shall write $$\EGL(n,\R)=\GL(n,\R)\times_{\Z_2}\U(1)$$
and $$\EO(n)=\Or(n)\times_{\Z_2}\U(1)$$ as well as $$\EO(n,\C)=\Or(n,\C)\times_{\Z_2}\C^*.$$
The ``E'' stands for enhanced or extended.
\end{notation}

Applying again Definition \ref{definition of Higgs bundle}, we give now a concrete definition of $\EGL(n,\R)$-Higgs bundle.
Notice that, if $\overline{G}=\EGL(n,\R)$, then a maximal compact subgroup of $\overline{G}$ is
$\overline{H}=\EO(n)$, so $\overline{H}^{\C}=\EO(n,\C)$. Also,
$\mathfrak{\overline{g}}^{\C}=\mathfrak{\overline{h}}^\C\oplus\mathfrak{\overline{m}}^{\C}$
where
$\mathfrak{\overline{g}}^\C=\mathfrak{gl}(n,\C)\oplus\C$, $\mathfrak{\overline{h}}^\C=\mathfrak{o}(n,\C)\oplus\C$ and
$\overline\liem^\C=\{(A,0)\in\overline\lieg^\C\st A=A^T\}$ is naturally isomorphic to the space of symmetric matrices.
\begin{definition}
A \emph{$\EGL(n,\R)$-Higgs bundle} over $X$ is a pair $(E,\Phi)$,
where $E$ is a holomorphic principal $\EO(n,\C)$-bundle and $\Phi\in
H^0(X,E\times_{\EO(n,\C)}\overline\liem^\C\otimes K)$, where $\overline\liem^\C=\{(A,0)\in\overline\lieg^\C\st A=A^T\}$.
\end{definition}

\begin{proposition}\label{pgl1}
Every $\PGL(n,\R)$-Higgs bundle $(E,\Phi)$ on $X$ lifts to a
$\EGL(n,\R)$-Higgs bundle
$(\overline{E},\overline\Phi)$.
\end{proposition}
\begin{proof}
We have the following short exact sequence of groups
$$0\longrightarrow\C^*\stackrel{i}{\longrightarrow}\EO(n,\C)\stackrel{p}{\longrightarrow}\PO(n,\C)\longrightarrow 0$$
where $i(\lambda)=[(I_n,\lambda)]$ and $p([(w,\lambda)])=[w]$.

Consider the sheaf $\EO(n,\mathcal O)$ of holomorphic functions on $X$ with values in
$\EO(n,\C)$. The above short exact sequence induces the following exact
sequence
\begin{equation}\label{seq}
H^1(X,\mathcal{O}^*)\longrightarrow H^1(X,\EO(n,\mathcal O))\stackrel{p_*}{\longrightarrow}H^1(X,\PO(n,\mathcal{O}))\longrightarrow 0
\end{equation}
hence, we see that there is no obstruction to lifting $E$ to a principal
$\EO(n,\C)$-bundle
$$\overline{E}\in H^1(X,\EO(n,\mathcal O)).$$

Write $G=\PGL(n,\R)$, $H=\PO(n)$, $\overline{G}=\EGL(n,\R)$ and $\overline{H}=\EO(n)$. We have
$\mathfrak{\overline{g}}^{\C}=\mathfrak{\overline{h}}^\C\oplus\mathfrak{\overline{m}}^{\C}$,
where
$$\mathfrak{\overline{g}}^\C=\mathfrak{gl}(n,\C)\oplus\C\supset\mathfrak{sl}(n,\C)\oplus\C=\liegc\oplus\C$$
$$\mathfrak{\overline{h}}^\C=\mathfrak{o}(n,\C)\oplus\C\supset\mathfrak{so}(n,\C)\oplus\C=\liehc\oplus\C$$
and
$$\overline\liem^\C=\{(A,0)\in\overline\lieg^\C\st A=A^T\}.$$ If we identify $\overline\liem^\C$ with $\{A\in\mathfrak{gl}(n,\C)\st A=A^T\}$, then $\liemc=\mathfrak{so}(n,\C)^\perp$ is the subspace of matrices in $\overline\liem^\C$ with trace equal to zero.

Now, the isotropy action of $H^\C=\PO(n,\C)$ in $\liemc$ is given by (where $[A]\in\PO(n,\C)$ and $B\in\liemc$)
\begin{equation}\label{01}
 \Ad([A])(B)=ABA^{-1}=ABA^T
\end{equation}
and the isotropy action of $\overline H^\C$ in $\overline\liem^\C$ is given by (where $[(A,\lambda)]\in\overline H^\C$ and $(B,0)\in\overline\liem^\C$)
\begin{equation}\label{02}
\overline\Ad([(A,\lambda)])(B,0)=(ABA^{-1},0)=(ABA^T,0).
\end{equation}

We have the bundle $\overline E$, which is a lift of $E$ and we have a map $\pi:\overline E\to E$ induced by the projection $\overline H^\C\to H^\C$.
Now, $\Phi\in H^0(X,E\times_{H^\C}\liemc\otimes K)$ can be thought as a $H^\C\times\C^*$-equivariant map $E\times_X E_K\to\liemc$ where $E_K$ is the $\C^*$-principal bundle associated to $K$ (i.e., the frame bundle associated to $K$), and $E\times_X E_K$ is the fibred product over $X$ of $E$ and $E_K$. Let $\overline\Phi:\overline E\times_X E_K\to\overline\liem^\C$ be defined by
\begin{equation}\label{03}
\overline\Phi=(i\otimes 1_K) \Phi(\pi\times_X 1_{E_K})
\end{equation}
where $i:\liemc\hookrightarrow\overline\liem^\C$ is the inclusion. So, we have the commutative diagram 
$$\xymatrix{&\overline E\times_X E_K\ar[d]_{\pi\times_X 1_{E_K}}\ar[r]^(0,6){\overline\Phi}&\overline\liemc\\
    &E\times_X E_K\ar[r]_(0,6){\Phi}&\liemc\ar[u]_{i}.}$$ From (\ref{01}) and (\ref{02}) follows that $\overline\Phi$ is $\overline H^\C\times\C^*$-equivariant.
Hence $\Phi\in H^0(X,E\times_{H^\C}\liemc\otimes K)$ induces, in a natural way, a Higgs field $\overline\Phi\in H^0(X,\overline E\times_{\overline{H}^\C}\overline\liem^\C\otimes K)$ given by (\ref{03}).

It follows that
$(\overline{E},\overline\Phi)$ is an $\EGL(n,\R)$-Higgs
bundle.
\end{proof}

Consider the actions of $\EO(n,\C)$ on $\C^n$ and on $\C$ induced, respectively, by the group homomorphisms
\begin{equation}\label{actCn}
\EO(n,\C)\longrightarrow\GL(n,\C),\ \ \ \ [(w,\lambda)]\mapsto\lambda w 
\end{equation}
and
\begin{equation}\label{actC}
\EO(n,\C)\longrightarrow\C^*,\ \ \ \ [(w,\lambda)]\mapsto\lambda^2.
\end{equation}
and by the corresponding standard actions of $\GL(n,\C)$ and $\C^*$.

\begin{proposition}\label{pgl2}
Let $(\overline{E},\overline\Phi)$ be an $\EGL(n,\R)$-Higgs bundle on $X$.
Through the actions (\ref{actCn}) and (\ref{actC}) of $\EO(n,\C)$ on $\C^n$ and on $\C$, associated
to $(\overline{E},\overline\Phi)$ there is
a quadruple $(V,L,Q,\overline\Phi)$ on $X$, where $V$ is a
holomorphic rank $n$ vector bundle, $L$ is a holomorphic line
bundle, $Q$ is a nowhere degenerate quadratic form on $V$ with values in
$L$ and $\overline\Phi\in H^0(X,S^2_QV\otimes K)$ where $S^2_QV$ denotes the
bundle of endomorphisms of $V$ which are symmetric with respect to
$Q$.
\end{proposition}
\begin{proof}  Keeping the notation of the proof of Proposition \ref{pgl1}, let $\overline{H}^{\C}=\EO(n,\C)$.
From the actions (\ref{actCn}) and (\ref{actC}) we define, respectively,
the vector bundle $V=\overline{E}\times_{\overline{H}^\C}\C^n$ and
the line bundle $L=\overline{E}\times_{\overline{H}^\C}\C$.

With these two bundles we have a $\overline{H}^\C$-equivariant map
$$Q:\overline{E}\times_{\overline{H}^\C}(\C^n\otimes\C^n)\longrightarrow\overline{E}\times_{\overline{H}^\C}\C$$
given fibrewise by $$v\otimes u\mapsto \sum v_iu_i=\langle
v,u\rangle$$ where $\overline{H}^\C$ acts on $\C^n\otimes\C^n$ by
$[(w,\lambda)]\mapsto\lambda w\otimes\lambda w$ and on $\C$ as above.
In other words $$Q:V\otimes V\longrightarrow L$$ is a
nowhere degenerate quadratic form on $V$ with values in $L$.

Since
$\mathfrak{gl}(n,\C)=\mathfrak{o}(n,\C)\oplus\overline{\liem}^\C$,
we have
$\overline{E}(\overline{\liem}^\C)=
\overline{E}\times_{\overline{H}^\C}\overline{\liem}^\C\subset\overline{E}
\times_{\overline{H}^\C}\mathfrak{gl}(n,\C)=\text{End}(V)$
and, indeed, $\overline{E}(\overline{\liem}^\C)=S^2_QV$.
Thus $\overline\Phi\in H^0(X,S^2_QV\otimes K)$ hence $\overline\Phi=\overline\Phi^*$ where
$\overline\Phi^*:V\to V\otimes K$ is such that $Q(\overline\Phi u,
v)=Q(u,\overline\Phi^*v)\in LK$. This means that the diagram
$$\xymatrix{V\ar[r]^(.4){q}\ar@<0.4ex>[d]_{\overline\Phi}&V^*\otimes L\ar[d]^{\overline\Phi^t\otimes 1_K\otimes 1_L}\\
    V\otimes K\ar@<-0.3ex>[r]^(.43){q\otimes 1_K}&V^*\otimes LK}$$
commutes where $q:V\to V^*\otimes L$ is the isomorphism
associated to $Q$, such that $q^t=q\otimes 1_L$.\end{proof}

The outcome of these results is that one can work with
$\EGL(n,\R)$-Higgs bundles instead of
$\PGL(n,\R)$-Higgs bundles with the advantage that in the former
case we work with the objects $(V,L,Q,\Phi)$, involving holomorphic vector bundles. That is what we
will do from now on.

We shall also call \emph{$\EGL(n,\R)$-Higgs bundles} the quadruples $(V,L,Q,\Phi)$ mentioned in the previous proposition.

Given an $\EGL(n,\R)$-Higgs bundle $(V,L,Q,\Phi)$, we associate a
$\PGL(n,\R)$-Higgs bundle $(E,\Phi_0)$, where $E$ is given by
the projection in sequence (\ref{seq}) and $\Phi_0$ is obtained by
projecting $\Phi$ to its traceless part.

\begin{proposition}\label{degl01}
Given a $\PGL(n,\R)$-Higgs bundle $(E,\Phi_0)$, it is possible to choose a lift of $(E,\Phi)$ to an $\EGL(n,\R)$-Higgs bundle $(V,L,Q,\Phi)$ such that $L$ is trivial or $\deg(L)=1$.
\end{proposition}
\begin{proof}
From (\ref{seq}) and from the actions (\ref{actCn}) and (\ref{actC}) defining
$V$ and $L$, two $\EGL(n,\R)$-Higgs bundles $(V,L,Q,\Phi)$ and $(V',L',Q',\Phi')$
give rise to the same $\PGL(n,\R)$-Higgs bundle if and only if
$V'=V\otimes F$ and $L'=L\otimes F^2$ where $F$ is a holomorphic
line bundle and $\Phi'_0=\Phi_0$.

Suppose we have an $\EGL(n,\R)$-Higgs bundle
$(V,L,Q,\Phi)$ associated to $(E,\Phi_0)$. Since $V$ and $V^*\otimes L$ are isomorphic, we have
$$\deg(V)=n\deg(L)/2.$$

If $\deg(L)$ is even we can choose a square
root $F$ of $L^{-1}$ and, from above, $(V\otimes F,\mathcal{O},Q',\Phi\otimes
1_F)$ also projects to $(E,\Phi_0)$.

 If $\deg(L)$ is odd then there is no such line bundle $F$.
Anyway, we can take $F$ such that $\deg(F)=(1-\deg(L))/2$ and
$(V\otimes F,L\otimes F^2, Q',\Phi\otimes 1_F)$ is also a lift of
$(E,\Phi_0)$ and the degree of the line bundle $L\otimes F^2$ is 1.
\end{proof}

From \cite{bradlow-garcia-prada-gothen:2004}, a $\GL(n,\R)$-Higgs bundle is a triple $(V,Q,\Phi)$ where $V$ is a rank $n$ holomorphic vector bundle, equipped with a nowhere degenerate quadratic form, and $\Phi$ is symmetric endomorphism of $V$.

\begin{corollary}\label{lift to GL(n,R)-Higgs}
Let $(E,\Phi_0)$ be a $\PGL(n,\R)$-Higgs bundle. Let $(V,L,Q,\Phi)$ be an $\EGL(n,\R)$-Higgs bundle which is a lift of $(E,\Phi_0)$. Then $(E,\Phi_0)$ lifts to a $\GL(n,\R)$-Higgs bundle if and only if $\deg(L)$ is even.
\end{corollary}
\begin{proof}
This follows directly from the proof of the above proposition: $\deg(L)$ is even if and only if we can change the lift $(V,L,Q,\Phi)$ to 
$(V\otimes F,\mathcal{O},Q',\Phi\otimes 1_F)$ (where $F^2=L^{-1}$) and this corresponds to a $\GL(n,\R)$-Higgs bundle.
\end{proof}

\begin{definition}\label{def of isom of quadruples}
Two $\EGL(n,\R)$-Higgs bundles $(V,L,Q,\Phi)$ and
$(V',L',Q',\Phi')$ are \emph{isomorphic} if there is a pair $(f,g)$ of
isomorphisms $f:V\to V'$ and $g:L\to L'$
such that the diagrams
$$\xymatrix{
    V\ar[r]^{f}\ar[d]_{\Phi}&V'\ar[d]^{\Phi'}\\
    V\otimes K\ar[r]^{f\otimes 1_K}&V'\otimes K} \hspace{1cm} \text{ and } \hspace{1cm}
\xymatrix{
    V\ar[r]^{f}\ar[d]_{q}&V'\ar[d]^{q'}\\
    V^*\otimes L\ar@<0.3ex>[r]^{(f^t)^{-1}\otimes g}&V'^*\otimes L'}$$
commute, where $q$ and $q'$ are the isomorphisms associated to $Q$
and $Q'$, respectively.
\end{definition}

Now we consider \emph{twisted orthogonal bundles}, i.e., triples $(V,L,Q)$ where $V$ is a holomorphic rank $n$ vector bundle equipped with a nowhere degenerate $L$-valued quadratic form $Q$. Of course, two twisted quadratic pairs $(V,L,Q)$ and $(V',L',Q')$ are isomorphic if there is a pair $(f,g)$ of isomorphisms $f:V\to V'$ and $g:L\to L'$ such that $((f^t)^{-1}\otimes g)q=q'f$.

Let $E$ and $E'$ be two principal $\EO(n,\C)$-bundles over $X$ and let $(V,L,Q)$ and $(V',L',Q')$ be the corresponding twisted orthogonal bundles through the actions (\ref{actCn}) and (\ref{actC}). It is easy to see that $E$ and $E'$ are isomorphic if and only if $(V,L,Q)$ and $(V',L',Q')$ are isomorphic.
Now, consider the notion of isomorphism between two $\EGL(n,\R)$-Higgs bundles, by applying Definition \ref{isomorphism of G-Higgs bundles}. Consider also Definition \ref{def of isom of quadruples} of isomorphism of $\EGL(n,\R)$-Higgs bundles written in terms of vector bundles, through Proposition \ref{pgl2}. We have then that, when applied to $\EGL(n,\R)$-Higgs bundles, the isomorphism notion of Definition \ref{isomorphism of G-Higgs bundles} is equivalent to the one of Definition \ref{def of isom of quadruples}, and that Proposition \ref{pgl2} gives a bijection between isomorphism classes of $\EGL(n,\R)$-Higgs bundles and of the objects $(V,L,Q,\Phi)$ which, because of this bijection, we also called $\EGL(n,\R)$-Higgs bundles.

Furthermore, the construction of Proposition \ref{pgl2} can be naturally applied to families of $\EGL(n,\R)$-Higgs bundles parametrized by varieties. Hence, the bijection between isomorphism classes of $\EGL(n,\R)$-Higgs bundles and of the objects $(V,L,Q,\Phi)$, naturally extends to families.

%%%%%%%%%%%%%%%%%%%%%%%%%%%%%%%%%%%%%%%%%%%%%%%%%%%%%%%%%%%%%%%%%%%%%%%%%%%%%
\section{Moduli space of $\EGL(n,\R)$-Higgs bundles}
%%%%%%%%%%%%%%%%%%%%%%%%%%%%%%%%%%%%%%%%%%%%%%%%%%%%%%%%%%%%%%%%%%%%%%%%%%%%%

%%%%%%%%%%%%%%%%%%%%%%%%%%%%%%%%
\subsection{Moduli space of $\EGL(n,\R)$-Higgs bundles}
%%%%%%%%%%%%%%%%%%%%%%%%%%%%%%%%

Recall from Subsection \ref{G-Higgs bundles} the definition of (poly,semi)stability of ($\GL(n,\C)$-)Higgs bundles, and let $$\M_{\GL(n,\C)}(d)$$ denote the moduli space of isomorphism
classes of polystable Higgs bundles of rank $n$ and degree $d$.
$\M_{\GL(n,\C)}(d)$ is \cite{nitsure:1991} a quasi-projective variety of complex
dimension $2n^2(g-1)+2$ which is smooth at the stable locus.

Given an $\EGL(n,\R)$-Higgs bundle $(V,L,Q,\Phi)$, we have a natural way
of associating to it a Higgs bundle $(V,\Phi)$ by simply forgetting the line bundle $L$ and the quadratic form $Q$.

\begin{definition}\label{poly EGL}
We say that an $\EGL(n,\R)$-Higgs bundle $(V,L,Q,\Phi)$ is \emph{polystable} if the corresponding
Higgs bundle $(V,\Phi)$ is polystable.
\end{definition}

The comparison of stability and strict polystability of $\EGL(n,\R)$-Higgs bundles and
that of the associated Higgs bundles is a delicate question, since the conditions might not correspond. Nevertheless, this problem does not occur for polystability due to the correspondence with the solutions to Hitchin's equations.

For a Lie group $G$, and a maximal compact subgroup $H$, the Hitchin's equations are equations for a pair $(A,\Phi)$ where $A$ is a $H$-connection on a $C^\infty$ $H^\C$-principal bundle $E_{H^\C}$ and $\Phi\in\Omega^{1,0}(X,E_{H^\C}(\liemc))$. The equations are
\begin{equation}\label{projectively flat Hitchin's equations}
 \begin{cases}
F(A)-[\Phi,\tau(\Phi)]=\lambda\,\omega\\
\bar\partial_A\Phi=0\\
 \end{cases}
\end{equation} where $F(A)\in\Omega^2(E_H,\lieh)$ is the curvature of $A$, $\tau$ is the involution on $G$ which defines $H$, $\lambda\in Z(\lieh)$ and $\omega$ is the normalized volume form on $X$ so that $\vol(X)=2\pi$. Furthermore, $\bar\partial_A$ is the unique $\bar\partial$ operator on $E_{H^\C}$, corresponding to the $H$-connection $A$ ($\bar\partial_A$ is then the unique holomorphic structure on $E_{H^\C}$ induced from the $(0,1)$-form $A^{0,1}$) and the second equation on (\ref{projectively flat Hitchin's equations}) says that $\Phi$ is holomorphic with respect to this holomorphic structure. The details of this theory can be found in \cite{hitchin:1987,garcia-prada-gothen-mundet:2008 II,bradlow-garcia-prada-mundet:2003}.

Now, given a $G$-Higgs bundle, $(E_{H^\C},\Phi)$ one associates to it a pair $(A,\Phi)$ given by a $H$-connection $A$ on the  $C^\infty$ $H^\C$-principal bundle $E_{H^\C}$ and $\Phi\in\Omega^{1,0}(X,E_{H^\C}(\liemc))$ (see, for instance, \cite{hitchin:1987,garcia-prada-gothen-mundet:2008 II}).
The Hitchin-Kobayashi correspondence says that a $G$-Higgs bundle $(E,\Phi)$ is polystable if and only if the associated pair $(A,\Phi)$ is a solution to the $G$-Hitchin's equations (see \cite{garcia-prada-gothen-mundet:2008 II} for the details).

Let us check that, the homomorphism
\begin{equation}\label{hom EGL to GL}
j:\EGL(n,\R)\longrightarrow\GL(n,\C)
\end{equation}
given by $j([(M,\alpha)])=\alpha M$ (which precisely corresponds to forgetting $L$ and $Q$ in $(V,L,Q,\Phi)$) induces a correspondence from solutions to the $\EGL(n,\R)$-Hitchin's equations to solutions to the $\GL(n,\C)$-Hitchin's equations.
Now, for $G=\EGL(n,\R)$ we have $H=\EO(n)$, hence $\lieh=\mathfrak{o}(n)\oplus\mathfrak{u}(1)$ (so  $Z(\lieh)=0\oplus\mathfrak{u}(1)$) and $\liemc=\mathfrak{o}(n,\C)^\perp\oplus 0\cong\mathfrak{o}(n,\C)^\perp$ (the space of symmetric matrices). For $G'=\GL(n,\C)$, we have $H'=\U(n)$ and $\lieh'=\mathfrak{u}(n)$ (hence $Z(\lieh)=\mathfrak{u}(1)$) and $(\liemc)'=\mathfrak{gl}(n,\C)$.
The homomorphism $j:G\to G'$, defined above, restricts to $j:H\to H'$ and therefore yields the map $j_*:\lieh\to\lieh'$ given by $j_*(M,\alpha)=\alpha+M$ and also $j_*:\liemc\to(\liemc)'$ given by $j_*(M,0)\mapsto M$. Hence, given a connection $(A,\beta)\in\Omega^1(E_H,\lieh)$ (where $A\in\Omega^1(E_H,\mathfrak{o}(n))$ and $\beta\in\Omega^1(E_H,\mathfrak{u}(1))$), we obtain the connection $\beta+A\in\Omega^1(E_{H'},\lieh')$ in $E_{H'}$.

\begin{proposition}
Let $(E_{\EO(n,\C)},\Phi)$ be an $\EGL(n,\R)$-Higgs bundle and $(E_{\GL(n,\C)},\Phi)$ the corresponding $\GL(n,\C)$-Higgs bundle obtained from $E_{\EO(n,\C)}$ by extending the structure group through the homomorphism $j$ defined in (\ref{hom EGL to GL}). Let $((A,\beta),\Phi)$ be the pair associated to $(E_{\EO(n,\C)},\Phi)$ and $(A',\Phi)$ be the pair associated to $(E_{\GL(n,\C)},\Phi)$. Then $((A,\beta),\Phi)$ is a solution to the $\EGL(n,\R)$-Hitchin's equations if and only if $(A',\Phi)$ is a solution to the $\GL(n,\C)$-Hitchin's equations.
\end{proposition}
\proof
The previous discussion shows that $j^*A'=\beta+A$ where $j$ is the homomorphism in (\ref{hom EGL to GL}).
Let $(M,\alpha)\in\Omega^2(E_H,\lieh)$ the curvature of $(A,\beta)$. Then, from (\ref{projectively flat Hitchin's equations}), the $\EGL(n,\R)$-Hitchin's equations are
$$\begin{cases}
(M,\alpha)+[\Phi,\Phi]=(0,\lambda)\,\omega\\
\bar\partial_{(A,\beta)}\Phi=0\\
 \end{cases}$$
where $\lambda\in\mathfrak{u}(1)$. Notice that, in this case, $\tau(X)=-X^t$. Moreover $[\Phi,\Phi]\in\Omega^2(E_H,\mathfrak{o}(n))$, and, from the definition of $\liemc$, we have $\bar\partial_{(A,\beta)}\Phi=\bar\partial_A\Phi$. Hence the above equations are indeed
$$\begin{cases}
M+[\Phi,\Phi]=0\\
\alpha=\lambda\,\omega\\
\bar\partial_A\Phi=0\\
 \end{cases}.$$
From this we obtain $$\begin{cases}
\alpha+M+[\Phi,\Phi]=\lambda\,\omega\\
\bar\partial_A\Phi=0\\
 \end{cases}$$ which are the $\GL(n,\C)$-Hitchin's equations. Furthermore, $\alpha+M$ is the curvature of $\beta+A$, and notice again that $\bar\partial_{\beta+A}\Phi=\bar\partial_A\Phi$ due to the definition of $\liemc$ and to the map $\liemc\to(\liemc)'$ defined above. So, the equations in $E_{\EO(n,\C)}$ and in the $\GL(n,\C)$ bundle $E_{\GL(n,\C)}$ obtained from $E_{\EO(n,\C)}$ by extending the structure group through the homomorphism $j$, are equivalent, and this proves the result.
\endproof

 Hence, the previous proposition and the Hitchin-Kobayashi correspondence show that Definition \ref{poly EGL} is consistent with the notion of polystability for $\EGL(n,\R)$-Higgs bundles.

\begin{notation}
Write $$\M_d$$ for the set of isomorphism classes of
polystable $\EGL(n,\R)$-Higgs bundles $(V,L,Q,\Phi)$ with $\mathrm{rk}(V)=n$ and
$\deg(L)=d$ (hence $\deg(V)=nd/2$).
\end{notation}

The group $\EGL(n,\R)$ can be seen as a closed subgroup of $\GL(n,\C)\times\U(1)$
through
$$[(A,\lambda)]\mapsto(\lambda A,\lambda^2).$$
The moduli space of $\GL(n,\C)\times\U(1)$-Higgs bundles is
$\M_{\GL(n,\C)}(d_1)\times\mathrm{Jac}^{d_2}(X)$, where
$\mathrm{Jac}^{d_2}(X)$ is the subspace of the Picard group of the
compact Riemann surface $X$ which parametrizes holomorphic line
bundles of degree $d_2$. It is isomorphic to the Jacobian of $X$.

We shall realize $\M_d$ as a subspace of
$\M_{\GL(n,\C)}(nd/2)\times\Jac^d(X)$ and it will be in this
subspace that we will work.

\begin{lemma}\label{lema} The map $$\begin{array}{rccl}
  i:&\M_d & \longrightarrow & \M_{\GL(n,\C)}(nd/2)\times\mathrm{Jac}^d(X) \\
  &(V,L,Q,\Phi) & \longmapsto & ((V,\Phi),L).
\end{array}$$ is injective.
\end{lemma}
\begin{proof} First of all we see that we only have to take care with the
form $Q$. Indeed, if $i(V,L,Q,\Phi)=i(V',L',Q',\Phi')$, then there
are isomorphisms $f:V\to V'$ such that $\Phi'f=(f\otimes
1_K)\Phi$ and $g:L\to L'$. Therefore $(f,g)$ is an
isomorphism between $(V,L,Q'',\Phi)$ and $(V',L',Q',\Phi')$ where
$Q''$ is given by $q''=(f^t\otimes g^{-1})q'f$.

Consider then the $\EGL(n,\R)$-Higgs bundles $(V,L,Q,\Phi)$ and $(V,L,Q',\Phi)$. These are mapped to $((V,\Phi),L)$
and we have to see that $(V,L,Q,\Phi)$ and $(V,L,Q',\Phi)$ are isomorphic.

Suppose that $(V,\Phi)$ is stable. The automorphism
$(q')^{-1}q$ of $V$ is $\Phi$-equivariant hence, from
stability, $(q')^{-1}q=\lambda\in\C^*$, so
$(\sqrt{\lambda},1_L)$ is an isomorphism between $(V,L,Q,\Phi)$
and $(V,L,Q',\Phi)$.

Suppose now that $(V,\Phi)$ is strictly polystable, with
$$(V,\Phi)=(V_1,\Phi_1)\oplus\dots\oplus(V_k,\Phi_k).$$ Here
$\Phi_i:V_i\to V_i\otimes K$ and all the Higgs bundles
$(V_i,\Phi_i)$ are stable and with the same slope $\mu=\deg(L)/2$.

Consider the decomposition of $q:V\to V^*\otimes L$ as a
matrix $(q_{ij})$ compatible with that of $V$ and of
$$V^*\otimes L=V_1^*\otimes L\oplus\dots\oplus V_k^*\otimes L.$$ Hence
$$q_{ij}\in\Hom(V_j,V_i^*\otimes L)$$ (note that
$q|_{V_j}=q_{1j}\oplus\dots\oplus q_{kj}$). Suppose that $q_{ij}$
is non-zero, for some $i,j$. Since $$(\Phi_i^t\otimes 1_K\otimes
1_L)q_{ij}=(q_{ij}\otimes 1_K)\Phi_j$$ then $q_{ij}$ is a
homomorphism between $(V_j,\Phi_j)$ and $(V_i^*\otimes
L,\Phi_i^t\otimes 1_K\otimes 1_L)$. These are stable Higgs bundles
and $\mu(V_j)=\mu(V_i^*\otimes L)$, therefore $q_{ij}$ must be an
isomorphism. Hence for each pair $i,j$, if $q_{ij}$ is non-zero,
then it is an isomorphism.

We will consider now three cases.

In the first case we suppose that $(V,\Phi)$ is a direct sum of
isomorphic copies of the same Higgs bundle $(W,\Phi_W)$, with
$\Phi_W:W\to W\otimes K$. Let then
$$(V,\Phi)=\underbrace{(W,\Phi_W)\oplus\dots\oplus(W,\Phi_W)}_{k \text{ summands}}.$$
We have $(W,\Phi_W)\cong(W^*\otimes L,\Phi_W^t\otimes 1_K\otimes
1_L)$ but we can have more than one isomorphism on each column of
$(q_{ij})$.

Choose $i_0$ and $j_0$ such that $q_{i_0j_0}:(W,\Phi_W)\to(W^*\otimes
L,\Phi_W^t\otimes 1_K\otimes 1_L)$ is non-zero, being therefore an
isomorphism. If $q_{ij}:(W,\Phi_W)\to(W^*\otimes
L,\Phi_W^t\otimes 1_K\otimes 1_L)$ is any homomorphism, then
$(q_{i_0j_0}^{-1})q_{ij}$ is an endomorphism of $(W,\Phi_W)$ and,
since $(W,\Phi_W)$ is stable, $q_{ij}=\alpha_{ij}q_{i_0j_0}$
where $\alpha_{ij}\in\C$ and, moreover, $\alpha_{ij}=0$ if and only if
$q_{ij}=0$. Hence the choice of $q_{i_0j_0}$ gives a way to represent  $(q_{ij})$ by a symmetric
$k\times k$ matrix where each entry is $\alpha_{ij}$ and which
can be diagonalized through a $k\times k$ matrix $(\lambda_{ij})$.
Define the automorphism $g$ of $V$ given, with respect to the
decomposition of $V$, by a $k\times k$ matrix $(g_{ij})$ where $g_{ij}=\lambda_{ij}:W\to W$. Thus $g$ is such that
$$(g^t\otimes 1_L)qg=\tilde{q}$$ where $\tilde{q}:V\to V^*\otimes L$ is an
isomorphism which is diagonal, by
$\mathrm{rk}(W)\times\mathrm{rk}(W)$ blocks, with respect to the
given decomposition of $V$.

Note also that $g$ is $\Phi$-equivariant. Hence, if $\tilde{Q}$ is the
quadratic form associated to $\tilde{q}$, $(g,1_L):(V,L,\tilde{Q},\Phi)\to(V,L,Q,\Phi)$ is an
isomorphism.
So we can suppose that $(q_{ij})$ and $(q'_{ij})$ are diagonal and,
an argument analogous to the case where $(V,\Phi)$ was stable shows then that
$(V,L,Q,\Phi)$ and $(V,L,Q',\Phi)$ are isomorphic, the isomorphism being $(f,1_L)$ where $f$ is given, according to
the decomposition of $V$, by
$$f=\begin{pmatrix}
    \sqrt{\lambda_1} & 0 & 0 & \dots & 0 \\
    0 &  \sqrt{\lambda_2} & 0 & \dots & 0 \\
    0 & 0 &  \ddots &  \dots & 0\\
    \vdots & \vdots & \vdots & \ddots & \vdots \\
    0 & 0 & 0 & \dots & \sqrt{\lambda_k}\
  \end{pmatrix}.$$
Here, each $\lambda_i\in\C^*$ is such that $q_{ii}=\lambda_i q'_{ii}$.

In the second case we consider
$$V=\underbrace{(W,\Phi_W)\oplus\dots\oplus
(W,\Phi_W)}_{l\text{ summands}}\oplus\underbrace{(W^*\otimes
L,\Phi_W^t\otimes 1_K\otimes 1_L)\oplus\dots\oplus (W^*\otimes
L,\Phi_W^t\otimes 1_K\otimes 1_L)}_{l\text{ summands}}$$ with
$W\not\cong W^*\otimes L$. In this case $q_{ij}=0$ if $-l\leq
i-j\leq l$. So $(q_{ij})$ splits into four $l\times l$ blocks in
the following way
$$(q_{ij})=\begin{pmatrix}
 0 & q_2 \\
 q_1 & 0
\end{pmatrix}$$
where $q_1$ represents $$q|_{W\oplus\dots\oplus
W}:W\oplus\dots\oplus W\longrightarrow W\oplus\dots\oplus W$$ and
$q_2$ represents $$q|_{W^*\otimes L\oplus\dots\oplus W^*\otimes
L}:W^*\otimes L\oplus\dots\oplus W^*\otimes L \longrightarrow
W^*\otimes L\oplus\dots\oplus W^*\otimes L.$$ Again, using the
stability of $(W,\Phi_W)$ and of $(W^*\otimes L,\Phi^t_W\otimes
1_K\otimes 1_L)$ and the fact that $\Phi$ is symmetric with
respect to $q$, we see that each entry of $q_1$ is given by a
scalar. The same happens with $q_2$. Hence we can write
$$(q_{ij})=\begin{pmatrix}
 0 & A^t \\
 A & 0
\end{pmatrix}$$
where $A$ is a non-singular $l\times l$ matrix. Now, if we write in
the same way,
$$(q'_{ij})=\begin{pmatrix}
 0 & B^t \\
 B & 0
\end{pmatrix}$$
then consider the isomorphism of $V$ given by $$f=\begin{pmatrix}
 B^{-1}A & 0 \\
 0 & I_l
\end{pmatrix}$$
where we mean by this that each entry of $B^{-1}A$ represents a
scalar automorphism of $(W,\Phi_W)$ and $f$ is the identity over $W^*\otimes
L\oplus\dots\oplus W^*\otimes L$. So $(f,1_L)$ is an
isomorphism between $(V,L,Q,\Phi)$ and $(V,L,Q',\Phi)$.

 The last case is the generic one, where we consider a combination of the previous cases.
We can always write
\begin{equation}\label{generic decomposition}
V=(V_1\oplus\dots\oplus V_i)\oplus(V_{i+1}\oplus\dots\oplus
V_j)\oplus\dots\oplus(V_{l+1}\oplus\dots\oplus V_k)
\end{equation}
and
$$\Phi=(\Phi_1\oplus\dots\oplus\Phi_i)\oplus(\Phi_{i+1}\oplus\dots\oplus
\Phi_j)\oplus\dots\oplus(\Phi_{l+1}\oplus\dots\oplus\Phi_k)$$ where
$(V_a,\Phi_a)$ and $(V_b,\Phi_b)$ are inside the same parenthesis in
(\ref{generic decomposition}) if and only if are isomorphic or
$(V_b,\Phi_b)\cong (V_a^*\otimes L,\Phi_a^t\otimes 1_K\otimes 1_L)$.
If $V_a$ and $V_b$ are not
inside the same parenthesis in (\ref{generic decomposition}) then
$q_{ab}=0=q'_{ab}$. Hence we have an isomorphism $f$ between $(V,L,Q,\Phi)$
and $(V,L,Q',\Phi)$ where $f$ is diagonal by blocks (not all of the
same size), each corresponding to an isomorphism of one of the
previous cases.
\end{proof}

We identify $\M_d$ with its image by the map $i$ and therefore
consider $\M_d$ as a subspace of
$\M_{\GL(n,\C)}(nd/2)\times\mathrm{Jac}^d(X)$.

Note that $\EGL(n,\R)$ is a reductive (not
semisimple) Lie group. Therefore, in view of Theorem \ref{fundamental correspondence for semisimple G}, $\M(c)$
is homeomorphic to $\cR_{\Gamma,\EGL(n,\R)}(c)$, for each topological class $c$. 

\begin{notation}\label{Rgamma}
From now on, we
shall write $\cR$ instead of
$\cR_{\Gamma,\EGL(n,\R)}$. 
\end{notation}

\begin{remark} Let $n\geq 4$ be even. Notice that a representation
$$\rho\in\Hom_{\rho(J)\in (Z(\PGL(n,\R))\cap\PO(n))_0}^\text{red}(\Gamma,\PGL(n,\R))$$ (cf. (\ref{homgammaG})), which is the same as
$\rho\in\Hom^\text{red}(\pi_1X,\PGL(n,\R))$, does not always lift to a
representation $\rho'\in\Hom^\text{red}(\Gamma,\GL(n,\R))$ with the
condition $\rho'(J)\in (Z(\GL(n,\R))\cap\Or(n))_0=I_n$, (i.e., to a representation $\rho'\in\Hom^\text{red}(\pi_1X,\GL(n,\R))$) as we have seen in Proposition \ref{nonempty}. But it always lifts to a representation
of $\Gamma$ in $\EGL(n,\R)$ such that $\rho'(J)\in (Z(\EGL(n,\R))\cap\EO(n))_0\cong\U(1)$, since
$Z(\EGL(n,\R))\cap\EO(n)$ is connected. This is an instance of the fact that a $\PGL(n,\R)$-Higgs bundle lifts to an $\EGL(n,\R)$-Higgs bundle, but not to a $\GL(n,\R)$-Higgs bundle.
\end{remark}

The following lemma will be useful in Section \ref{morse quadruples}.

\begin{lemma}\label{fechado}
$\M_d$ is a closed subspace of
$\M_{\GL(n,\C)}(nd/2)\times\mathrm{Jac}^d(X)$.
\end{lemma}
\begin{proof} Let $\cR_d$ be the subspace of $\cR$ which
corresponds to $\M_d$, under Theorem \ref{fundamental correspondence for semisimple G}.

Consider the following commutative diagram
$$\xymatrix{\cR_{\Gamma,\GL(n,\C)\times\U(1)}(nd/2,d)\ar[r]^(.45){\simeq}&\M_{\GL(n,\C)}(nd/2)\times\Jac^d(X)\\
\cR_d\ar[u]^{j}\ar[r]^(.5){\simeq}&\M_d\ar[u]_{i}}$$ where
$\cR_{\Gamma,\GL(n,\C)\times\U(1)}(nd/2,d)$ is the subspace
of $\cR_{\Gamma,\GL(n,\C)\times\U(1)}$ which corresponds, via again Theorem \ref{fundamental correspondence for semisimple G},
to $\M_{\GL(n,\C)}(nd/2)\times\Jac^d(X)$.
So, the top and bottom maps are the homeomorphisms given by Theorem
\ref{fundamental correspondence for semisimple G}.
The map $j$ is induced by the injective map $[(A,\lambda)]\mapsto
(\lambda A,\lambda^2)$ of $\EGL(n,\R)$ into
$\GL(n,\C)\times\U(1)$. The diagram commutes due to the fact that the actions
of $\EO(n,\C)$ on $\C^n$ and on $\C$
defining an $\EGL(n,\R)$-Higgs bundle are given by $[(w,\lambda)]\mapsto\lambda w$ and
$[(w,\lambda)]\mapsto\lambda^2$,
which are then compatible with the inclusion of
$\EGL(n,\R)$ into $\GL(n,\C)\times\U(1)$ and
therefore with the induced action of $\EGL(n,\R)$ on $\C^n\times\C$. Now, since, from Lemma \ref{lema}, $i$ is injective, it follows that $j$ is also injective, and, as we identify
$\M_d$ with $i(\M_d)$, we also identify $\cR_d$ with
$j(\cR_d)$. So $\cR_d$ can be seen as the space of
reductive homomorphisms of $\Gamma$ in $\GL(n,\C)\times\U(1)$ which
have their image in $\EGL(n,\R)$, modulo
$\GL(n,\C)\times\U(1)$-equivalence.

Now, since $\EGL(n,\R)$ is a closed subgroup
of $\GL(n,\C)\times\U(1)$, it follows that $\cR_d$ is closed in
$\cR_{\Gamma,\GL(n,\C)\times\U(1)}(nd/2,d)$, hence
$\M_d$ is closed in $\M_{\GL(n,\C)}(nd/2)\times\Jac^d(X)$.
\end{proof} 

%%%%%%%%%%%%%%%%%%%%%%%%%%%%%%%%%%%%%%
\subsection{Deformation theory of $\EGL(n,\R)$-Higgs bundles}
%%%%%%%%%%%%%%%%%%%%%%%%%%%%%%%%%%%%%%

In this section, we briefly recall the description of Biswas and Ramanan \cite{biswas-ramanan:1994} (see also \cite{nitsure:1991}) of the deformation theory of $G$-Higgs bundles and, in particular, the identification of the tangent space of $\M_G$ at the smooth points, and then apply it to the case of $\EGL(n,\R)$-Higgs bundles.

The spaces $\liehc$ and $\liemc$ in the Cartan decomposition of $\liegc$ verify the relation $$[\liehc,\liemc]\subset\liemc$$ hence, given $v\in\liemc$, there is
an induced map $\text{ad}(v)|_{\liehc}:\liehc\to\liemc$.
Applying this to a $G$-Higgs bundle $(E_{H^\C},\Phi)$, we obtain the
following complex of sheaves on $X$:
$$C^\bullet(E_{H^\C},\Phi):\mathcal{O}(E_{H^\C}(\liehc))\xrightarrow{\text{ad}(\Phi)}\mathcal{O}(E_{H^\C}(\liemc)\otimes K).$$

\begin{proposition}[\textbf{Biswas, Ramanan} \cite{biswas-ramanan:1994}]\label{deformation for Higgs}
Let $(E_{H^\C},\Phi)$ represent a $G$-Higgs bundle over the compact Riemann surface $X$.
\begin{enumerate}
    \item The
infinitesimal deformation space of $(E_{H^\C},\Phi)$ is isomorphic to the
first hypercohomology group $\mathbb{H}^1(X,C^\bullet(E_{H^\C},\Phi))$ of the
complex $C^\bullet(E_{H^\C},\Phi)$;
    \item There is a long exact sequence
\begin{equation*}
\begin{split}
0&\longrightarrow\mathbb{H}^0(X,C^\bullet(E_{H^\C},\Phi))\longrightarrow
H^0(X,E_{H^\C}(\liehc))\longrightarrow
H^0(X,E_{H^\C}(\liemc)\otimes
K)\longrightarrow\\
&\longrightarrow\mathbb{H}^1(X,C^\bullet(E_{H^\C},\Phi))\longrightarrow H^1(X,E_{H^\C}(\liehc))\longrightarrow
H^1(X,E_{H^\C}(\liemc)\otimes K)\longrightarrow\\
&\longrightarrow\mathbb{H}^2(X,C^\bullet(E_{H^\C},\Phi))\longrightarrow 0
\end{split}
\end{equation*}
where the maps $H^i(X,E_{H^\C}(\liehc))\to H^i(X,E_{H^\C}(\liemc)\otimes K)$ are induced by
$\mathrm{ad}(\Phi)$.
\end{enumerate}
\end{proposition}

Proposition \ref{deformation for Higgs} applied to the case of $\EGL(n,\R)$-Higgs bundles, yields:
\begin{proposition}\label{deformation complex for quadruples}
Let $(V,L,Q,\Phi)$ be an $\EGL(n,\R)$-Higgs bundle over $X$. There is a complex of
sheaves
$$C^\bullet(V,L,Q,\Phi):\Lambda^2_QV\oplus\mathcal{O}\xrightarrow{[\Phi,-]}S^2_QV\otimes
K$$ and
\begin{enumerate}
    \item the
infinitesimal deformation space of $(V,L,Q,\Phi)$ is isomorphic to
the first hypercohomology group
$\mathbb{H}^1(X,C^\bullet(V,L,Q,\Phi))$ of
$C^\bullet(V,L,Q,\Phi)$. In particular, if $(V,L,Q,\Phi)$
represents a smooth point of $\M_d$, then
    $$T_{(V,L,Q,\Phi)}\M\simeq\mathbb{H}^1(X,C^\bullet(V,L,Q,\Phi));$$
    \item there is an exact sequence
\begin{equation*}
\begin{split}
0&\longrightarrow\mathbb{H}^0(X,C^\bullet(V,L,Q,\Phi))\longrightarrow
H^0(X,\Lambda^2_QV\oplus\mathcal{O})\longrightarrow H^0(X,S^2_QV\otimes
K)\longrightarrow\\
&\longrightarrow\mathbb{H}^1(X,C^\bullet(V,L,Q,\Phi))\longrightarrow H^1(X,\Lambda^2_QV\oplus\mathcal{O})\longrightarrow
H^1(X,S^2_QV\otimes K)\longrightarrow\\
&\longrightarrow\mathbb{H}^2(X,C^\bullet(V,L,Q,\Phi))\longrightarrow 0
\end{split}
\end{equation*}
where
the maps $H^i(X,\Lambda^2_QV\oplus\mathcal{O})\to
H^i(X,S^2_QV\otimes K)$ are induced by the map $[\Phi,-]$.
\end{enumerate}
\end{proposition}

\subsection{Topological classification of $\EGL(n,\R)$-Higgs bundles}
Our calculations will be
performed on $\M_d$ so we will also need the topological invariants
of $\EGL(n,\R)$-Higgs bundles.

We will define these discrete invariants using the relation between $\EGL(n,\R)$-Higgs bundles and $\PGL(n,\R)$-Higgs bundles. In fact, we already know from Theorem \ref{classtoppgl} that the invariants $\mu_1$ and $\mu_2$ completely classify $\PGL(n,\R)$-bundles over $X$, and also know from Proposition \ref{degl01} that if two $\EGL(n,\R)$-Higgs bundles project to the same $\PGL(n,\R)$-Higgs bundle, then the degree
of the corresponding line bundle $L$ is equal modulo $2$.

As we are dealing with topological classification of bundles, we can forget the Higgs field and consider only twisted orthogonal bundles $(V,L,Q)$ which correspond to elements of the set $H^1(X,\mathcal{C}(\EO(n)))$, and $\PO(n)$-bundles $E$ which are in bijection with $H^1(X,\mathcal{C}(\PO(n)))$. There is then a relation between $(V,L,Q)$ and $E$ which is similar to the one between $\EGL(n,\R)$-Higgs bundles and $\PGL(n,\R)$-Higgs bundles, but now forgetting the Higgs field.

Consider then the following commutative diagram: %\vspace{0,5cm}
\begin{equation}\label{diagram}
\xymatrix{&H^1(X,\Z_2)\ar[r]^{\cong}\ar[d]&H^1(X,\Z_2)\ar[d]&\\
          H^1(X,\mathcal{C}(\U(1)))\ar[r]\ar[d]_{\cong}&H^1(X,\mathcal{C}(\mathrm{O}(n)\times\U(1)))\ar[r]\ar[d]^{p_2'}&H^1(X,\mathcal{C}(\mathrm{O}(n)))\ar[r]\ar[d]^{p_2}&0\\
          H^1(X,\mathcal{C}(\U(1)))\ar[r]&H^1(X,\mathcal{C}(\EO(n)))\ar[r]\ar[d]&H^1(X,\mathcal{C}(\PO(n)))\ar[r]\ar[d]&0\\
          &H^2(X,\Z_2)\ar[r]^{\cong}&H^2(X,\Z_2).}
\end{equation}
%\vspace{0,5cm}
The map $$p_2':H^1(X,\mathcal{C}(\mathrm{O}(n)\times\U(1)))\longrightarrow H^1(X,\mathcal{C}(\EO(n)))$$ is induced from the projection $\Or(n)\times\U(1)\to\EO(n)$ and hence defined, in terms of vector bundles, by $$p_2'((W,Q_W),M)=(W\otimes M,M^2,Q_W\otimes 1_{M^2}).$$ Once again we see that $(V,L,Q)$ is in the image of $p_2'$ if and only if $\deg(L)$ is even. Moreover, if this is the case, the pre-image of $(V,L,Q)$ under $p_2'$ is the following set of $2^{2g}$ pairs:
\begin{equation}\label{pre image of p2'}
p_2'^{-1}(V,L,Q)=\{((V\otimes L^{-1/2}F,Q\otimes Q_F\otimes 1_{L^{-1}}),L^{1/2}F)\st F^2\cong\mathcal O\}
\end{equation}
where $L^{1/2}$ is a fixed square root of $L$ (notice that saying that $F$ is a $2$-torsion point of the Jacobian is equivalent to say that $F$ has a nowhere degenerate quadratic form $Q_F:F\otimes F\to\mathcal O$).

The map $p_2$ is the one induced from the projection $\Or(n)\to\PO(n)$.

\begin{definition}\label{overlinemu1}
Given a twisted orthogonal bundle $(V,L,Q)$, let $E$ be the corresponding $\PO(n)$-bundle under the map $H^1(X,\mathcal{C}(\EO(n)))\to H^1(X,\mathcal{C}(\PO(n)))$. Define the \emph{first invariant} $\overline\mu_1$ of $(V,L,Q)$ as $$\overline{\mu}_1(V,L,Q)=\mu_1(E)\in(\Z_2)^{2g}.$$
\end{definition}

$\overline{\mu}_1(V,L,Q)$ is then the obstruction to
reducing the structure group of $(V,L,Q)$ to $\ESO(n)=\SO(n)\times_{\Z_2}\U(1)$. Hence, this happens if and only if
$E$ reduces to a $\mathrm{PSO}(n)$-bundle.

If $\deg(L)$ is even, then
\begin{equation}\label{linemu1=w1 if degL even}
\overline{\mu}_1(V,L,Q)=w_1(V\otimes L^{-1/2},Q\otimes 1_{L^{-1}})
\end{equation}
the first Stiefel-Whitney class of the real orthogonal bundle $V\otimes L^{-1/2}$ (the value of $w_1$ is independent of the choice of the square root of $L$ because $n$ is even - notice that $w_1(W\otimes F)=w_1(W)+\rk(W)_2w_1(F)$ where $\rk(W)_2=\rk(W)\text{ mod }2$).

\vspace{0,5cm}
\begin{definition}\label{overlinemu2}
Let $(V,L,Q)$ be a twisted orthogonal bundle with $\rk(V)=n\geq 4$. Define the \emph{second invariant} $\overline\mu_2$ of $(V,L,Q)$ as follows:
\begin{enumerate}
 \item If $\overline\mu_1(V,L,Q)=0$, $$\overline{\mu}_2(V,L,Q)=\begin{cases}
    (w_2(V\otimes L^{-1/2}),\deg(L))\in\Z_2\times 2\Z & \text{if }\deg(L)\text{ even}\\
    \deg(L)\in 2\Z+1 & \text{if }\deg(L)\text{ odd}
\end{cases}$$ where $w_2(V\otimes L^{-1/2})$ is the second Stiefel-Whitney class of $V\otimes L^{-1/2}$ and $2\Z$ represents the set of even integers and $2\Z+1$ the set of odd integers.
 \item If $\overline\mu_1(V,L,Q)\neq 0$, $$\overline\mu_2(V,L,Q)=\deg(L)\in\Z.$$
\end{enumerate}
\end{definition}
Notice that on the first and third items, $w_2(V\otimes L^{-1/2})$ does not depend on the choice of the square root of $L$ due to the vanishing of $\overline{\mu}_1(V,L,Q)$ (cf. Remark \ref{SWclasses}).

Let $E$ be a $\PO(n)$-bundle and $(V,L,Q)$ be a twisted orthogonal bundle which maps to $E$. From (\ref{diagram}), $E$ lifts to a $\mathrm{O}(n)$-bundle if and only if
$(V,L,Q)$ lifts to a $\mathrm{O}(n)\times\U(1)$-bundle and this occurs if and only if $\deg(L)$ is even (recall also Corollary \ref{lift to GL(n,R)-Higgs}). Hence, using Proposition \ref{obs} and the fact that by definition $\overline\mu_1(V,L,Q)=\mu_1(E)$, the following proposition is immediate:
\begin{proposition}\label{relation overlinemu2 and mu2}
Let $n\geq 4$ be even. Let $E$ be a $\PO(n)$-bundle and $(V,L,Q)$ be a twisted orthogonal bundle which maps to $E$.
\begin{enumerate}
 \item If $\overline\mu_1(V,L,Q)=0$, then:
\begin{itemize}
 \item $\overline\mu_2(V,L,Q)=(0,\deg(L))\text{ with }\deg(L)\text{ even }\Longleftrightarrow\mu_2(E)=0$;
 \item $\overline\mu_2(V,L,Q)=(1,\deg(L))\text{ with }\deg(L)\text{ even }\Longleftrightarrow\mu_2(E)=1$;
 \item $\overline\mu_2(V,L,Q)=\deg(L)\text{ with }\deg(L)\text{ odd }\Longleftrightarrow\mu_2(E)=\omega_n$.
\end{itemize}
 \item If $\overline\mu_1(V,L,Q)\neq 0$, then:
\begin{itemize}
 \item $\overline\mu_2(V,L,Q)=\deg(L)\text{ even }\Longleftrightarrow\mu_2(E)=0$;
 \item $\overline\mu_2(V,L,Q)=\deg(L)\text{ odd }\Longleftrightarrow\mu_2(E)=\omega_n$.
\end{itemize}
\end{enumerate}
\end{proposition}

From Theorem \ref{classtoppgl}, we know that $\mu_1$ and $\mu_2$ completely classify $\PO(n)$-bundles. Moreover, since we also know that the difference between two $(V,L,Q)$ and $(V',L',Q')$ mapping to the same $\PO(n)$-bundle lies in the degree of $L$, we have then the following:
\begin{theorem}\label{topquad}
Let $X$ be a closed oriented surface of genus $g\geq 2$ and let $n\geq 4$ be even. Then twisted orthogonal bundles over $X$ are  topologically classified by the invariants
$$(\overline{\mu}_1,\overline{\mu}_2)\in\left(\{0\}\times\left(\left(\Z_2\times2\Z\right)\cup(2\Z+1)\right)\right)\cup\left(\left(\Z_2^{2g}\setminus\{0\}\right)\times\Z\right).$$
\end{theorem}

Now, returning to our principal objects - $\EGL(n,\R)$-Higgs bundles - we see that Theorem \ref{topquad} also gives the topological classification of $\EGL(n,\R)$-Higgs bundles.

\begin{notation}\label{submodquad}
Let $$\M(\overline{\mu}_1,\overline{\mu}_2)$$
denote the subspace of the space of $\EGL(n,\R)$-Higgs bundles in
which the $\EGL(n,\R)$-Higgs bundles have invariants
$(\overline{\mu}_1,\overline{\mu}_2)$.
\end{notation}

\begin{remark}\label{identif}
If $\overline{\mu}_1=0$ and if $d_1$ and $d_2$ are even, then
$$\M(0,(w_2,d_1))\simeq\M(0,(w_2,d_2))$$
and, if $d_1$ and $d_2$ are odd,
$$\M(0,d_1)\simeq\M(0,d_2).$$
If $\overline{\mu}_1\neq 0$ and $d_1=d_2\text{ mod }2$, then $$\M(\overline{\mu}_1,d_1)\simeq\M(\overline{\mu}_1,d_2).$$
In all cases the bijection is given by $(V,L,Q,\Phi)\mapsto
(V\otimes F,L\otimes F^2,Q\otimes 1_{F^2},\Phi\otimes 1_F)$, where $F$ is
a holomorphic line bundle of suitable degree.
\end{remark}

Again, we define the same invariants for the space of representations
$\cR$ (recall Notation \ref{Rgamma}). $\cR(\overline{\mu}_1,(w_2,d))$ corresponds
to $\M(\overline{\mu}_1,(w_2,d))$ if $d$ is even, and
$\cR(\overline{\mu}_1,d)$ corresponds to
$\M(\overline{\mu}_1,d)$ if $d$ is odd.

From Proposition \ref{pgl1}, the surjective map taking an $\EGL(n,\R)$-Higgs bundle to the
corresponding $\PGL(n,\R)$-Higgs bundle induces a surjective continuous map
$p:\cR\to \cR_{\PGL(n,\R)}$. Using Propositions \ref{degl01} and \ref{relation overlinemu2 and mu2}, the following
is immediate.
\begin{proposition}\label{projrep}
The map $p:\cR\to \cR_{\PGL(n,\R)}$ satisfies the following identities:
\begin{enumerate}
 \item If $\overline{\mu}_1=0$, then $$p(\cR(0,(0,0)))=\cR_{\PGL(n,\R)}(0,0)$$
\vspace{-0,3cm}$$p(\cR(0,(1,0)))=\cR_{\PGL(n,\R)}(0,1)$$ and
$$p(\cR(0,1))=\cR_{\PGL(n,\R)}(0,\omega_n).$$
 \item If $\overline{\mu}_1\neq 0$, then $$p(\cR(\overline{\mu}_1,0))=\cR_{\PGL(n,\R)}(\mu_1,0)$$ and
$$p(\cR(\overline{\mu}_1,1))=\cR_{\PGL(n,\R)}(\mu_1,\omega_n).$$
\end{enumerate}
\end{proposition}

 From this and from Proposition \ref{nonempty}, we have:
\begin{corollary}\label{nonempty cR}
$\cR(\overline{\mu}_1,\overline{\mu}_2)$ is
non-empty for any choice of invariants $\overline{\mu}_1$ and $\overline{\mu}_2$.
\end{corollary}

%%%%%%%%%%%%%%%%%%%%%%%%%%%%%%%%%%%%%%%%%%%%%%%%%%%%%%%%%%%%%%%%%%%%%%%%%%%%%%%%%%%%%%%%%%%%%%%%%
\section{The Hitchin proper functional}\label{morse quadruples}
%%%%%%%%%%%%%%%%%%%%%%%%%%%%%%%%%%%%%%%%%%%%%%%%%%%%%%%%%%%%%%%%%%%%%%%%%%%%%%%%%%%%%%%%%%%%%%%%%

Here we use the method introduced by Hitchin in \cite{hitchin:1987} to
study the topology of moduli space $\M_G$ of $G$-Higgs bundles.

Define $$f:\M_G(c)\longrightarrow\R$$ by
\begin{equation}\label{proper function}
f(E_{H^\C},\Phi)=\|\Phi\|_{L^2}^2=\int_{X}|\Phi|^2\mathrm{dvol}.
\end{equation}
This function $f$ is usually called the \emph{Hitchin functional}.

Here we are using the \emph{harmonic metric} (cf. \cite{corlette:1988,donaldson:1987}) on $E_{H^\C}$ to define
$\|\Phi\|_{L^2}$. So we are using the identification between $\M_G(c)$ with the space of gauge-equivalent solutions
of Hitchin's equations. We opt to work with $\M_G(c)$, because in this case we have more algebraic tools at our disposal. We shall make use of the tangent space of $\M_G(c)$, and we know from \cite{hitchin:1987} that the above identification induces a diffeomorphism between the corresponding tangent spaces.

We have the following result:
\begin{proposition}[\textbf{Hitchin} \cite{hitchin:1987}]\label{properness of f}\mbox{}
 \begin{enumerate}
 \item The function $f$ is proper.
 \item If $\M_G(c)$ is smooth, then $f$ is a non-degenerate perfect Bott-Morse function.
\end{enumerate}
\end{proposition}

Since $f$ is proper, it attains a minimum on each connected component of $\M_G(c)$. Moreover, we have the following result from general topology:
\begin{proposition}\label{proper}
If the subspace of local minima of $f$ is connected, then so is $\M_G(c)$.
\end{proposition}

Now, fix $L\in\Jac^d(X)$ and consider the space
$$\M_{\GL(n,\C)}(nd/2)\cong\M_{\GL(n,\C)}(nd/2)\times\{L\}\subset\M_{\GL(n,\C)}(nd/2)\times\Jac^d(X).$$
In our case, the Hitchin functional $$f:\M_{\GL(n,\C)}(nd/2)\longrightarrow\R$$
is given by
\begin{equation}\label{Hitchin proper function quadruples}
f(V,\Phi)=\|\Phi\|_{L^2}^2=\frac{\sqrt{-1}}{2}\int_{X}\mathrm{tr}(\Phi\wedge\Phi^*)\mathrm{dvol}.
\end{equation}

\subsection{Smooth minima} 

Away from the singular locus of
$\M_{\GL(n,\C)}(nd/2)$, the Hitchin functional $f$ is a moment map for the Hamiltonian
$S^1$-action on $\M_{\GL(n,\C)}(nd/2)$ given by
\begin{equation}\label{circle action Higgs}
(V,\Phi)\mapsto(V,e^{\sqrt{-1}\theta}\Phi).
\end{equation}
From this it follows immediately that a stable point of $\M_{\GL(n,\C)}(nd/2)$ is a
critical point of $f$ if and only if is a fixed point of the
$S^1$-action. Let us then study the fixed point set of the given
action (this is analogous to \cite{hitchin:1992} and \cite{bradlow-garcia-prada-gothen:2004}).

Let $(V,\Phi)$ represent a stable fixed point. Then either $\Phi=0$
or (since the action is on $\M_{\GL(n,\C)}(nd/2)$) there is a
one-parameter family of gauge transformations
$g(\theta)$ such that
$g(\theta)\cdot(V,\Phi)=(V,e^{\sqrt{-1}\theta}\Phi)$. In the latter case, let
$$\psi=\frac{d}{d\theta}g(\theta)|_{\theta=0}$$ be the infinitesimal gauge
transformation generating this family. Simpson shows in \cite{simpson:1992}
that this $(V,\Phi)$ is what is called a \emph{complex variation of
Hodge structure}. This means that
$$V=\bigoplus_j F_j$$ where the $F_j$'s are the eigenbundles of
the infinitesimal gauge transformation $\psi$: over $F_j$,
\begin{equation}\label{psi over Fj}
\psi=\sqrt{-1}j\in\C.
\end{equation}
$\Phi_j=\Phi|_{F_j}$ is a map
$$\Phi_j:F_j\longrightarrow F_{j+1}\otimes K$$ which is non-zero
for all $j$ except the maximal one.

Set
$$\M_L=\{(V,L',Q,\Phi)\in\M_d\ |\ L'\cong L\}\subset\M_d.$$
From Lemma \ref{fechado}, we know that $\M_d$ is a closed subspace of
$\M_{\GL(n,\C)}(nd/2)\times\Jac^d(X)$ hence, for each $L$, $\M_L$ is
closed in
$\M_{\GL(n,\C)}(nd/2)\times\{L\}\cong\M_{\GL(n,\C)}(nd/2)$.
The following proposition is a direct consequence of this and of the
properness of the $f$ given in (\ref{Hitchin proper function quadruples}).

\begin{proposition}
The restriction of the Hitchin functional $f$ to $\M_L$ is a proper and
bounded below function.
\end{proposition}

 From now on we will consider the restriction
of $f$ to $\M_L$. This fact will be important in the counting of
components of each $\M_d$, as we shall see in Section \ref{compM}.

The circle action (\ref{circle action Higgs}) restricts to $\M_L$. So, if $(V,L,Q,\Phi)$ is
an $\EGL(n,\R)$-Higgs bundle such that $(V,\Phi)$ is stable and is a fixed point of the $S^1$-action (i.e.,
is a critical point of $f$), then it is a variation of Hodge
structure. In this case, $g(\theta)\in
H^0(X,\EO(n,\mathcal O))$ and, since the Lie algebra of $\EO(n,\C)$ is $\mathfrak{o}(n,\C)\oplus\C$, we have
$\psi\in H^0(X,\mathfrak{o}(n,\mathcal{O})\oplus\mathcal{O})$, therefore being
skew-symmetric with respect to $Q$. Thus, using (\ref{psi over Fj}) we have that, if $v_j\in F_j$ and
$v_l\in F_l$,
$$\sqrt{-1}jQ(v_j,v_l)=Q(\psi v_j,v_l)=-Q(v_j,\psi v_l)=-\sqrt{-1}lQ(v_j,v_l).$$
$F_j$ and $F_l$ are therefore orthogonal under $Q$ if $l\neq -j$,
and $q:V\to V^*\otimes L$ yields an isomorphism
\begin{equation}\label{FjcongF-j*L}
q|_{F_j}:F_j\stackrel{\cong}{\longrightarrow}F_{-j}^*\otimes L.
\end{equation}
This means that
$$V=F_{-m}\oplus\dots\oplus F_m$$ for some $m$ integer or half-integer.

Using these isomorphisms and the fact that $\Phi$ is symmetric under
$Q$, we see that $$(q\otimes 1_K)\Phi_j=(\Phi_{-j-1}^t\otimes
1_K\otimes 1_L)q$$ for $j\in\{-m,\dots,m\}$.

The Cartan decomposition of $\liegc$ induces a decomposition of vector bundles
$$E_{H^\C}(\liegc)=E_{H^\C}(\liehc)\oplus E_{H^\C}(\liemc)$$
where $E_{H^\C}(\liegc)$ (resp. $E_{H^\C}(\liehc)$) is the adjoint bundle, associated to the adjoint representation of $H^\C$ on $\liegc$ (resp. $\liehc$).
For the group $\EGL(n,\R)$, we have
$E_{H^\C}(\liegc)=\text{End}(V)\oplus\mathcal{O}$ where
$\mathcal{O}=\text{End}(L)$ is the trivial line bundle on $X$ and
we already know that
$E_{H^\C}(\liehc)=\Lambda^2_QV\oplus\mathcal{O}$ and
$E_{H^\C}(\liemc)=S^2_QV$ where $\Lambda^2_QV$ is the bundle of
skew-symmetric endomorphisms of $V$ with respect to the form $Q$.
The involution in $\text{End}(V)\oplus\mathcal{O}$ defining
 the above decomposition
is $\theta\oplus 1_\mathcal{O}$ where $\theta$ is the
involution on $\text{End}(V)$ defined by
\begin{equation}\label{involution theta}
\theta(A)=-(qAq^{-1})^t\otimes 1_L.
\end{equation}
Its $+1$-eigenbundle is $\Lambda^2_QV\oplus\mathcal{O}$ and its
$-1$-eigenbundle is $S^2_QV$.

We
also have a decomposition of this vector bundle as
$$\text{End}(V)\oplus\mathcal{O}=\bigoplus_{k=-2m}^{2m}U_k\oplus\mathcal{O}$$
where
$$U_k=\bigoplus_{i-j=k}\Hom(F_j,F_i).$$ From (\ref{psi over Fj}), this is the $\sqrt{-1}k$-eigenbundle for the adjoint action
$\text{ad}(\psi):\text{End}(V)\oplus\mathcal{O}\to\text{End}(V)\oplus\mathcal{O}$
of $\psi$. We say that $U_k$ is the subspace of
$\text{End}(V)\oplus\mathcal{O}$ with \emph{weight} $k$.

Write
$$U_{i,j}=\Hom(F_j,F_i).$$ The restriction of the involution
$\theta$, defined in (\ref{involution theta}), to $U_{i,j}$ gives an isomorphism
\begin{equation}\label{theta sends Uij to U-i,-j}
\theta:U_{i,j}\to U_{-j,-i}
\end{equation}
so $\theta$ restricts to
$$\theta:U_k\longrightarrow U_k.$$

Write $$U^+=\Lambda^2_QV\ \text{ and }\ U^-=S^2_QV$$ so that
$$E_{H^\C}(\liehc)=U^+\oplus\mathcal{O}$$ and
$$E_{H^\C}(\liemc)=U^-.$$ Let also $$U_k^+=U_k\cap U^+$$ and
$$U_k^-=U_k\cap U^-$$ so that $U_k=U_k^+\oplus U_k^-$ is the
corresponding eigenbundle decomposition. Hence $$U^+=\bigoplus_k
U_k^+$$ and $$U^-=\bigoplus_k U_k^-.$$ Observe that $\Phi\in
H^0(X,U_1^-\otimes K)$.

The map $\text{ad}(\Phi)=[\Phi,-]$ interchanges $U^+$ with $U^-$ and
therefore maps $U_k^\pm$ to $U_{k+1}^\mp\otimes K$. So, for each
$k$, we have a weight $k$ subcomplex of the complex $C^\bullet(V,L,Q,\Phi)$ defined in Proposition \ref{deformation complex for quadruples}:
$$C^\bullet_k(V,L,Q,\Phi):U_k^+\oplus\mathcal{O}\xrightarrow{[\Phi,-]}U_{k+1}^-\otimes K.$$
In fact, since $\text{ad}(\psi)|_\mathcal{O}=0$, $C^\bullet_k(V,L,Q,\Phi)$ is given by $$C^\bullet_0(V,L,Q,\Phi):U^+_0\oplus\mathcal{O}\xrightarrow{[\Phi,-]}U^-_1\otimes K$$
and, for $k\neq 0$, by
$$C^\bullet_k(V,L,Q,\Phi):U^+_k\xrightarrow{[\Phi,-]}U^-_{k+1}\otimes K.$$

From Proposition \ref{deformation complex for quadruples}, if an $\EGL(n,\R)$-Higgs bundle $(V,L,Q,\Phi)$ is such that $(V,\Phi)$ is stable, its infinitesimal deformation
space is
$$\mathbb{H}^1(X,C^\bullet(V,L,Q,\Phi))=\bigoplus_k\mathbb{H}^1(X,C^\bullet_k(V,L,Q,\Phi)).$$
We say that $\mathbb{H}^1(X,C^\bullet_k(V,L,Q,\Phi))$ is the subspace of
$\mathbb{H}^1(X,C^\bullet(V,L,Q,\Phi))$ with \emph{weight} $k$.

By Hitchin's computations in \cite{hitchin:1992}, we
have the following result which gives us a way to compute the
eigenvalues of the Hessian of the Hitchin functional $f$ at a smooth (here we mean smooth
in $\M_{\GL(n,\C)}(nd/2)$) critical point.

\begin{proposition}\label{eigen}
Let $f$ be the Hitchin functional. Let $(V,L,Q,\Phi)$ be an $\EGL(n,\R)$-Higgs bundle with $(V,\Phi)$ stable and which represents a critical
point of $f$. The eigenspace of the Hessian of $f$
corresponding to the eigenvalue $k$ is
$$\mathbb{H}^1(X,C^\bullet_{-k}(V,L,Q,\Phi)).$$ In particular,
$(V,L,Q,\Phi)$ is a local minimum of $f$ if and only if
$\mathbb{H}^1(X,C^\bullet(V,L,Q,\Phi))$ has no subspaces with
positive weight.
\end{proposition}

For the moment we will only care about the stable points of $\M_L$.

Using Proposition \ref{eigen}, one can prove the following result
by an argument analogous to the proof of Corollary
4.15 of \cite{bradlow-garcia-prada-gothen:2003} (see also Remark $4.16$ in the same paper and Lemma $3.11$ of \cite{bradlow-garcia-prada-gothen:2008}). It is the fundamental result which makes possible the
description of the stable local minima of $f$.

\begin{theorem}\label{ad}
Let $(V,L,Q,\Phi)\in\M_L$ be a critical point of $f$ with $(V,\Phi)$ stable. Then
$(V,L,Q,\Phi)$ is a local minimum if and only if either $\Phi=0$ or
$$\mathrm{ad}(\Phi)|_{U_k^+}:U_k^+\longrightarrow U_{k+1}^-\otimes K$$
is an isomorphism for all $k\geq 1$.
\end{theorem}

The following theorem is quite similar to the corresponding one in \cite{hitchin:1992} and in \cite{bradlow-garcia-prada-gothen:2004} as one would naturally expect.
Indeed, the proof of this theorem is inspired in the one of Theorem 4.3 of \cite{bradlow-garcia-prada-gothen:2004}.

\begin{theorem}\label{min}
Let the $\EGL(n,\R)$-Higgs bundle $(V,L,Q,\Phi)$ be a a critical point of the Hitchin functional
$f$ such that $(V,\Phi)$ is stable. Then $(V,L,Q,\Phi)$ represents a local minimum if and
only if one of the following conditions occurs:
\begin{enumerate}
  \item $\Phi=0$.
  \item For each $i$, $\mathrm{rk}(F_i)=1$ and
  $\Phi_i$ is an isomorphism, for $i\neq m$.
\end{enumerate}
\end{theorem}
\begin{proof} The proof that a local minimum of $f$ must be of one of the above types is very similar to the one presented in the proof of Theorem 4.3 of \cite{bradlow-garcia-prada-gothen:2004}, so we skip it.

To prove the converse, let $(V,L,Q,\Phi)$ represent a point of type
(2). Then 
\begin{equation}\label{decomp of V}
V=\bigoplus_{i=-m}^m F_i
\end{equation}
with $\mathrm{rk}(F_i)=1$, so
$n=2m+1$, and 
\begin{equation}\label{decomp of Phi}
\Phi=\bigoplus_{i=-m}^m\Phi_i
\end{equation}
with $\Phi_i:F_i\to F_{i+1}\otimes K$ isomorphism, if $i\neq m$.

For each $k\in\{1,\dots,2m\}$, $\mathrm{rk}(U_k)=2m-k+1$ hence
$$\mathrm{rk}(U_k^+)=\frac{2m-k+1}{2}=\mathrm{rk}(U_{k+1}^-\otimes K)$$
if $n=k\ \text{mod}\ 2$, and
$$\mathrm{rk}(U_k^+)=\frac{2m-k}{2}=\mathrm{rk}(U_{k+1}^-\otimes K)$$
if $n\neq k\ \text{mod}\ 2$.
Therefore, if we prove that $\text{ad}(\Phi):U_k^+\to
U_{k+1}^-\otimes K$ is injective, we conclude that it is an
isomorphism and, from Theorem \ref{ad}, that $(V,L,Q,\Phi)$
represents a local minimum of $f$.

Let $g\in U_k^+=U_k\cap U^+=\bigoplus_{i-j=k}\Hom(F_j,F_i)\cap\Lambda^2_QV$. We can write $g$ as
\begin{equation}\label{decomp of g}
g=g_{-m}\oplus g_{-m+1}\oplus\dots\oplus g_{m-k}
\end{equation}
where
$g_j:F_j\to F_{j+k}$ and
$g_j=-q^{-1}(g_{-j-k}^t\otimes 1_L)q$. Now,
$$\text{ad}(\Phi)(g)=[\Phi,g]=\Phi g-(g\otimes
1_K)\Phi$$ and, using the decompositions (\ref{decomp of V}), (\ref{decomp of Phi}) and (\ref{decomp of g}),
this yields
$$[\Phi,g]=(\Phi_{-m+k} g_{-m}-(g_{-m+1}\otimes 1_K)\Phi_{-m})\oplus\dots\oplus(\Phi_{m-1}g_{m-k-1}-(g_{m-k}\otimes 1_K)\Phi_{m-k-1}).$$
The summands lie in different $U_{i,j}^-\otimes K$, hence
$[\Phi,g]=0$ is equivalent to the following system of equations
\begin{equation}\label{sist}
\begin{cases}
\Phi_{-m+k} g_{-m}-(g_{-m+1}\otimes 1_K)\Phi_{-m}=0\\
\Phi_{-m+k+1} g_{-m+1}-(g_{-m+2}\otimes 1_K)\Phi_{-m+1}=0\\
\vdots\\
\Phi_{m-1} g_{m-k-1}-(g_{m-k}\otimes 1_K)\Phi_{m-k-1}=0.
\end{cases}
\end{equation}

Take any fibre of $V$ and choose suitable basis of $V$ and
$V^*\otimes L$ such that, with respect to these basis,
$$\Phi=\begin{pmatrix}
    0 & \dots & \dots & \dots & 0 \\
    1 &  \ddots & & & 0 \\
    0 & 1 &  &  & \vdots \\
    \vdots & \ddots & \ddots &  & \vdots \\
    0 & \dots & 0 & 1 & 0 \
  \end{pmatrix},\hspace{2cm}q=\begin{pmatrix}
    0 & \dots & \dots & 0 & 1 \\
    \vdots &  & & \dots & 0 \\
    \vdots &  & 1 &  & \vdots \\
    0 & \dots &  &  & \vdots \\
    1 & 0 & \dots & \dots & 0 \
  \end{pmatrix}$$ and 
\begin{equation}\label{gj=-g-j-k}
g_j=-g_{-j-k}
\end{equation}
over the corresponding fibre of $U_k^+$. Then (\ref{sist}) implies that, over this fibre,
  $g_i=g_j$ for all $i,j$. In particular, 
\begin{equation}\label{gj=g-j-k}
g_j=g_{-j-k}
\end{equation}
for all $j$. From (\ref{gj=-g-j-k}) and (\ref{gj=g-j-k}), we must then have $g=0$. Since we considered any fibre, the result follows.
\end{proof}

\begin{remark}\label{rank}
Although we are always assuming $\mathrm{rk}(V)\geq 4$ even, we will need
during the proof of Proposition \ref{poly} below, to consider
$\EGL(n,\R)$-Higgs bundles of rank $1$ and $2$ and also of rank bigger or equal than $3$ odd. In the first two cases it is straightforward to see
that the minima of the Hitchin functional
$f(V,\Phi)=\|\Phi\|_{L^2}^2$, with $(V,\Phi)$ stable, in the corresponding moduli spaces are the following:
\begin{itemize}
\item If $\mathrm{rk}(V)=1$, $(V,L,Q,\Phi)$ is a minimum of $f$ if and
only if $\Phi=0$;
\item If $\mathrm{rk}(V)=2$, $(V,L,Q,\Phi)$ is a minimum of $f$ if and
only if either $\Phi=0$ or $V=F\oplus (F^*\otimes L)$ with
$\mathrm{rk}(F)=1$ and
$$\Phi=\left(\begin{array}{cc}
0 & 0 \\
\Phi' & 0 \\
\end{array}
\right)$$ with $\Phi':F\to F^*\otimes L$ non-zero (not
necessarily isomorphism).
\end{itemize}
For $\rk(V)\geq 3$ odd, $(V,L,Q,\Phi)$ is a minimum of $f$ if and
only if either $\Phi=0$ or $V=\oplus_{i=-m}^m F_i$ with $\mathrm{rk}(F_i)=1$ and
  $\Phi_i$ is an isomorphism, for $i\neq m$. This case is completely analogous to the even case considered here. The details can be found in \cite{oliveira:2008}.
\end{remark}

 Let $(V,L,Q,\Phi)$ represent a local minimum of
$f$ of type (2) of Theorem \ref{min}. Then,
\begin{equation}\label{mineven}
V=F_{-m}\oplus\dots\oplus F_{-1/2}\oplus F_{1/2}\oplus\dots\oplus
F_m
\end{equation}
where $m$ is an half-integer.

\begin{corollary}\label{min2}
Let $(V,L,Q,\Phi)$ represent a local minimum of $f$ of type
$\mathrm{(2)}$.
\begin{enumerate}
  \item Then $F_{-1/2}^2\cong LK$ and the others $F_i$ are uniquely
  determined by the choice of this square root of $LK$ as
  $F_{-1/2+i}\cong F_{-1/2}K^{-i}$.
  \item Then $(V,L,Q,\Phi)$ is isomorphic to an $\EGL(n,\R)$-Higgs bundle
  where
  $$q=\begin{pmatrix}
    0 & \dots & \dots & 0 & 1 \\
    \vdots &  & & \dots & 0 \\
    \vdots &  & 1 &  & \vdots \\
    0 & \dots &  &  & \vdots \\
    1 & 0 & \dots & \dots & 0 \
  \end{pmatrix}\ \text{ and }\ \Phi=\begin{pmatrix}
    0 & \dots & \dots & \dots & 0 \\
    1 &  \ddots & & \dots & 0 \\
    0 & 1 &  &  & \vdots \\
    \vdots & \ddots & \ddots &  & \vdots \\
    0 & \dots & 0 & 1 & 0 \
  \end{pmatrix}$$
  with respect to the decomposition $V=F_{-m}\oplus\dots\oplus F_m$.
\end{enumerate}
\end{corollary}

\subsection{Singular minima}

We must now show that Theorem \ref{min} gives us all non-zero minima of the Hitchin proper function $f$.

Let $(V,L,Q,\Phi)$ be an $\EGL(n,\R)$-Higgs bundle such that $(V,\Phi)$ is strictly polystable, with
$(V,\Phi)=\bigoplus_i(V_i,\Phi_i)$. Suppose moreover that $Q$ also
splits accordingly $Q=\bigoplus_i Q_i$ so that we have $\EGL(n,\R)$-Higgs bundles
$(V_i,L,Q_i,\Phi_i)$. We have $$f(V,L,Q,\Phi)=\sum_i
f(V_i,L,Q_i,\Phi_i)$$ so, if $(V,L,Q,\Phi)$ is a local minimum of
$f$, each of its stable summands is also a local minimum of $f$ in
the corresponding lower rank space $\M_L$. Hence each
$(V_i,L,Q_i,\Phi_i)$ is a fixed point of the circle action and
therefore the same happens to $(V,L,Q,\Phi)$. So $(V,L,Q,\Phi)$ is
a complex variation of Hodge structure $$V=\bigoplus_\alpha
W_\alpha$$ where each $W_\alpha$ is the
$\sqrt{-1}\alpha$-eigenbundle for an infinitesimal
$\EO(n,\C)$-gauge transformation $\psi$
and where $\Phi_\alpha:W_\alpha\to W_{\alpha+1}\otimes
K$, with the possibility that $\Phi_\alpha=0$. We can then also write
$$\text{End}(V)\oplus\mathcal{O}=\bigoplus_\lambda
U_\lambda\oplus\mathcal{O}$$ where $U_\lambda$ is the
$\sqrt{-1}\lambda$ eigenbundle of $\text{ad}(\psi)$. Let
$U_\lambda^{\pm}=U_\lambda\cap U^{\pm}$, where $U^+=\Lambda^2_QV$
and $U^-=S^2_QV$, and define the following complex of sheaves associated to
$(V,L,Q,\Phi)$:
\begin{equation}\label{complex >0}
C^\bullet_{>0}(V,L,Q,\Phi):\bigoplus_{\lambda>0}U_\lambda^+\xrightarrow{[\Phi,-]}\bigoplus_{\lambda>1}U_\lambda^-\otimes K.
\end{equation}

Hitchin's computations in \cite{hitchin:1992} for showing that a
given fixed point of the circle action is not a local minimum yield
the following proposition.

\begin{proposition}\label{notmin}
Let $(V,L,Q,\Phi)$ be a fixed point of the $S^1$-action on $\M_L$. Let $(V_t,L,Q_t,\Phi_t)$ be a one-parameter family of polystable
$\EGL(n,\R)$-Higgs bundles such that $(V_0,L,Q_0,\Phi_0)=(V,L,Q,\Phi)$. If there is
a non-trivial tangent vector to the family at $0$ which lies in the
subspace $$\mathbb{H}^1(X,C^\bullet_{>0}(V,L,Q,\Phi))$$ of the
infinitesimal deformation space
$\mathbb{H}^1(X,C^\bullet(V,L,Q,\Phi))$, then
$(V,L,Q,\Phi)$ is not a local minimum of $f$.
\end{proposition}

In other words, if $(V,\Phi)$ is strictly polystable, Hitchin's arguments in
\cite{hitchin:1992} are also valid: if there is a non-empty subspace
of $\mathbb{H}^1(X,C^\bullet(V,L,Q,\Phi))$ which gives directions
in which $f$ decreases and if these directions are integrable into a
one-parameter family in $\M_L$, then $(V,L,Q,\Phi)$ is not a local
minimum of $f$.

The following result, adapted from \cite{hitchin:1992}, shows that there are no more non-zero minima of $f$ besides the ones of Theorem \ref{min}.

\begin{proposition}\label{poly}
Let $(V,L,Q,\Phi)$ represent a point of $\M_L$ such that $(V,\Phi)$
is strictly polystable. If $\Phi\neq 0$, then $(V,L,Q,\Phi)$ is not a
local minimum of $f$.
\end{proposition}
\begin{proof} Suppose $V=V_1\oplus V_2$, $\Phi=\Phi_1\oplus\Phi_2$ and
$(V,\Phi)$ represents a local minimum of $f$ in $\M_L$, with
$\Phi_1\neq 0\neq\Phi_2$.

Consider first the case where $V_1$ and $V_2$ are not isomorphic and
$V_1\cong V_1^*\otimes L$ and $V_2\cong V_2^*\otimes L$. Then the
quadratic form $Q$ also splits as $Q=Q_1\oplus Q_2$ with $Q_i:V_i\otimes
V_i\to L$, $i=1,2$. We have therefore the $\EGL(n,\R)$-Higgs bundles
$(V_1,L,Q_1,\Phi_1)$ and $(V_2,L,Q_2,\Phi_2)$ which are local
minima of $f$ on the corresponding lower rank moduli space.
Let $n_1=\rk(V_1)$ and $n_2=\rk(V_2)$ so that $n=n_1+n_2$ (here,
the cases $n_1=2$ or $n_2=2$ or $n_i\geq 3$ odd are included). So we have
$$V_1=F_{-m}\oplus\dots\oplus F_m$$ and
$$V_2=G_{-k}\oplus\dots\oplus G_k.$$

Consider the complex
$$C^\bullet_{m+k}(V,L,Q,\Phi):U^+_{m+k}\xrightarrow{[\Phi,-]}
U^-_{m+k+1}\otimes K.$$ Since $\Phi\neq 0$, we have $m+k>0$ and
$C^\bullet_{m+k}(V,L,Q,\Phi)$ is a subcomplex of the complex
$C^\bullet_{>0}(V,L,Q,\Phi)$ defined in (\ref{complex >0}).

Consider the space
$H^1(X,\Hom(G_{-k},F_m))=H^1(X,F_mG_k
L^{-1})$. For $i=1,2$, $$\deg(V_i)=n_i\deg(L)/2$$ and, since $F_m$
(resp. $G_k$) is a $\Phi_1$ (resp. $\Phi_2$)-invariant subbundle
of $V_1$ (resp. $V_2$), we have, from the stability of $(V_1,\Phi_1)$
and of $(V_2,\Phi_2)$,
 $$\deg(F_m G_k
L^{-1})=\deg(F_m)+\deg(G_k)-\deg(L)<0.$$ It follows, by Riemann-Roch,
that $H^1(X,\Hom(G_{-k},F_m))$ is non-zero. Choose then $$0\neq
h\in H^1(X,\Hom(G_{-k},F_m))$$ and let
\begin{equation}\label{sigma lies in positive subspace}
\sigma=(h,\theta_*(h))\in
H^1(X,\Hom(G_{-k},F_m)\oplus\Hom(F_{-m},G_k)\cap\Lambda_Q^2
V)\subset H^1(X,U_{m+k}^+)
\end{equation}
where
$\theta_*:H^1(X,\text{End}(V))\to H^1(X,\text{End}(V))$
is the map induced by the involution $\theta$ on $\mathrm{End}(V)$
previously defined. $\sigma$ is obviously non-zero and, moreover, it is
annihilated by $$\text{ad}(\Phi)=[\Phi,-]:H^1(X,U_{m+k}^+)\longrightarrow
H^1(X,U_{m+k+1}^-\otimes K)$$ hence it defines an element in
$\mathbb{H}^1(X,C^\bullet_{m+k}(V,L,Q,\Phi))$, which we also
denote by $\sigma$.

Now, $\sigma$ defines extensions
$$0\longrightarrow F_m\stackrel{\ i_\sigma\ }{\longrightarrow}U_\sigma\stackrel{\ p_\sigma\ }{\longrightarrow}G_{-k}\longrightarrow 0$$
and
$$0\longrightarrow G_k\xrightarrow{\ p_\sigma^t\otimes 1_L\ }U_\sigma^*\otimes L\xrightarrow{\ i_\sigma^t\otimes 1_L\ }F_{-m}\longrightarrow 0.$$

Let
\begin{equation}\label{deformation of V}
V_\sigma=\bigoplus_{i=-m+1}^{m-1}F_i\oplus U_\sigma\oplus\bigoplus_{j=-k+1}^{k-1}G_j\oplus (U_\sigma^*\otimes
L)
\end{equation}
and $\Phi_\sigma:V_\sigma\to V_\sigma\otimes K$
given by
\begin{equation}\label{deformation of Higgs field}
\begin{split}
  &\Phi_\sigma(v_{-m+1},\dots,v_{m-1},u_\sigma,w_{-k+1},\dots,w_{k-1},u_\sigma^*\otimes
l)=\\
=&\,(\Phi_1v_{-m+1},\dots,(i_\sigma\otimes
1_K)\Phi_1v_{m-1},\Phi_2p_\sigma
u_\sigma,\\
&\Phi_2w_{-k+1},\dots,(p_\sigma^t\otimes 1_L\otimes
1_K)\Phi_2w_{k-1},\Phi_1(i_\sigma^t\otimes 1_L)(u_\sigma^*\otimes
l)).
\end{split}
\end{equation}

Let us see that $(V_\sigma,\Phi_\sigma)$ is stable. If $W$ is
a proper $\Phi_\sigma$-invariant subbundle of $V_\sigma$ then $W$ is
one of the following:
\begin{itemize}
    \item $W=F_m$;
    \item $W=G_k$;
    \item $W=\bigoplus_{i=-m+a}^{m-1}F_i\oplus F_m,\text{ with }1\leq
a\leq 2m-1$;
    \item $W=\bigoplus_{j=-k+b}^{k-1}G_j\oplus G_k,\text{ with }1\leq
b\leq 2k-1$;
    \item $W=\bigoplus_{i=-m+a}^{m-1}F_i\oplus U_\sigma\oplus\bigoplus_{j=-k+1}^{k-1}G_j\oplus G_k,\text{ with
}1\leq a\leq 2m-1$;
    \item $W=U_\sigma\oplus\bigoplus_{j=-k+1}^{k-1}G_j\oplus G_k$.
\end{itemize} Using the stability of $(V_1,\Phi_1)$ or of
$(V_2,\Phi_2)$ and the fact that $\mu(V_i)=\mu(V)=\mu(V_\sigma),\
i=1,2$, it follows that $\mu(W)<\mu(V_\sigma)$,
$(V_\sigma,\Phi_\sigma)$ being therefore stable.

The summands $\bigoplus_{i=-m+1}^{m-1}F_i$ and
$\bigoplus_{j=-k+1}^{k-1}G_j$ in $V_\sigma$ have a quadratic form
coming from $Q$, and we also have the canonical $L$-valued quadratic form on
$U_\sigma\oplus(U_\sigma^*\otimes L)$. These give a $L$-valued
quadratic form $Q_\sigma$ on $V_\sigma$. 

So we have seen that $V_\sigma$ defined in (\ref{deformation of V}), $\Phi_\sigma$ defined in (\ref{deformation of Higgs field}) and $Q_\sigma$ just defined give rise to a stable $\EGL(n,\R)$-Higgs bundle $(V_\sigma,L,Q_\sigma,\Phi_\sigma)$. 

Notice that, if $\sigma=0$, then $(V_0,L,Q_0,\Phi_0)=(V,L,Q,\Phi)$. Now,
consider the family $(V_{t\sigma},L,Q_{t\sigma},\Phi_{t\sigma})$
of $\EGL(n,\R)$-Higgs bundles. The induced infinitesimal deformation is given by
$\sigma$ which, from (\ref{sigma lies in positive subspace}), lies in a positive weight subspace of
$\mathbb{H}^1(X,C^\bullet(V,L,Q,\Phi))$. Taking Proposition
\ref{notmin} in account, this proves that $(V,L,Q,\Phi)$ is not a local
minimum of $f$.

Suppose now that $V_1\not\cong V_2$, but the form $Q$ does not
decompose.  From the stability of $(V_1,\Phi_1)$ and of
$(V_2,\Phi_2)$ we must have $$q=\left(%
\begin{array}{cc}
  0 & q_{12} \\
  q_{21} & 0 \\
\end{array}%
\right)$$ where $q_{12}:V_2\to V_1^*\otimes L$ and
$q_{21}:V_1\to V_2^*\otimes L$ are isomorphisms and
$q_{21}=q_{12}^t\otimes 1_L$.

Hence we can write $$V=V_1\oplus\left(V_1^*\otimes L\right)$$ and
$$\Phi=\Phi_1\oplus\left(\Phi_1^t\otimes 1_K\otimes 1_L\right).$$ Consider the
point in $\M_{\GL(n,\C)}(nd/2)$ represented by $(V_1,\Phi_1)$. Since $\Phi_1\neq
0$, we know from \cite{hitchin:1987} that $(V_1,\Phi_1)$ is not a local
minimum of $f$ in $\M_{\GL(n,\C)}(nd/2)$ (this is because the group $\GL(n,\C)$ is complex). Therefore one can find a
family $(V_{1,s},\Phi_{1,s})$ of stable Higgs bundles near
$(V,\Phi)$ such that $f(V_{1,s},\Phi_{1,s})< f(V,\Phi)$ for all
$s$, i.e., 
\begin{equation}\label{Phi1s<Phi}
\|\Phi_{1,s}\|_{L^2}^2<\|\Phi\|_{L^2}^2.
\end{equation}

Consider now the family of $\EGL(n,\R)$-Higgs bundles in $\M_L$ given by
$$(V_{1,s}\oplus\left(V_{1,s}^*\otimes
L\right)
L,L,Q_s,\Phi_{1,s}\oplus\Phi_{1,s}^t\otimes 1_K\otimes 1_L)$$ where
$Q_s$ is the canonical quadratic form in $V_{1,s}\oplus\left(V_{1,s}^*\otimes
L\right) L$. We have 
\begin{equation}\label{normPhi=..oplus..}
\|\Phi_{1,s}\oplus\left(\Phi_{1,s}^t\otimes 1_K\otimes
1_L\right)\|_{L^2}^2=\|\Phi_{1,s}\|_{L^2}^2+\|\Phi_{1,s}^t\otimes
1_K\otimes 1_L\|_{L^2}^2
\end{equation}
where we are using the harmonic metric on
$V_{1,s}^*$ and on $V_{1,s}\oplus V_{1,s}^*$ induced by the
one on $V_{1,s}$. We have
$\mathrm{tr}((\Phi_{1,s}^t\otimes 1_K\otimes
1_L)(\Phi_{1,s}^t\otimes 1_K\otimes
1_L)^*)=\mathrm{tr}(\Phi_{1,s}\Phi_{1,s}^*)$ therefore (\ref{normPhi=..oplus..}) is equivalent to
$$\|\Phi_{1,s}\oplus\Phi_{1,s}^t\otimes 1_K\otimes 1_L\|_{L^2}^2=
2\|\Phi_{1,s}\|_{L^2}^2$$
and from (\ref{Phi1s<Phi}) we conclude that
$$\|\Phi_{1,s}\oplus\Phi_{1,s}^t\otimes 1_K\otimes 1_L\|_{L^2}^2<2\|\Phi_1\|_{L^2}^2=\|\Phi\|_{L^2}^2$$ for all $s$. Hence $(V_1\oplus
V_2,L,Q,\Phi_1\oplus\Phi_2)$ is not a local minimum of $f$.

If $V_1\cong V_2$, then we saw in the proof of Lemma \ref{lema} that
we can decompose $Q=Q_1\oplus Q_2$ so that we can decompose the $\EGL(n,\R)$-Higgs bundles
$(V,L,Q,\Phi)$ as $(V_1,L,Q_1,\Phi'_1)\oplus(V_1,L,Q_2,\Phi'_2)$.
Hence we use the same argument as the first case to prove that
$(V,L,Q,\Phi)$ is not a minimum of $f$.

If, for example, $\Phi_1\neq 0$ and $\Phi_2=0$ then, due to the
symmetry of $\Phi$ relatively to $Q$, the quadratic form must split
into $Q_1\oplus Q_2$, so that we have $(V_1,L,Q_1,\Phi_1)$ and
$(V_2,L,Q_2,0)$ and in a similar manner to the first case
considered, we prove that $(V,L,Q,\Phi)$ is not a local minimum
of $f$.

For $\EGL(n,\R)$-Higgs bundles such that $(V,\Phi)$ has more than two summands, just consider
the first two and use one of the above arguments.
\end{proof}

%%%%%%%%%%%%%%%%%%%%%%%%%%%%%%%%%%%%%%%%%%%%%%%%
\section{Connected components of the space of $\EGL(n,\R)$-Higgs bundles}\label{compM}
%%%%%%%%%%%%%%%%%%%%%%%%%%%%%%%%%%%%%%%%%%%%%%%%

In this section we compute the number of components of the subspaces of the moduli space of $\EGL(n,\R)$-Higgs bundles such that
the degree of $L$ is $0$ and $1$. Denote these subspaces by $\M_0$ and 
$\M_1$, respectively. In other words, using Notation \ref{submodquad}, we write $\M_0$ as a disjoint union $$\M_0=\bigsqcup_{w_2\in\Z_2}\M(0,(w_2,0))\sqcup\bigsqcup_{\overline{\mu}_1\in\Z_2^{2g}\setminus\{0\}}\M(\overline{\mu}_1,0).$$
On the other hand, $$\M_1=\bigsqcup_{\overline{\mu}_1\in\Z_2^{2g}}\M(\overline{\mu}_1,1).$$

Of course, the space $\M$ of isomorphism classes of polystable $\EGL(n,\R)$-Higgs bundles
has an infinite number of components
because the invariant given by the degree of $L$ can be any integer.
But our computation will also give the number of components of any
subspace of $\M$ with the degree of $L$ fixed, due to the
identifications given in Remark \ref{identif}.

Before proceeding with the computation, we need some results which will be used.
Let $\N_{\EO(n,\C)}$ be the
moduli space of holomorphic semistable principal $\EO(n,\C)$-bundles on $X$ and
$\N_{\EO(n,\C)}(\overline{\mu}_1,\overline{\mu}_2)$ be the subspace with
invariants $(\overline{\mu}_1,\overline{\mu}_2)$.

The following is an adaptation of Proposition 4.2 of \cite{ramanathan:1975}.

\begin{proposition}
$\N_{\EO(n,\C)}(\overline{\mu}_1,\overline{\mu}_2)$ is connected.
\end{proposition}
\begin{proof} Let $E'$ and $E''$ represent two classes in
$\N_{\EO(n,\C)}(\overline{\mu}_1,\overline{\mu}_2)$. Let $P$ be the
underlying $C^\infty$ principal bundle, and let
$\overline{\partial}_{A'}$ and $\overline{\partial}_{A''}$ be the
operators on $P$ defining, respectively, $E'$ and $E''$ and given by
unitary connections $A'$ and $A''$.

Let $\mathbb D$ be an open disc in $\C$ containing 0 and 1. Consider the
$C^\infty$ principal-$\EO(n,\C)$ bundle
$\mathbb E\to\mathbb D\times X$, where $\mathbb E=\mathbb D\times P$. Define the
connection form on $\mathbb{E}$ by
$$A_z(v,w)=zA''(w)+(1-z)A'(w)\in\Omega^1(\mathbb D\times P,\mathfrak{o}(n,\C)\oplus\C)$$
where $v$ is tangent to $\mathbb D$ at $z$ and $w$ is tangent to $P$ at some
point $p$. If we consider the holomorphic bundle $E_z$ given by
$\mathbb{E}|_{\{z\}\times X}$ with the holomorphic structure given
by $A_z$, then we have that $E_0\cong E'$ and $E_1\cong E''$.

Since semistability is an open condition with respect to the Zariski
topology, $\mathbb D\setminus D'$ is connected where $D'=\{z\in \mathbb D:E_z\text{ is not
semistable}\}$. Hence $\{E_z\}_{z\in \mathbb D\setminus D'}$ is a connected family of
semistable $\EO(n,\C)$-principal bundles joining $E_0$ and $E_1$. Since $E_0\cong E'$ and
$E_1\cong E''$, using the universal property of the coarse moduli space
$\N_{\GL(n,\C)}$ of $\GL(n,\C)$-principal bundles, there is
a connected family in $\N_{\GL(n,\C)}$ joining $E'$ and
$E''$. But, of course this connected family lies in
$\N_{\EO(n,\C)}(\overline{\mu}_1,\overline{\mu}_2)$.
\end{proof}

Let $$\M'_L$$ be the subspace of $\M_L$ consisting of those components of
$\M_L$ such that the minimum of the Hitchin function $f$ attained on these components is $0$.
Hence the local minima on $\M'_L$ are those with $\Phi=0$.

Proposition \ref{proper} and the previous one yield the following:

\begin{corollary}\label{Zero Higgs}
For each $(\overline{\mu}_1,\overline{\mu}_2)$, the space $\M'_L(\overline{\mu}_1,\overline{\mu}_2)$ is (if non-empty) connected.
\end{corollary}

Recall from Corollary \ref{nonempty cR} that $\M'_L(\overline{\mu}_1,\overline{\mu}_2)$ is empty precisely when $n$ and $\deg(L)$ are both odd. Excluding this case, $\M'_L(\overline{\mu}_1,\overline{\mu}_2)$ is hence connected.

All the analysis of the proper function $f$ carried in Section \ref{morse quadruples} was done over
each $\M_L$, hence one can use Proposition \ref{proper} to compute the number of components of $\M_L$, and then
compute the number of components of $\M_0$ and of $\M_1$.

For each $L$, define $$\M''_L=\M_L\setminus\M'_L$$ so that we have a disjoint union
\begin{equation}\label{disj}
\M_L=\M'_L\sqcup\M''_L.
\end{equation}

Let us now concentrate attentions on $\M_0$.

For each $L\in\Jac^0(X)=\Jac(X)$, the Jacobian of $X$, $\M_L$ is a subspace of $\M_0$ and it is the fibre over $L$ of
the map
\begin{equation}\label{nu0}
\nu_0:\M_0\longrightarrow\Jac(X)
\end{equation}
given by $$\nu_0(V,L,Q,\Phi)=L.$$
To emphasize the fact that now $\M_L\subset\M_0$, we shall write
$\M_{L,0}$ instead of $\M_L$.
Any two fibres $\M_{L,0}$ and $\M_{L',0}$ of $\nu_0$ are isomorphic through the map
$$(V,L,Q,\Phi)\mapsto(V\otimes L^{-1/2}\otimes L'^{1/2},L',Q\otimes 1_{L^{-1}}\otimes 1_{L'},\Phi\otimes 1_{L^{-1/2}}\otimes 1_{L'^{1/2}}).$$
In particular, any fibre is isomorphic to $\M_{\mathcal{O}}$.

More precisely, after lifting to a finite cover, (\ref{nu0}) becomes a product. This is a similar situation to the one which occurs on the moduli of vector bundles with fixed determinant (cf. \cite{atiyah-bott:1982}). Indeed, we have the following commutative diagram:
\begin{equation}\label{lift cover becomes product}
\xymatrix{\M_{\mathcal O}\times\Jac(X)\ar[d]_{\pi}\ar[r]^(0.55){\text{pr}_2}&\Jac(X)\ar[d]^{\pi'}\\
          \M_0\ar[r]^{\nu_0}&\Jac(X)}
\end{equation}
where $\pi((W,\mathcal O,Q,\Phi),M)=(W\otimes M,M^2,Q\otimes 1_{M^2},\Phi\otimes 1_M)$ and $\pi'(M)=M^2$. Hence $\nu_0$ is a fibration.

Recall that an $\EGL(n,\R)$-Higgs bundle $(V,\mathcal{O},Q,\Phi)$ is topologically classified by the invariants $(\overline\mu_1,\overline\mu_2)$ where $\overline{\mu}_1=w_1(V,Q,\Phi)\in\Z_2^{2g}$ and, if $\overline{\mu}_1\neq 0$, then $\overline{\mu}_2=0=\deg(\mathcal{O})$, and, if $\overline{\mu}_1=0$, then $\overline{\mu}_2=(w_2(V,Q,\Phi),0=\deg(\mathcal{O}))$.

Now, if $\M_{\GL(n,\R)}$ denotes the moduli space of $\GL(n,\R)$-Higgs bundles \cite{bradlow-garcia-prada-gothen:2004}, which are classified by the first and second Stiefel-Whitney classes, there is a surjective map
\begin{equation}\label{map from MGL to MO}
\M_{\GL(n,\R)}\longrightarrow\M_{\mathcal{O}}
\end{equation}
given by $(W,Q,\Phi)\mapsto (W,\mathcal O,Q,\Phi)$ and such that:
\begin{itemize}
 \item $\M_{\GL(n,\R)}(0,w_2)$ is mapped onto $\M_{\mathcal{O}}(0,(w_2,0))$;
 \item  if $w_1\neq 0$, $\M_{\GL(n,\R)}(w_1,w_2)$ is mapped onto $\M_{\mathcal{O}}(w_1,0)$.
\end{itemize}

The following result is proved in Proposition 4.6 of \cite{bradlow-garcia-prada-gothen:2004} and
gives a more detailed information about the structure of $\M_{\GL(n,\R)}$.
\begin{proposition}\label{minneq0GL(n,R)}
Let $(V,\mathcal{O},Q,\Phi)\in\M_{\GL(n,\R)}$ be a local minimum of $f$ with
$\Phi\neq 0$. Then, $$w_1(V,\mathcal{O},Q,\Phi)=0$$ and $$w_2(V,\mathcal{O},Q,\Phi)=(g-1)n^2/4\ \mathrm{mod}\ 2.$$
\end{proposition}

Therefore, using the surjection (\ref{map from MGL to MO}) and the fact that any fibre of $\nu_0$ is isomorphic to $\M_{\mathcal{O}}$, we obtain:

\begin{proposition}
Let $(V,L,Q,\Phi)\in\M_{L,0}$ be a local minimum of $f$ with
$\Phi\neq 0$. Then, $$\overline{\mu}_1(V,L,Q,\Phi)=0$$ and $$\overline{\mu}_2(V,L,Q,\Phi)=((g-1)n^2/4\ \mathrm{mod}\ 2,0).$$
\end{proposition}

From now on we shall write $$z_0=(g-1)\frac{n^2}{4}\ \mathrm{mod}\ 2.$$
From this proposition and from what we saw above follows that,
\begin{equation}\label{L0mu1neq0}
 \M_{L,0}(\overline{\mu}_1)=\M'_{L,0}(\overline{\mu}_1)
\end{equation}
if $\overline{\mu}_1\neq 0$,
\begin{equation}\label{L0mu1=0w2neqz0}
\M_{L,0}(0,w_2)=\M'_{L,0}(0,w_2)
\end{equation}
if $w_2\neq z_0$, and
\begin{equation}\label{L0mu1=0w2=z0}
 \M_{L,0}(0,z_0)=\M'_{L,0}(0,z_0)\sqcup\M''_{L,0}(0,z_0).
\end{equation}
In other words,
\begin{equation}\label{even0}
\M_{L,0}=\bigsqcup_{\overline{\mu}_1\in(\Z_2)^{2g}\setminus\{0\}}\M'_{L,0}(\overline{\mu}_1)\sqcup\bigsqcup_{w_2\in\Z_2}\M'_{L,0}(0,w_2)\sqcup\M''_{L,0}(0,z_0).
\end{equation}

\begin{proposition}\label{ccq}
Let $n\geq 4$ be even and $L\in\Jac(X)$ be given. Then $\M_{L,0}$ has $2^{2g+1}+1$ connected components. More precisely,
\begin{enumerate}
        \item $\M_{L,0}(\overline{\mu}_1)$ with $\overline{\mu}_1\neq 0$, is connected;
	\item $\M_{L,0}(0,w_2)$ with $w_2\neq z_0$, is connected;
        \item $\M_{L,0}(0,z_0)$ has $2^{2g}+1$ components.
\end{enumerate}
\end{proposition}

This result follows immediately from Theorem $5.2$ of \cite{bradlow-garcia-prada-gothen:2004} and from the existence of the map $\M_{\GL(n,\R)}\to\M_{\mathcal{O}}$ described in (\ref{map from MGL to MO}). However, for completeness, we are still going to give a proof.

\begin{proof} Let $L\in\Jac(X)$. Fix
$\overline{\mu}_1\neq 0$ and consider the subspace
$$\M_{L,0}(\overline{\mu}_1)\subset\M_{L,0}.$$
This space is connected by (\ref{L0mu1neq0}) and by Corollary \ref{Zero Higgs}. So there are $2^{2g}-1$ components of $\M_{L,0}$
of this kind.

For the same reason but using (\ref{L0mu1=0w2neqz0}), we see that $\M_{L,0}(0,w_2)$ with $w_2\neq z_0$, is connected.

For the space $\M_{L,0}(0,z_0)$ we have the decomposition (\ref{L0mu1=0w2=z0}). The space $\M'_{L,0}(0,z_0)$ is connected from Corollary \ref{Zero Higgs}. Let us then analyse the space $\M''_{L,0}(0,z_0)$.
Consider the non-zero local minima of the Hitchin functional $f$. From Corollary \ref{min2}, these are such that
\begin{equation}\label{projeven}
  V=F_{-1/2}\otimes\bigoplus_{i=-r-1}^r K^i
\end{equation}
where $r=m-1/2$ and $F_{-1/2}$ is a square root of $LK$. There are
$2^{2g}$ different choices for $F_{-1/2}$ thus the space of local minima of
this kind consists of $2^{2g}$ isolated points. Therefore
$\M''_{L,0}(0,z_0)$ has $2^{2g}$ connected components. All these are
homeomorphic to a vector space and constitute the so-called
\emph{Hitchin} or \emph{Teichmüller components} of $\M_{L,0}$ \cite{hitchin:1992}. So
$\M_{L,0}(0,z_0)$ has $2^{2g}+1$ components.

If follows from (\ref{even0}) and from the count above that $\M_{L,0}$ has $2^{2g+1}+1$
connected components.
\end{proof}

We have computed the components of each fibre of $\nu_0$. Let us see that the space $\M_0$ has less components than $\M_{L,0}$.

\begin{theorem}\label{ccM0}
The space $\M_0$ has $2^{2g}+2$ components.
\end{theorem}
\begin{proof}
From Theorem \ref{topquad}, there are $2^{2g}+1$ topological invariants of $\EGL(n,\R)$-Higgs bundles in $\M_0$, hence $\M_0$ has at least $2^{2g}+1$ components.

Let $(V,L,Q,\Phi),(V',L',Q',\Phi')\in\M_0(0,z_0)$ such that each $\EGL(n,\R)$-Higgs bundle is a local minimum
of type (2) on the corresponding fibre of $\nu_0$ (see (\ref{nu0})). Hence
$$V=F_{-m}\oplus\dots\oplus F_{-1/2}\oplus F_{1/2}\oplus\dots\oplus F_m$$ and $$V'=F'_{-m}\oplus\dots\oplus F'_{-1/2}\oplus F'_{1/2}\oplus\dots\oplus F'_m$$
where $F_{-1/2}$ (resp. $F'_{-1/2}$)
is a square root of $LK$ (resp. $L'K$).
Since $\Jac(X)$ is connected, there is a path $L_t$ in $\Jac(X)$ joining $L$ to $L'$. Set
$$V_t=F_{-m,t}\oplus\dots\oplus F_{-1/2,t}\oplus F_{1/2,t}\oplus\dots\oplus F_{m,t}$$
where $F_{-1/2,t}^2\cong L_tK$ and $F_{-1/2+i,t}\cong F_{-1/2,t}K^{-i}$. With
$$q_t=\begin{pmatrix}
    0 & \dots & \dots & 0 & 1 \\
    \vdots &  & & 1 & 0 \\
    \vdots &  & 1 &  & \vdots \\
    0 & \dots &  &  & \vdots \\
    1 & 0 & \dots & \dots & 0 \
  \end{pmatrix}\ \text{ and }\ \Phi_t=\begin{pmatrix}
    0 & \dots & \dots & \dots & 0 \\
    1 &  \ddots & & \dots & 0 \\
    0 & 1 &  &  & \vdots \\
    \vdots & \ddots & \ddots &  & \vdots \\
    0 & \dots & 0 & 1 & 0 \
  \end{pmatrix}$$ $(V_t,L_t,Q_t,\Phi_t)_t$ is a path in $\M_0$ joining $(V,L,Q,\Phi)$ and $(V',L',Q',\Phi')$ and such that, for every $t$, $(V_t,L_t,Q_t,\Phi_t)$
  is a minimum of $f$ in $\M_{L_t,0}$ of type (2). Hence we conclude that all the $2^{2g}$
  Hitchin components of all fibres of $\nu_0$ join together to form a unique component of $\M_0$: $\M''_0(0,z_0)=\bigcup_{L\in\Jac(X)}\M''_{L,0}(0,z_0)$. Note that this is not a Hitchin component. Indeed, the group $\EGL(n,\R)$ is not a split real form (due to $\U(1)$), so the moduli space of $\EGL(n,\R)$-Higgs bundles on $X$ was not expected to have a Hitchin component (cf. \cite{hitchin:1992}).

On the other hand, $\M'_0(\overline{\mu}_1)=\bigcup_{L\in\Jac(X)}\M'_{L,0}(\overline{\mu}_1)$ is connected because $\nu_0|_{\M'_0(\overline{\mu}_1)}:\M'_0(\overline{\mu}_1)\to\Jac(X)$ is surjective and with connected fibre from item (1)(a) of Proposition \ref{ccq} and $\Jac(X)$ is connected. For an analogous reason, we also conclude that $\M'_0(0,w_2)=\bigcup_{L\in\Jac(X)}\M'_{L,0}(0,w_2)$ is connected.

Finally, $\M'_0(0,z_0)$ and $\M''_0(0,z_0)$ are two different connected components of $\M_0(0,z_0)$.

Concluding, we have one component for each $\M'_0(0,z_0)$, $\M''_0(0,z_0)$ and $\M'_0(0,w_2)$ with $w_2\neq z_0$, and $2^{2g}-1$ components coming from $\M'_0(\overline{\mu}_1)$. These yield the $2^{2g}+2$ components of $\M_0$.
\end{proof}

Let us now deal with the space $\M_1$.

We have again a map $\nu_1:\M_1\to\Jac^1(X)$ and, if $\deg(L)=1$, $\M_{L,1}=\nu_1^{-1}(L)$. In fact, when we fix a line bundle $L_0\in\Jac^1(X)$, we have a analogous diagram to (\ref{lift cover becomes product}):
$$\xymatrix{\M_{L_0}\times\Jac(X)\ar[d]_{\pi}\ar[r]^(0.55){m}&\Jac^1(X)\ar[d]^{\pi'}\\
          \M_1\ar[r]^{\nu_1}&\Jac^1(X)}$$
where $m((W,L_0,Q,\Phi),M)=ML_0$, $\pi((W,L_0,Q,\Phi),M)=(W\otimes M,L_0M^2,Q\otimes 1_{M^2},\Phi\otimes 1_M)$ and $\pi'(L)=L^2L_0^{-1}$. Hence $\nu_1$ is also a fibration.

If an $\EGL(n,\R)$-Higgs bundle $(V,L,Q,\Phi)$, with $(V,\Phi)$ stable, is a non-zero local
minimum of $f$ in $\M_L$ then it follows from
Corollary \ref{min2} that $\deg(L)$ is even.
Hence, if $\deg(L)=1$,
$$\M_{L,1}(\overline{\mu}_1)=\M'_{L,1}(\overline{\mu}_1)$$ thus,
\begin{equation}\label{even1}
\M_{L,1}=\bigsqcup_{\overline{\mu}_1\in(\Z_2)^{2g}}\M'_{L,1}(\overline{\mu}_1).
\end{equation}

\begin{proposition}\label{ccq1}
Let $n\geq 4$ be even and let $L\in\Jac^1(X)$. Then $\M_{L,1}$ has $2^{2g}$
connected components. More precisely, each
$\M_{L,1}(\overline{\mu}_1)$ is connected.
\end{proposition}
\begin{proof}
The result follows from (\ref{even1}) and from Corollary \ref{Zero Higgs},
just like in the proof of Proposition \ref{ccq}.
\end{proof}

Now we compute the components of $\M_1$.

\begin{theorem}\label{ccM1}
$\M_1$ has $2^{2g}$ components.
\end{theorem}
\begin{proof}
$\M_1(\overline\mu_1)=\bigcup_{L\in\Jac(X)}\M_{L,1}(\overline\mu_1)$ is connected since $\Jac^1(X)$ is connected and $\nu_1|_{\M_1(\overline\mu_1)}:\M_1(\overline\mu_1)\to\Jac^1(X)$ is a fibration with connected fibre $\M_{L,1}(\overline\mu_1)$, from Proposition \ref{ccq1}. The result follows since $\M_1=\bigsqcup_{\overline\mu_1\in(\Z_2)^{2g}}\M_1(\overline\mu_1)$.
\end{proof}

%%%%%%%%%%%%%%%%%%%%%%%%%%%%%%%%%%%%
\section{Topology of $\M_{\SL(3,\R)}$}
%%%%%%%%%%%%%%%%%%%%%%%%%%%%%%%%%%%%

In this subsection we shall consider the lower rank case of $\SL(3,\R)$-Higgs bundles. In holomorphic terms these are triples $(V,Q,\Phi)$ where $V$ is holomorphic vector bundle equipped with a nowhere degenerate quadratic form $Q$ and with trivial determinant, and $\Phi$ is a traceless $K$-twisted endomorphism of $V$, symmetric with respect to $Q$.

Let $\M_{\SL(3,\R)}$ be the moduli space of $\SL(3,\R)$-Higgs bundles. These objects are classified by the second Stiefel-Whitney class $w_2\in\{0,1\}$, and let $\M_{\SL(3,\R)}(w_2)$ be the subspace of $\M_{\SL(3,\R)}$ whose elements have the given $w_2$.

The moduli space $\M_{\SL(3,\R)}$ was considered in \cite{hitchin:1992} where the minimum subvarieties of the Hitchin functional were studied. There it was shown that if $(V,Q,\Phi)$ represents a fixed point of the circle action (\ref{circle action Higgs}), with $\Phi\neq 0$, then $V$ is of the form
$$V=F_{-m}\oplus\dots\oplus F_m$$ with $F_{-j}\cong F_j^*$, hence $\mathrm{rk}(F_j)=\mathrm{rk}(F_{-j})$ for all $j$. From this and since $\mathrm{rk}(V)=3$, we conclude that
fixed points with non-zero Higgs field are precisely
those such that $$V=F_{-1}\oplus\mathcal O\oplus F_1$$ with
$\mathrm{rk}(F_j)=1$ and, if $j\neq 1$, $\Phi_j:F_j\to F_{j+1}\otimes K$ is an isomorphism. These are local minima of the Hitchin function $f$.
The corresponding connected component, the Hitchin component, being isomorphic to a vector
space, is contractible.

For each $w_2\in\{0,1\}$, let $$\M'_{\SL(3,\R)}(w_2)$$ be the subspace of $\M_{\SL(3,\R)}(w_2)$ such that the minima on each of its connected components have $\Phi=0$.
Given $(V,Q,\Phi)\in\M_{\SL(3,\R)}(w_2)$, we know from
\cite{simpson:1994b}, that $$\lim_{t\to
0}(V,Q,t\Phi)$$ exists on $\M_{\SL(3,\R)}(w_2)$ and it is a fixed point of
the $\C^*$-action
$(V,Q,\Phi)\mapsto(V,Q,t\Phi)$ on $\M_{\SL(3,\R)}(w_2)$, being therefore a minimum of $f$. Hence, if
$(V,Q,\Phi)\in\M'_{\SL(3,\R)}(w_2)$, it follows that $\lim_{t\to
0}(V,Q,t\Phi)=(V',Q',0)\in\N_{\mathrm{SO}(3,\C)}(w_2)$ which is the space of local minima with zero Higgs field. Note that, in principle, it may happen that $\lim_{t\to 0}(V,Q,t\Phi)\neq (V,Q,0)=(V,Q)$, as $(V,Q)$ may be unstable as an ordinary orthogonal vector bundle.

Let us then consider the map
$$F:\M'_{\SL(3,\R)}(w_2)\times[0,1]\longrightarrow\M'_{\SL(3,\R)}(w_2)$$ given by
\begin{equation}\label{F}
F((V,Q,\Phi),t)=\begin{cases}
(V,Q,t\Phi) & \text{if }\ t\neq 0\\
\lim_{t\to 0}(V,Q,t\Phi) & \text{if }\ t=0.
\end{cases} 
\end{equation}

This map, together with the previous discussion, provides the
following result (recall that, from \cite{hitchin:1992}, we know that $\M_{\SL(3,\R)}$ has $3$ components).

\begin{theorem}\label{top3}
The space $\M_{\SL(3,\R)}$ has
one contractible component and the space consisting of the other two
components is homotopically equivalent to $\N_{\mathrm{SO}(3,\C)}$.
\end{theorem}
\begin{proof} The first part has already been discussed.
For the second part, we have to see that the map $F$ defined in (\ref{F}) is continuous, providing then a
retraction from
$\M'_{\SL(3,\R)}(w_2)$ into $\N_{\mathrm{SO}(3,\C)}(w_2)$, for each value of $w_2$. When $t\neq
0$, the continuity of $F$ is obvious. We will take care of the case
$t=0$.

The space $\M_{\GL(3,\C)}$ is a
quasi-projective, algebraic variety and $\C^*$ acts algebraically on
it as $(V,\Phi)\mapsto(V,t\Phi)$. Linearise this action with respect to an ample
line bundle $N$ (such that $N^s$ is very ample) over $\M_{\GL(3,\C)}$. This $\C^*$-action induces one on $N^s$ and, therefore we obtain a $\C^*$-action on $H^0(\M_{\GL(3,\C)},N^s)$ given by $$(t\cdot s)(V,\Phi)=t\cdot (s(V,t^{-1}\Phi)).$$ One can choose a rank $r+1$, $\C^*$-invariant subspace $W\subseteq H^0(\M_{\GL(3,\C)},N^s)$ and hence $\C^*$ acts on $W$. From this action we obtain a $\C^*$-action on $\mathbb{P}^r\cong\mathbb{P}(W)$, and there is a
$\C^*$-equivariant locally closed embedding
\begin{equation}\label{eqembb}
\iota:\M_{\GL(3,\C)}\hookrightarrow\mathbb{P}^r.
\end{equation}
If we linearise the given $\C^*$-action on $\mathbb{P}^r$ with respect to the very ample $\mathcal O_{\mathbb{P}^r}(1)$, then this is compatible with the morphism (\ref{eqembb}) and with the isomorphism $N^s\cong \mathcal O_{\mathbb{P}^r}(1)$.

Now, we can decompose $W$ as
$$W=\bigoplus_{i=1}^kW_{r_i}$$ where $r_i=\rk(W_{r_i})$ and $\C^*$ acts over each $W_{r_i}$
as $v\mapsto t^{\alpha_i}v$, $t\in\C^*$, $\alpha_i\in\Z$ and
$\alpha_i<\alpha_j$ whenever $i<j$. So, for each $r_i$, we have a
subspace of $\mathbb{P}^r$ given by
$\mathbb{P}^{r_i-1}=\mathbb{P}(W_{r_i})$.
With respect to the above decomposition of $W$, $\C^*$ acts as
\begin{equation}\label{action}
(v_1,\ldots,v_k)\mapsto (t^{\alpha_1}v_1,\ldots,t^{\alpha_k}v_k).
\end{equation}

Then, we also have the induced $\C^*$-action on the closed subspace $\M'_{\SL(3,\R)}(w_2)$ and a $\C^*$-equivariant
topological embedding 
\begin{equation}\label{topembb}
\iota|_{\M'_{\SL(3,\R)}(w_2)}:\M'_{\SL(3,\R)}(w_2)\hookrightarrow\mathbb{P}^r 
\end{equation}
and we denote the image in $\mathbb{P}^r$ of $\M'_{\SL(3,\R)}(w_2)$ through $\iota|_{\M'_{\SL(3,\R)}(w_2)}$ also by $\M'_{\SL(3,\R)}(w_2)$.
So we view $\M'_{\SL(3,\R)}(w_2)$ not as a
subvariety of $\mathbb{P}^r$, but simply a closed subspace (for the complex topology).

From (\ref{action}), the fixed point set of the $\C^*$-action on $\mathbb{P}^r$ is
$$\mathrm{Fix}_{\C^*}(\mathbb{P}^r)=\bigcup_{i=1}^k\mathbb{P}^{r_i-1}$$ so the
fixed point set of this action on $\M'_{\SL(3,\R)}(w_2)$ is
$$\mathrm{Fix}_{\C^*}(\M'_{\SL(3,\R)}(w_2))=\M'_{\SL(3,\R)}(w_2)\cap\bigcup_{i=1}^k\mathbb{P}^{r_i-1}.$$
But we already know that
$\mathrm{Fix}_{\C^*}(\M'_{\SL(3,\R)}(w_2))=\N_{\mathrm{SO}(3,\C)}(w_2)$ which is an
irreducible variety, by Theorem 5.9 of \cite{ramanathan:1996b}. So we conclude
that
\begin{equation}\label{fix}
\mathrm{Fix}_{\C^*}(\M'_{\SL(3,\R)}(w_2))=\M'_{\SL(3,\R)}(w_2)\cap\mathbb{P}^{r_{i_0}-1}
\end{equation}
for
some $i_0\in\{1,\dots,k\}$.

Actually, 
\begin{equation}\label{i0}
i_0=\min\{i\in\{1,\dots,k\}\;|\;v_i\neq 0,\,\text{ for some }(v_1,\dots,v_k)\in\M'_{\SL(3,\R)}(w_2)\}. 
\end{equation}
In fact, and let
$j=\min\{i\in\{1,\dots,k\}\;|\;v_i\neq 0,\,\text{ for some }(v_1,\dots,v_k)\in\M'_{\SL(3,\R)}(w_2)\}$ and let $(v_1,\dots,v_k)\in\M'_{\SL(3,\R)}(w_2)$ so that we can write
it as $(0,\dots,0,v_j,\dots,v_k)$. We have
\begin{equation*}
 \begin{split}
  \lim_{t\to 0}t(0,\dots,0,v_j,\dots,v_k)
&=\lim_{t\to 0}(0,\dots,0,t^{\alpha_j}v_j,\ldots,t^{\alpha_k}v_k)\\
&=\lim_{t\to
0}(0,\dots,0,v_j,t^{\alpha_{j+1}-\alpha_j}v_{j+1},\ldots,t^{\alpha_k-\alpha_j}v_k)\\
&=(0,\dots,v_j,\dots,0)\in\M'_{\SL(3,\R)}(w_2)\cap\mathbb{P}^{r_j-1}.
\end{split}
\end{equation*}
But, since we already know that $\lim_{t\to
0}t(v_1,\dots,v_k)\in\mathrm{Fix}_{\C^*}(\M'_{\SL(3,\R)}(w_2))$, we have
from (\ref{fix}) that $i_0=j$ and this settles (\ref{i0}).

If we take the map
$\tilde{F}:\mathbb{P}^r\times[0,1]\to\mathbb{P}^r$ given
by
$$\tilde{F}((v_1,\dots,v_k),t)=\begin{cases}
t(v_1,\dots,v_k)=(t^{\alpha_1}v_1,\ldots,t^{\alpha_k}v_k)\ &\text{if}\ t\neq 0\\
\lim_{t\to 0}t(v_1,\dots,v_k)=(0,\dots,0,v_{i_0},0,\dots,0)\
&\text{if}\ t=0 \end{cases}$$ then it is well-defined by the definition of $i_0$ in (\ref{i0}) and it is clearly continuous because $i_0$ is constant. By the compatibility of the actions, we have that
$F$ corresponds, under (\ref{topembb}), to $\tilde{F}|_{\M'_{\SL(3,\R)}(w_2)\times[0,1]}$, so $F$ is also continuous.\end{proof}

%%%%%%%%%%%%%%%%%%%%%%%%%%%%%%%%%%%%%%%%%%%%%%%%%%%%%%%%%%%%%
\section{Connected components of spaces of representations}
%%%%%%%%%%%%%%%%%%%%%%%%%%%%%%%%%%%%%%%%%%%%%%%%%%%%%%%%%%%%%

Recall that our main goal is to compute the number of components of
$\cR_{\PGL(n,\R)}$ for $n\geq 4$ even, but we had to work with the group
$\EGL(n,\R)$. The work done also gives a way
to count the components of the subspace of
$\cR=\cR_{\Gamma,\EGL(n,\R)}$
given by the disjoint union $\cR_0\sqcup\cR_1$. Denote
this subspace by $\cR_{0,1}$.

\begin{proposition}\label{ccR}
Let $n\geq 4$ be even. Then, $\cR_{0,1}$ has $2^{2g+1}+2$ connected components. More precisely,
\begin{enumerate}
        \item $\cR_0(\overline{\mu}_1)$ is connected, if $\overline{\mu}_1\neq 0$;
	\item $\cR_0(0,w_2)$ is connected, if $w_2\neq z_0$;
        \item $\cR_0(0,z_0)$ has $2$ components;
        \item $\cR_1(\overline{\mu}_1)$ is connected.
\end{enumerate}
\end{proposition}
\begin{proof}
By Theorem \ref{fundamental correspondence for semisimple G}, $\cR_0(\overline{\mu}_1)\cong\M_0(\overline{\mu}_1)$, $\cR_0(0,w_2)\cong\M_0(0,w_2)$ and $\cR_1(\overline{\mu}_1)\cong\M_1(\overline{\mu}_1)$.
The result follows directly from the analysis of the components of $\M_0$
and $\M_1$ in Theorems \ref{ccM0} and \ref{ccM1}.
\end{proof}

Now our main result follows as a corollary.

\begin{theorem}\label{princ}
Let $n\geq 4$ be even, and $X$ a closed oriented surface of genus $g\geq 2$. Then the moduli space $\cR_{\PGL(n,\R)}$ of reductive representations of $\pi_1X$ in $\PGL(n,\R)$ has
$2^{2g+1}+2$ connected components. More precisely,
\begin{enumerate}
        \item $\cR_{\PGL(n,\R)}(\mu_1,0)$ is connected, if $\mu_1\neq 0$;
        \item $\cR_{\PGL(n,\R)}(0,w_2)$ is connected, if $w_2\neq z_0$;
	\item $\cR_{\PGL(n,\R)}(0,z_0)$ has $2$ components;
        \item $\cR_{\PGL(n,\R)}(\mu_1,\omega_n)$ is connected.
\end{enumerate}
\end{theorem}
\begin{proof} The result follows immediately from the existence of the surjective continuous map
$p:\cR\to\cR_{\PGL(n,\R)}$ satisfying the identities of Proposition
\ref{projrep} and from the previous proposition. Note that the two components of $\cR_0(0,z_0)$ are not mapped into only
one in $\cR_{\PGL(n,\R)}(0,z_0)$ because if that were the case, every representation in $\PGL(n,\R)$ with $(0,z_0)$ as invariants could deform to a representation into $\PO(n)$, the maximal compact, and then the same would occur for the group
$\EGL(n,\R)$. We know however that this is not possible because of the analysis of the minima with invariants $(0,z_0)$: the component with minima with $\Phi\neq 0$ corresponds precisely to those representations which do not deform to a representation in $\EO(n)$. On the other hand, $\PGL(n,\R)$ is a split real form so by \cite{hitchin:1992} the space $\cR_{\PGL(n,\R)}$ should have a Hitchin component which in this case corresponds to the representations which do not deform to $\PO(n)$.
\end{proof}

\begin{remark}
For the proof of Theorem \ref{princ} is not essential to have Proposition \ref{ccR}. We could have used
Propositions \ref{ccq} and \ref{ccq1} and noticed that the vector bundles corresponding to minima of $f$ of type (2)
are projectively equivalent.
This would give us the number of components of $\M_{\PGL(n,\R)}(\mu_1,\mu_2)$, therefore of $\cR_{\PGL(n,\R)}(\mu_1,\mu_2)$ from
Theorem \ref{fundamental correspondence for semisimple G}.
\end{remark}

\begin{remark}
If $\mu_1=0$, then we might expect to get the same components as
Hitchin did in \cite{hitchin:1992} but that does not happen. We computed $4$ components while Hitchin's
result was $6$. The difference is that we are considering
$\PGL(n,\R)$-equivalence (cf. Remark \ref{PGL+}), while Hitchin considered
$\PSL(n,\R)$-equivalence.
\end{remark}

%%%%%%%%%%%%%%%%%%%%%%%%%%%%%%%%%%%%%%%%%%%%
\section{Topology of $\cR_{\SL(3,\R)}$}
%%%%%%%%%%%%%%%%%%%%%%%%%%%%%%%%%%%%%%%%%%%%

We finish with a corollary of Theorem \ref{top3}.
When $n$ is odd, $\PGL(n,\R)\cong\SL(n,\R)$, so $\cR_{\PGL(3,\R)}=\cR_{\SL(3,\R)}$.
Furthermore, from \cite{hitchin:1992} we know that $\cR_{\SL(3,\R)}$ has three components.

\begin{theorem}\label{princ2}
Let $X$ be a closed oriented surface of genus $g\geq 2$. The moduli space $\cR_{\SL(3,\R)}$ of reductive representations of $\pi_1X$ in $\SL(3,\R)$ has
one contractible component (the Hitchin component) and the space consisting of the other two
components is homotopically equivalent to $\cR_{\SO(3)}$.
\end{theorem}
\begin{proof}
The moduli space $\cR_{\SL(3,\R)}$ is isomorphic, via Theorem \ref{fundamental correspondence for semisimple G},
to $\M_{\SL(3,\R)}$. The result follows from Theorem \ref{top3}.
\end{proof}

Very recently, in \cite{ho-liu:2007}, Ho and Liu have computed, among other things, the Poincaré polynomials
of the spaces $\cR_{\SO(2n+1)}(w_2)$, $w_2=0,1$. For $n=3$, their result is (Theorem $5.5$ and Example $5.7$ of \cite{ho-liu:2007})
\begin{equation}\label{Ptw20}
P_t(\cR_{\SO(3)}(0))=\frac{-(1+t)^{2g}t^{2g+2}+(1+t^3)^{2g}}{(1-t^2)(1-t^4)} 
\end{equation}
and
\begin{equation}\label{Ptw21}
P_t(\cR_{\SO(3)}(1))=\frac{-(1+t)^{2g}t^{2g}+(1+t^3)^{2g}}{(1-t^2)(1-t^4)}. 
\end{equation}

From this result and from Theorem \ref{princ2}, we have:
\begin{theorem}
 The Poincaré polynomials of $\cR_{\SL(3,\R)}(w_2)$, $w_2=0,1$, are given by
$$P_t(\cR_{\SL(3,\R)}(0))=\frac{-(1+t)^{2g}t^{2g+2}+(1+t^3)^{2g}}{(1-t^2)(1-t^4)}+1$$
and $$P_t(\cR_{\SL(3,\R)}(1))=\frac{-(1+t)^{2g}t^{2g}+(1+t^3)^{2g}}{(1-t^2)(1-t^4)}.$$
\end{theorem}

\section*{Acknowledgments}
This paper is part of my PhD thesis \cite{oliveira:2008} and I would like to thank my supervisor, Peter
Gothen, for introducing me to the subject and for so many patient
explanations and fruitful discussions.
I thank Gustavo Granja for very helpful discussions about the topological classification of real projective bundles on surfaces.
Finally, I also thank an anonymous referee for pointing out several aspects which could be improved and for providing a much simpler argument for the topological classification of $G$-principal bundles on surfaces.

This work was partially supported by CMUP - Centro de Matemática da
Universidade do Porto, CMUTAD - Centro de Matemática da Universidade de Trás-os-Montes e Alto Douro 
and by the grant SFRH/BD/23334/2005 and the
project POCTI/MAT/58549/2004, financed by Fundação para a Ciência e
a Tecnologia (Portugal) through the programmes POCTI and POSI of the
QCA III (20002006) with European Community (FEDER) and national
funds.

\clearpage
\thispagestyle{empty}
\section{Erratum} The purpose of this note is to point out a flawed argument in \cite{oliveira:2011} which invalidates Theorem 9.1. This also affects the results in Section 11, namely Theorems 11.1 (stated as Theorem 1.3 in the Introduction) and 11.2, but the main results of the paper --- Theorems 1.1 and 1.2, concerning the number of connected components of the moduli space $\mathcal{M}_{\mathrm{PGL}(n,\R)}$ of $\mathrm{PGL}(n,\R)$-Higgs bundles over a compact Riemann surface $X$ --- are completely unaffected by this flaw. Indeed, Sections 9 and 11 only deal with the particular case of $\SL(3,\R)\cong\mathrm{PGL}(3,\R)$, where we consider the topology (more precisely, the Poincaré polynomial) of $\mathcal{M}_{\SL(3,\R)}$.

By the non-abelian Hodge correspondence, $\mathcal{M}_{\SL(3,\R)}$ is homeomorphic to the $\SL(3,\R)$-character variety $\mathcal{M}_{\SL(3,\R)}$ of representations $\pi_1X\to \SL(3,\R)$.

 $\SL(3,\R)$-Higgs bundles $(V,Q,\Phi)$ are topologically classified by the second Stiefel-Whitney class $w_2\in H^2(X,\Z_2)\cong\Z_2$ of the underlying rank $3$ orthogonal bundle $(V,Q)$. Denote by $\mathcal{M}_{\SL(3,\R)}(w_2)$ the subspace corresponding to $\SL(3,\R)$-Higgs bundles with $w_2(V,Q)=w_2$. Then 
 \[\mathcal{M}_{\SL(3,\R)}=\mathcal{M}_{\SL(3,\R)}(0)\sqcup\mathcal{M}_{\SL(3,\R)}(1).\]
 Moreover, $\mathcal{M}_{\SL(3,\R)}(1)$ is connected, while \[\mathcal{M}_{\SL(3,\R)}(0)=\mathcal{M}'_{\SL(3,\R)}(0)\sqcup \mathcal{M}''_{\SL(3,\R)}(0)\] has two connected components, one of them --- denoted here by $\mathcal{M}''_{\SL(3,\R)}(0)$ --- being the celebrated Hitchin component. Call the other two components, namely $\mathcal{M}'_{\SL(3,\R)}(0)$ and $\mathcal{M}_{\SL(3,\R)}(1)$, the \emph{non-Hitchin components}.
 
 In Theorem 9.1 we claim that each of the two non-Hitchin components admits a deformation retraction onto the corresponding components of the moduli space of rank $3$ orthogonal vector bundles (i.e.\ $\mathrm{SO}(3)$-Higgs bundles), which are also distinguished by the second Stiefel-Whitney class. However, \emph{this claim is false}. The reason is that, contrary to what is claimed on page 273, there are $\C^*$-fixed points in $\mathcal{M}'_{\SL(3,\R)}(0)$ and in $\mathcal{M}_{\SL(3,\R)}(1)$ which are not local minima of the Hitchin functional. 
 
 Indeed, a fixed point $(V,Q,\Phi)$, with non-zero Higgs field, must be of the form 
 \begin{equation}\label{eq:orth-bundle}
 V=F\oplus \mathcal{O}\oplus F^{-1}
 \end{equation} where $F$ is a line bundle, $Q$ is the obvious anti-diagonal quadratic form, and the Higgs field must be of one of the following forms:
 \begin{equation}\label{eq:Higgsfield1}
 \Phi=\left(\begin{array}{ccc}0 & 0 & 0 \\\Phi_1 & 0 & 0 \\0 & \Phi_1 & 0\end{array}\right)
 \end{equation} with $\Phi_1:F\to K$ non-zero, or 
 \begin{equation}\label{eq:Higgsfield2}
 \Phi=\left(\begin{array}{ccc}0 & 0 & 0 \\0 & 0 & 0 \\\Phi_1 & 0 & 0\end{array}\right)
 \end{equation} with $\Phi_1:F\to F^{-1}K$ non-zero.
In the case of \eqref{eq:orth-bundle} and \eqref{eq:Higgsfield1}, we must have 
\[0\leq \deg(F)\leq 2g-2,\]
while in the case of \eqref{eq:orth-bundle} and \eqref{eq:Higgsfield2}, 
\[0\leq \deg(F)\leq g-1.\]

Among the above fixed points, only the ones of type \eqref{eq:orth-bundle} and \eqref{eq:Higgsfield1} with $\deg(F)=2g-2$ (thus $F\cong K$) are local minima. These lie in the Hitchin component, which does not have any further fixed points (thus it is contractible). Hence the other fixed points lie in $\mathcal{M}_{\SL(3,\R)}'(0)$ or $\mathcal{M}_{\SL(3,\R)}(1)$. Indeed, it is easy to check that these fixed points have $w_2(V,Q)=\deg(F)\ \mathrm{mod}\, 2$, hence they exist both in $\mathcal{M}_{\SL(3,\R)}'(0)$ and $\mathcal{M}_{\SL(3,\R)}(1)$, depending on the degree of $F$.

This implies that the moduli space of rank $3$ special orthogonal bundles (which has two components distinguished by $w_2$), \emph{is not a deformation retraction of the space $\mathcal{M}'_{\SL(3,\R)}(0)\sqcup \mathcal{M}_{\SL(3,\R)}(1)$ given by the non-Hitchin components}. Hence Theorems 9.1, 11.1 (and 1.3) and 11.2. do not hold.

\end{document}